\documentclass[12pt,a4paper]{article}
\usepackage{amssymb}

\usepackage{amsmath}

\newcounter{resultnum}[section]
\newtheorem{conclusion}{Conclusion}[section]

\newcounter{conclusionnum}[section]

\newcounter{conditionnum}[section]

\newcounter{conjecturenum}[section]
\newtheorem{example}{Example}[section]

\newcounter{examplenum}[section]

\newcounter{exercisenum}[section]
\newtheorem{lemma}{Lemma}[section]

\newcounter{lemmanum}[section]

\newcounter{notationnum}[section]
\newtheorem{theorem}{Theorem}[section]

\newcounter{theoremnum}[section]
\newtheorem{definition}{Definition}[section]

\newcounter{definitionnum}[section]
\newtheorem{corollary}{Corollary}[section]

\newcounter{corollarynum}[section]
\newtheorem{remark}{Remark}[section]

\newcounter{remarknum}[section]

\newcounter{propositionnum}[section]

\newcounter{acknowledgementnum}[section]

\newcounter{algorithmnum}[section]

\newcounter{axiomnum}[section]

\newcounter{casenum}[section]

\newcounter{claimnum}[section]

\newcounter{summarynum}[section]

\newcounter{problemnum}[section]
\newenvironment{proof}[1][]{\textbf{Proof.} }{}

\begin{document}

\title{Nonholonomic Ricci Flows:\\
III. Curve Flows and Solitonic Hierarchies}
\date{April 18, 2007}
\author{ Sergiu I. Vacaru\thanks{%
sergiu$_{-}$vacaru@yahoo.com, svacaru@fields.utoronto.ca } \\
%EndAName
{\quad} \\
\textsl{The Fields Institute for Research in Mathematical Science} \\
\textsl{222 College Street, 2d Floor, } \textsl{Toronto \ M5T 3J1, Canada} }
\maketitle

\begin{abstract}
The geometric constructions are elaborated on (semi) Rieman\-ni\-an manifolds
and vector bundles provided with nonintegrable distributions defining
nonlinear connection structures induced canonically by metric tensors. Such
spaces are called nonholonomic manifolds and described by two equivalent
linear connections also induced in unique forms by a metric tensor (the Levi
Civita and the canonical distinguished connection, d--connec\-ti\-on). The
lifts of geometric objects on tangent bundles are performed for certain
classes of d--connections and frame transforms when the Riemann tensor is
parametrized by constant matrix coefficients. For such configurations, the
flows of non--stretching curves and corresponding bi--Hamilton and solitonic
hierarchies encode information about Ricci flow evolution, Einstein spaces
and exact solutions in gravity and geometric mechanics. The applied methods
were elaborated formally in Finsler geometry and allows us to develop the
formalism for generalized Riemann--Finsler and Lagrange spaces.
Nevertheless, all geometric constructions can be equivalently re--defined
for the Levi Civita connections and holonomic frames on (semi) Riemannian
manifolds.


\end{abstract}

%\tableofcontents

\section{ Introduction}

Both the theory of Ricci flows and the theory of integrable partial
differential equations have deep links to the geometry of curves and
surfaces, generalized Riemann--Finsler spaces and geometric analysis:

Originally, the Ricci flow theory has addressed geometrical and topological
issues, and certain applications in physics, for Riemannian manifolds \cite%
{3ham1,3ham2,3per1,3per2,3per3} (we cite here some reviews on
Hamilton--Perelman theory \cite{3caozhu,3cao,3kleiner,3rbook}). In parallel,
it was found that various classes of solitonic equations (for instance, the
sine--Gordon, SG, and modified Korteveg--de Vries, mKdV, equations) and
along with their hierarchies of symmetries, conservation laws and associated
recursion operators can be encoded into the geometry of flows of
non--stretching curves in Riemannian symmetric spaces and related Lie
algebras and Klein spaces \cite{3ath,3gp,3nsw,3sw,3chou1,3mbsw,3saw,3serg},
see also reviews and new results in Refs. \cite{3anc2,3anc3}.

In \cite{3vsh}, it was proven that a more general class of (pseudo)
Riemannian spaces can be encoded into bi--Hamilton structures and related
solitonic hierarchies. The key construction was to deform nonholonomically
the frame and linear connection structures in order to get constant matrix
curvature coefficients, with respect to certain classes of nonholonomic
frames. Such frames are adapted to a nonlinear connection (N--connection)
structure induced by some generic off--diagonal metric coefficients . We
also concluded that having generated the corresponding solitonic hierarchies
for the so--called canonical distinguished connection (d--connection), we
can re--define equivalently the geometric objects, conservation laws and
basic equations and solutions in terms of the usual Levi Civita connection.

The formalism applied in \cite{3vsh}, based on the geometry of moving
nonholonomic frames with associated N--connection structure, was originally
developed in Finsler--Lagrange geometry and generalizations \cite%
{3ma1,3ma2,3bej}. Our idea \cite{3vsgg} was to apply it to usual (semi)
Riemannian spaces, or to Riemann--Cartan ones (with nontrivial torsion),
prescribing certain nonholonomic distributions arising naturally if we
constrain partially some degrees of freedom. For such systems, there are
defined certain classes of preferred frames and symmetries for the
gravitational and matter field interactions. Following this approach, it was
possible to construct various classes of exact solutions in Einstein and
string gravity modelling Finsler like locally anisotropic structure,
possessing noncommutative symmetries and defining generically off--diagonal
metrics and nonlinear interactions of pp--waves, two and three dimensional
gravitational solitons and spinor fields \cite{3vsgg,3vncg,3vcfa}.

Nevertheless, if realistic theories of gravitational and gauge field
interactions and/or generalized Finsler models of geometric mechanics are
introduced into consideration, the solitonic encoding of metric, connection
and frame structures is possible for certain effective generalized Lagrange
spaces. In this case, we model the geometric constructions on couples of
symmetric Riemannian spaces provided with nonholonomic distributions. The
work \cite{3avw} concluded an approach when different theories of gravity
and geometric mechanics are treated in a unified geometric way as
nonholonomic (semi) Riemannian manifolds, or vector bundles. It was also
proven that the data for geometric objects and fundamental physical
equations (their symmetries and conservation laws) can be encoded into
bi--Hamilton structures and correspondingly derived solitonic hierarchies.

Integrable and nonintegrable (i.e. holonomic and nonholonomic / anholonomic)
flows of geometric structures were also recently considered in a series of
works on nonholonomic Ricci flows \cite{3nrf1,3nrf2} and applications in
physics \cite{3nrfes,3nrft1,3nrft2,3nfelf}. Some important results of those
works were the proofs that constrained Ricci flows of (semi) Riemannian
metrics result in Finsler like metrics and generalizations and, inversely,
Finsler--Lagrange type geometrical objects can be described equivalently by
flows on Riemann (in general, Riemann--Cartan) spaces.

The goal of this paper, the third one in a series defined together with \cite%
{3nrf1,3nrf2}, is to prove that solitonic hierarchies can be generated by
any (semi) Riemannian metric $g_{ij}$ on a manifold $V$ of dimension $\dim
V=n\geq 2$ if the the geometrical objects are lifted in the total space of
the tangent bundle $TV,$ or of a vector bundle $\mathcal{E}=(M,\pi ,E),$ $%
\dim E=m\geq n,$ by defining such frame transforms when constant matrix
curvatures are defined canonically with respect to certain classes of
preferred systems of reference. We shall also define the criteria when
families of bi--Hamilton structures and solitonic hierarchies encode (in
general, nonholonomic) Ricci flow evolutions of geometric objects and/ or
exact solutions of gravitational field equations.

The paper is organized as follows:

In section 2 we outline the geometry of nonholonomic manifolds and vector/
tangent bundles provided with nonlinear connection structure. We emphasize
the possibility to define fundamental geometric objects induced by a (semi)
Riemannian metric on the base space when the Riemannian curvature tensor has
constant coefficients with respect to a preferred nonholonomic basis. We
also present some results on evolution equations of nonholonomic Ricci flows
and exact solutions in gravity.

In section 3 we consider Ricci flow families of curve flows on nonholonomic
vector bundles. We sketch an approach to classification of such spaces
defined by conventional horizontal and vertical symmetric (semi) Riemannian
subspaces and provided with nonholonomic distributions defined by the
nonlinear connection structure. It is constructed a corresponding family of
nonholonomic Klein spaces for which the bi--Hamiltonian operators are
defined by canonical distinguished connections, adapted to the nonlinear
connection structure, for which the distinguished curvature coefficients are
constant.

Section 4 is devoted to the formalism of distinguished bi--Hamiltonian
operators and vector soliton equations for arbitrary (semi) Riemannian
spaces. Then we consider the properties of cosympletic and sympletic
operators adapted to the nonlinear connection structure. We define the basic
equations for nonholonomic curve flows and parametrize their possible Ricci
flows.

Section 5 is devoted to formulation of the Main Result: a proof that for any
nonholonomic Ricci flow system, one can be defined a natural family of
N--adapted bi--Hamiltonian flow hierarchies inducing anholonomic solitonic
configurations. There are constructed in explicit form the solitonic
hierarchies corresponding to the bi--Hamiltonian anholonomic curve flows.
Finally, there are speculated the conditions when from solitonic hierarchies
we can extract solutions of the Ricci flow and/or field equations.

We summarize and discuss the results in section 6. For convenience, we
outline the necessary definitions and formulas from the geometry of
nonholonomic manifolds in Appendix A. Then, in Appendix B, we consider the
geometry of N--anholonomic Klein spaces. A proof of the Main Theorem is
sketched in Appendix C.

\subsection*{Notation remarks:}

There are considered two types of flows of geometrical objects on manifolds
of necessary smooth class, induced by 1) non--stretching curve flows $\gamma
(\tau ,\mathbf{l}),$ defined by real parameters $\tau $ and $\mathbf{l,}$
and 2) Ricci flows of metrics $g_{ij}(\chi ),$ parametrized by a real $\chi
. $ The non--stretching flows of a curve are constrained by the condition $%
g_{ij}(\gamma _{\tau },\gamma _{\mathbf{l}})=1,$ which under Ricci flows
transforms into a family of such conditions, $g_{ij}(\gamma _{\tau },\gamma
_{\mathbf{l}},\chi )=1.$ For Ricci flows, we get evolutions of families of
non--stretching curves parametrized by hypersurfaces $\gamma (\tau ,\mathbf{%
l,}\chi )\doteqdot ~^{\chi }\gamma (\tau ,\mathbf{l}).$ It is convenient to
use in parallel two types of denotations for the geometric objects subjected
to both curve and Ricci flows: by emphasized all dependencies on parameters $%
\tau ,\mathbf{l}$ and $\chi $ or by introducing ''up/low'' labels like $%
~^{\chi }\gamma =\gamma (...,\chi )$ or $~_{\chi }A=A(...,\chi ).$

We shall also write ''boldface'' symbols for geometric objects and spaces
adapted to a noholonomic/ nonlinear connection structure, like $\mathbf{V,E}%
,...$ and write $V,E,...$ if the nonholonomic structure became trivial, i.e.
integrable/ holonomic. In order to investigate the properties of curve and
Ricci flow evolution equations it is convenient to use both abstract/global
denotations and coefficient formulas with respect to coordinate or
nonholonomic bases.

A nonholonomic distribution with associated nonlinear connection structure
splits the manifolds into conventional horizontal (h) and vertical (v)
subspaces. The geometric objects, for instance, a vector $\mathbf{X}$ can be
written in abstract form as $\mathbf{X}=(hX,vX)=(~^{h}X,~^{v}X),$ or in
coefficient forms as $\mathbf{X}^{\alpha }=(X^{i},X^{a})=(X^{\underline{i}%
},X^{\underline{a}}),$ where $\mathbf{X=X}^{\alpha }\mathbf{e}_{\alpha
}=X^{i}e_{i}+X^{a}e_{a}=X^{\underline{i}}\partial _{\underline{i}}+X^{%
\underline{a}}\partial _{\underline{a}}$ can be equivalently decomposed with
respect to a general nonholonomic frame $\mathbf{e}_{\alpha }=(e_{i},e_{a})$
or coordinate frame $\partial _{\underline{\alpha }}=(\partial _{\underline{i%
}},\partial _{\underline{a}})$ for local\ h- and v--coordinates $u=(x,y),$
or $u^{\alpha }=(x^{i},y^{a})$ when $\partial _{\underline{\alpha }%
}=\partial /\partial u^{\underline{\alpha }}$ and $\partial _{\underline{i}%
}=\partial /\partial x^{\underline{i}},\partial _{\underline{a}}=\partial
/\partial x^{\underline{a}},$ when indices will be underlined if it is
necessary to emphasize certain decompositions are defined for coordinate
bases. The h--indices $i,j,k,...=1,2,...n$ will be used for base/
nonholonomic vector objects and the v--indices $a,b,c...=n+1,n+2,...n+m$
will be used for fiber/ holonomic vector objects. Greek indices of type $%
\alpha ,\beta ,...$ will be used as cumulative ones.

Finally, we note that we shall omit labels, indices and parametric/
coordinate dependencies for some formulas if it does not result in
ambiguities.

\section{Nonholonomic Lifts and Ricci Flows}

In this section, we prove that for any family of (semi) Riemannian metrics $%
g_{ij}(\chi )$ on a manifold $V,$ parametrized by a real parameter $\chi ,$
it is possible to define lifts to the tangent bundle $TV$ provided with
canonical nonlinear connection (in brief, N--connection), Sasaki type
metrics and (linear) canonical distinguished connection (d--connection)
structures. We also outline some important formulas for nonholonomic Ricci
flow evolution equations of geometric structures. The reader is recommended
to consult Refs. \cite{3vsgg,3vncg,3ma1,3ma2,3nrf1,3nrf2} and Appendix A on
details on N--connection geometry and recent developments in modern gravity
and Ricci flow theory.

\subsection{N--connections induced by families of Riemannian metrics}

Let $\mathcal{E}=(E,\pi ,F,M)$ be a (smooth) vector bundle over base
manifold $M,$ $\dim M=n$ and $\dim E=(n+m),$ for $n\geq 2,$ and $m\geq n$
being the dimension of the typical fiber $F.$ It is defined a surjective
submersion $\pi :E\rightarrow M.$ In any point $u\in E,$ the total space $E$
splits into ''horizontal'', $M_{u},$ and ''vertical'', $F_{u},$ subspaces.
We denote the local coordinates in the form $u=(x,y),$ or $u^{\alpha
}=\left( x^{i},y^{a}\right) ,$ with horizontal indices $i,j,k,\ldots
=1,2,\ldots ,n$ and vertical indices $a,b,c,\ldots =n+1,n+2,\ldots ,n+m.$%
\footnote{%
In a particular case, we have a tangent bundle $E\mathbf{=}TM,$ when $n=m;$
for such bundles both type of indices run the same values but it is
convenient to distinguish the horizontal and vertical ones by using
different groups of small Latin indices. Here one should be noted that on $%
TM $ we are able to contract the vertical indices with the corresponding
horizontal ones, and inversely, but not on a general nonholonomic manifold ${%
\mathbf{V}},$ or ${\mathbf{E}}.$} The summation rule on the same ''up'' and
''low'' indices will be applied.

The base manifold $M$ is provided with a family of (semi) Riemannian
metrics, nondegenerate second rank tensors,\footnote{%
in physical literature, one uses the term (pseudo) Riemannian/Euclidean space%
} $h\underline{g}(\chi )=\underline{g}_{ij}(x,\chi )dx^{i}\otimes dx^{j},$
for $0\leq \chi \leq \chi _{0}\in \mathbb{R}.$ We also introduce a family of
vertical metrics $v\underline{g}(\chi )=\underline{g}_{ab}(x,\chi
)dy^{a}\otimes dy^{b}$ by completing the matrices $\underline{g}_{ij}(x,\chi
) $ diagonally with $\pm 1$ till any nondegenerate second rank tensor $%
\underline{g}_{ab}(x,\chi )$ if $m>n$ and then subjecting to any frame
transforms. This way, we define certain families of metrics $\underline{%
\mathbf{g}}(\chi )=[h\underline{g}(\chi ),v\underline{g}(\chi )]$ (we shall
also use the notation $\underline{g}_{\alpha \beta }(\chi )=[\underline{g}%
_{ij}(\chi ),\underline{g}_{ab}(\chi )])$ on $\mathcal{E}.$ Considering
frame (vielbein) transforms,%
\begin{equation}
g_{\alpha \beta }(x,y,\chi )=e_{\alpha }^{~\underline{\alpha }}(x,y,\chi
)~e_{\beta }^{~\underline{\beta }}(x,y,\chi )g_{\underline{\alpha }%
\underline{\beta }}(x,\chi ),  \label{3auxm}
\end{equation}%
where $\underline{g}_{\alpha \beta }(x,\chi )$ is written in equivalent form
$g_{\underline{\alpha }\underline{\beta }}(x,\chi ),$ we can deform the
metric structures, $\underline{g}_{\alpha \beta }\rightarrow g_{\alpha \beta
}=[g_{ij},g_{ab}]$ (we shall omit dependencies on coordinates and parameters
if it does not result in ambiguities). The coefficients $e_{\alpha }^{~%
\underline{\alpha }}(x,y,\chi )$ will be defined below (see formula (\ref%
{3aux4})) from the condition of generating curvature tensors with constant
coefficients with respect to certain preferred systems of reference.

For any $g_{ab}(\chi )$ from the set $g_{\alpha \beta }(\chi ),$ we can
construct a family of effective generation functions%
\begin{equation*}
\mathcal{L}(x,y,\chi )=g_{ab}(x,y,\chi )y^{a}y^{b}
\end{equation*}%
inducing families of vertical metrics
\begin{equation}
\tilde{g}_{ab}(x,y,\chi )=\frac{1}{2}\frac{\partial ^{2}\mathcal{L}(x,y,\chi
)}{\partial y^{a}\partial y^{b}}  \label{3ehes}
\end{equation}%
which is ''weakly'' regular if $\det |\tilde{g}_{ab}|\neq 0.$

By straightforward computations, we prove\footnote{%
see Refs. \cite{3ma1,3ma2} for details of a similar proof; here we note that
in our case, in general, \ $e_{\alpha }^{~\underline{\alpha }}\neq \delta
_{\alpha }^{~\underline{\alpha }}$}:

\begin{theorem}
\label{3teleq}The family of \ Lagrangians $L(\chi )=\sqrt{|\mathcal{L}(\chi
)|},$ where $y^{i}=\frac{dx^{i}}{d\tau }$ for paths $x^{i}(\tau )$ on $M,$
depending on parameter $\tau ,$ with weakly regular metrics induces a family
of Euler--Lagrange equations on $TM,$
\begin{equation*}
\frac{d}{d\tau }\left( \frac{\partial L(\chi )}{\partial y^{i}}\right) -%
\frac{\partial L(\chi )}{\partial x^{i}}=0,
\end{equation*}
which are equivalent to a corresponding family of ``nonlinear'' geodesic
equations
\begin{equation*}
\frac{d^{2}x^{i}}{d\tau ^{2}}+2\widetilde{G}^{i}(x^{k},\frac{dx^{j}}{d\tau }%
,\chi )=0
\end{equation*}%
defining paths of a canonical semispray $S(\chi )=y^{i}\frac{\partial }{%
\partial x^{i}}-2\widetilde{G}^{i}(x,y,\chi )\frac{\partial }{\partial y^{i}}%
,$ where
\begin{equation*}
2\widetilde{G}^{i}(x,y,\chi )=\frac{1}{2}\ \tilde{g}^{ij}(\chi )\left( \frac{%
\partial ^{2}L(\chi )}{\partial y^{i}\partial x^{k}}y^{k}-\frac{\partial
L(\chi )}{\partial x^{i}}\right)
\end{equation*}%
with $\tilde{g}^{ij}(\chi )$ being inverse to (\ref{3ehes}).
\end{theorem}

The Theorem \ref{3teleq} states the possibility to geometrize the regular
Lagrange mechanics by geometric objects on nonholonomic spaces and inversely:

\begin{conclusion}
For any family of (semi) Riemannian metrics $\underline{g}_{ij}(x,\chi )$ on
$M,$ we can associate canonically certain families of effective regular
Lagrange mechanical systems on $TM$ with the Euler--Lagrange equations
transformed into corresponding families of nonlinear (semispray) geodesic
equations.
\end{conclusion}

\begin{theorem}
Any family of (semi) Riemannian metrics $\underline{g}_{ij}(x,\chi )$ on $M$
induces a corresponding family of canonical N--connection structures on $TM.$
\end{theorem}

\begin{proof}
We sketch a proof by defining the coefficients of N--connection, see (\ref%
{3whitney}),
\begin{equation}
\tilde{N}_{\ j}^{i}(x,y,\chi )=\frac{\partial \tilde{G}^{i}(x,y,\chi )}{%
\partial y^{j}}  \label{3cnlce}
\end{equation}%
where
\begin{eqnarray}
\tilde{G}^{i}(\chi ) &=&\frac{1}{4}\tilde{g}^{ij}\left( \frac{\partial ^{2}%
\mathcal{L}}{\partial y^{i}\partial x^{k}}y^{k}-\frac{\partial \mathcal{L}}{%
\partial x^{j}}\right) =\frac{1}{4}\tilde{g}^{ij}(\chi )~g_{jk}(\chi
)~\gamma _{lm}^{k}(\chi )y^{l}y^{m},  \label{3aux2} \\
\gamma _{\ lm}^{i}(\chi ) &=&\frac{1}{2}g^{ih}(\chi )~\left[ \partial
_{m}g_{lh}(\chi )+\partial _{l}g_{mh}(\chi )-\partial _{h}g_{lm}(\chi )%
\right] ,\ \partial _{h}=\partial /\partial x^{h},  \notag
\end{eqnarray}%
with $g_{ah}(\chi )$ and $\tilde{g}_{ij}(\chi )$ defined respectively by
formulas (\ref{3auxm}) and (\ref{3ehes}). $\square $
\end{proof}

The families of N--adapted partial derivative and differential operators,
see Appendix for more general formulas (\ref{3dder}) and (\ref{3ddif}), are
defined by the N--connection coefficients (\ref{3cnlce}) and may be denoted
respectively $\mathbf{\tilde{e}}_{\nu }(\chi )=(\mathbf{\tilde{e}}_{i}(\chi
),e_{a})$ and $\mathbf{\tilde{e}}^{\mu }(\chi )=(e^{i},\mathbf{\tilde{e}}%
^{a}(\chi )).$

For any metric structure $\mathbf{g}$ on a manifold, there is the unique
metric compatible and torsionless Levi Civita connection $\nabla $ for which
$\ ^{\nabla }\mathcal{T}^{\alpha }=0$ and $\nabla \mathbf{g=0.}$ This
connection is not a d--connection because it does not preserve under
parallelism the N--connection splitting (\ref{3whitney}). One has to
consider less constrained cases, admitting nonzero torsion coefficients,
when a d--connection is constructed canonically for a d--metric structure. A
simple minimal metric compatible extension of $\nabla $ is that of canonical
d--connection $\widehat{\mathbf{D}},$ with $T_{\ jk}^{i}=0$ and $T_{\
bc}^{a}=0$ but, $\ $in general, nonzero $T_{\ ji}^{a}$ and $T_{\ bi}^{a},$
see (\ref{3dtors}). The coefficient formulas for such connections are given
in Appendix, see (\ref{3candcon}) and related discussion. It should be noted
that on tangent bundle $TM$ it is possible to define the torsionless
canonical d--connection (\ref{3candcontm}) which is completely similar to
the Levi Civita connection. For families of metrics $\mathbf{g}(\chi ),$ we
get certain families of connections $\nabla (\chi )$ and $\widehat{\mathbf{D}%
}(\chi ).$

\begin{theorem}
Any family of (semi) Riemannian metrics $\underline{g}_{ij}(x,\chi )$ on $M$
\ induces a parametrized by $\chi$ family of nonholonomic (semi) Riemannian
structures on $TM.$
\end{theorem}

\begin{proof}
The family $\underline{g}_{ij}(x,\chi )$ on $M$ induces a family of
canonical d--metric structures on $TM,$
\begin{equation}
\mathbf{\tilde{g}}(\chi )=\tilde{g}_{ij}(x,y,\chi )\ e^{i}\otimes e^{j}+\
\tilde{g}_{ij}(x,y,\chi )\ \mathbf{\tilde{e}}^{i}(\chi )\otimes \mathbf{%
\tilde{e}}^{j}(\chi ),  \label{3slme}
\end{equation}%
where $\mathbf{\tilde{e}}^{i}(\chi )$ are elongated as in (\ref{3ddif}), but
with $\tilde{N}_{\ j}^{i}(\chi )$ from (\ref{3cnlce}). Then, we note that\
there are canonical d--connections on $TM$ induced by $\underline{g}%
_{ij}(x,\chi ):$ we can construct them in explicit form by introducing $%
\tilde{g}_{ij}(\chi )$ and $\tilde{g}_{ab}(\chi )$ in formulas (\ref%
{3candcontm}), in order to compute the coefficients $\tilde{\Gamma}_{\ \beta
\gamma }^{\alpha }(\chi )=(\tilde{L}_{\ jk}^{i}(\chi ),\tilde{C}%
_{bc}^{a}(\chi )).\square $
\end{proof}

The corresponding curvature curvature tensor%
\begin{equation*}
\widetilde{R}_{\ \beta \gamma \tau }^{\alpha }(\chi )=\{\widetilde{R}_{\
hjk}^{i}(\chi ),\tilde{P}_{\ jka}^{i}(\chi ),\tilde{S}_{\ bcd}^{a}(\chi )\}
\end{equation*}%
can be computed by introducing respectively the values $\tilde{g}_{ij}(\chi
),\tilde{N}_{\ j}^{i}(\chi )$ and $\widetilde{\mathbf{e}}_{k}(\chi )$ into (%
\ref{3candcontm}), defining $\tilde{\Gamma}_{\ \beta \gamma }^{\alpha }(\chi
)=\left( \tilde{L}_{\ jk}^{i}(\chi ),\tilde{C}_{bc}^{a}(\chi )\right) $ and
then into formulas (\ref{3dcurvtb}). Here one should be noted that the
constructions on $TM$ depend on arbitrary vielbein coefficients $e_{\alpha
}^{~\underline{\alpha }}(x,y,\chi )$ in (\ref{3auxm}). We can restrict such
sets of coefficients in order to generate various particular classes of
(semi) Riemannian geometries on $TM,$ for instance, in order to generate
symmetric Riemannian spaces with constant curvature, see Refs. \cite%
{3helag,3kob,3sharpe}.

\begin{corollary}
There are lifts of a family of (semi) Riemannian metric $\underline{g}%
_{ij}(x,\chi )$ on $M,$ $\dim M=n,$ generating a corresponding family
Riemannian structures on $TM$ with the curvature coefficients of the
canonical d--connections coinciding (with respect to N--adapted bases) to
those for the families of Riemannian space of constant curvature of
dimension $n+n.$
\end{corollary}

\begin{proof}
For a given set $\underline{g}_{ij}(x,\chi )$ on $M,$ in (\ref{3auxm}), we
chose such coefficients $e_{\alpha }^{~\underline{\alpha }}(x,y,\chi
)=\left\{ e_{a}^{~\underline{a}}(x,y,\chi )\right\} $ that
\begin{equation*}
g_{ab}(x,y,\chi )=e_{a}^{~\underline{a}}(x,y,\chi )~e_{b}^{~\underline{b}%
}(x,y,\chi )g_{\underline{a}\underline{b}}(x,\chi )
\end{equation*}%
results in (\ref{3ehes}) of type%
\begin{equation}
\tilde{g}_{ef}(\chi )=\frac{1}{2}\frac{\partial ^{2}\mathcal{L}(\chi )}{%
\partial y^{e}\partial y^{f}}=\frac{1}{2}\frac{\partial ^{2}(e_{a}^{~%
\underline{a}}~e_{b}^{~\underline{b}}y^{a}y^{b})}{\partial y^{e}\partial
y^{f}}g_{\underline{a}\underline{b}}(x,\chi )=\ \mathring{g}_{ef}(\chi ),
\label{3aux4}
\end{equation}%
where $\ \mathring{g}_{ab}(\chi )$ are metrics of symmetric Riemannian space
(of constant curvature). Considering a prescribed set $~\ \mathring{g}%
_{ab}(\chi ),$ we have to integrate two times on $y^{e}$ in order to find
any solution for $e_{a}^{~\underline{a}}(\chi )$ defining a frame structure
in the vertical subspace. The next step is to construct the d--metric $~\
\mathring{g}_{\alpha \beta }(\chi )=[\ \mathring{g}_{ij}(\chi ),\ \mathring{g%
}_{ab}(\chi )]$ of type (\ref{3slme}), in our case, with respect to a
nonholonomic base elongated by $~\widetilde{\mathring{N}^{i}}_{~j}(\chi ),$
generated by $\underline{g}_{ij}(x,\chi )$ and $\tilde{g}_{ef}=\mathring{g}%
_{ab}(\chi ),$ like in (\ref{3cnlce}) and (\ref{3aux2}). This defines a
constant curvature Riemannian space of dimension $n+n.$ The coefficients of
the canonical d--connection, which in this case coincide with those for the
Levi Civita connection, and the coefficients of the Riemannian curvature can
be computed respectively by introducing $\tilde{g}_{ef}=\ \mathring{g}%
_{ab}(\chi )$ in formulas (\ref{3candcontm}) and (\ref{3dcurvtb}). Finally,
we note that the induced symmetric Riemannian spaces contain additional
geometric structures like the N--connection and anholonomy coefficients $%
W_{\alpha \beta }^{\gamma }(\chi ),$ see (\ref{3anhrel}).$\square $
\end{proof}

There are various possibilities to generate on $TM$ nonholonomic Riemannian
structures from a given set $\underline{g}_{ij}(x,\chi )$ on $M.$ They
result in different geometrical and physical models.

\begin{remark}
\label{rem11}We can simplify substantially the geometric constructions if
instead of families of constant coefficients $\ \mathring{g}_{ef}(\chi ),$
we consider only one set of constant coefficients $\ \mathring{g}_{ef}=\
\mathring{g}_{ef}(\chi _{0}).$ This is possible even $g_{\underline{a}%
\underline{b}}(x,\chi )$ have quite general dependencies on $\chi $ but
supposing that we can define such $e_{a}^{~\underline{a}}(x,y,\chi )$ when (%
\ref{3aux4}) can be solved for a fixed right side. For simplicity, in our
furhter considerations we shall fix any set $\mathring{g}_{ef}$ and
parametrize the dependencies on $\chi $ for $\underline{g}_{ij}$ and $%
e_{a}^{~\underline{a}}.$
\end{remark}

In this work, we emphasize the possibility of generating spaces with
constant curvature because for such symmetric spaces it was elaborated a
bi--Hamiltonian approach and corresponding solitonic hierarchies \cite%
{3mbsw,3saw,3anc2,3anc3}. For general Riemannian and/or Finsler--Lagrange
spaces it is not possible to get constant coefficient curvature coefficients
for the Levi Civita connection but for the corresponding lifts to the
canonical d--connection there are constructions generating constant
curvature coefficients \cite{3vsh,3avw}. This will allow us to construct the
corresponding solitonic hierarchies from which, by imposing the
corresponding constraints, it will be possible to extract the information
for very general classes of metrics.

\begin{example}
The simplest example when a Riemannian structure with constant matrix
curvature coefficients is generated on $TM$ is given by a d--metric induced
by $\tilde{g}_{ij}=\delta _{ij},$ i.e.%
\begin{equation}
\mathbf{\tilde{g}}_{[E]}=\delta _{ij}e^{i}\otimes e^{j}+\ \delta _{ij}\
\mathbf{\tilde{e}}^{i}\otimes \mathbf{\tilde{e}}^{j},  \label{3clgs}
\end{equation}%
with $\mathbf{\tilde{e}}^{i}$ defined by $\tilde{N}_{\ j}^{i}$ in their turn
induced by a given set $\underline{g}_{ij}(x)$ on $M.$ For families of
geometric objects, we consider
\begin{equation*}
\mathbf{\tilde{g}}_{[E]}(\chi )=\delta _{ij}e^{i}\otimes e^{j}+\ \delta
_{ij}\ \mathbf{\tilde{e}}^{i}(\chi )\otimes \mathbf{\tilde{e}}^{j}(\chi ),
\end{equation*}%
when $\tilde{N}_{\ j}^{i}(\chi )$ are defined by a given set $\underline{g}%
_{ij}(x,\chi ).$
\end{example}

For more general nonholonomic configurations on $TM,$ we can consider
families of metrics of type
\begin{equation}
\mathbf{g}(\chi )=g_{ij}(x,y,\chi )\ e^{i}\otimes e^{j}+\ g_{ab}(x,y,\chi )\
\mathbf{e}^{a}(\chi )\otimes \mathbf{e}^{b}(\chi ),  \label{3lme}
\end{equation}%
where
\begin{equation*}
g_{ij}(x,y,\chi )=\eta _{ij}(x,y,\chi )\tilde{g}_{ij}(x,y,\chi )\mbox{ and }%
g_{ab}(x,y,\chi )=\eta _{ab}(x,y,\chi )\tilde{g}_{ab}(x,y,\chi )
\end{equation*}%
and $\mathbf{e}^{a}(\chi )$ are elongated by $N_{\ j}^{i}(x,y,\chi )=\eta
_{\ j}^{i}(x,y,\chi )~\tilde{N}_{\ j}^{i}(x,y,\chi ).$ We note that in the
formulas defining the coefficients of the metrics (\ref{3lme}) one does not
consider summation on repeating indices which are not ''cross'' one, i.e. $%
\eta _{ij}$ $\tilde{g}_{ij}$ means a simple product $\eta _{ij}$ $\times
\tilde{g}_{ij}$ between deformation function $\eta _{ij}$ and metric
coefficient $\tilde{g}_{ij},$ for any fixed values $i,j,...$ or $a,b,...$ It
is possible to write down the metrics (\ref{3lme}) in ''generic
off--diagonal forms'', see (\ref{3metr}) and (\ref{3ansatz}), for any fixed
value of $\chi .$

\subsection{Nonholonomic Ricci flows and Einstein spaces}

In the theory of Ricci flows, the families of metrics (\ref{3lme}) must
satisfy certain evolution equations on parameter $\chi .$ For normalized
(holonomic) Ricci flows \cite{3per1,3caozhu,3kleiner,3rbook}, with respect
to a coordinate base $\partial _{\underline{\alpha }}=\partial /\partial u^{%
\underline{\alpha }},$ the evolution equations are postulated in the form
\begin{equation}
\frac{\partial }{\partial \chi }g_{\underline{\alpha }\underline{\beta }%
}=-2\ _{\shortmid }R_{\underline{\alpha }\underline{\beta }}+\frac{2r}{5}g_{%
\underline{\alpha }\underline{\beta }},  \label{3feq}
\end{equation}%
where the normalizing factor $r=\int \ _{\shortmid }RdV/dV$ is introduced in
order to preserve the volume $V$ and the metric coefficients $g_{\underline{%
\alpha }\underline{\beta }}$ are parametrized in the form (\ref{3ansatz}),
and $\ _{\shortmid }R_{\underline{\alpha }\underline{\beta }}$ is the Ricci
tensor for the the Levi Civita connection $\nabla .$ In Refs. \cite%
{3nrf1,3nrf2,3nfelf} we discuss in details the N--anholonomic Ricci flows
and prove that the nonholonomic version of (\ref{3feq}) can be proven by a
N--adapted calculus from the Perelman's functionals,
\begin{eqnarray}
\frac{\partial }{\partial \chi }g_{ij} &=&-2\left[ \widehat{R}_{ij}-\lambda
g_{ij}\right] -\frac{\partial }{\partial \chi }(g_{ab}N_{i}^{a}N_{j}^{b}),
\label{3eq1} \\
\frac{\partial }{\partial \chi }g_{ab} &=&-2\ \left( \widehat{R}%
_{ab}-\lambda g_{ab}\right) ,\   \label{3eq2} \\
\ \widehat{R}_{\alpha \beta } &=&0\mbox{ for }\ \alpha \neq \beta ,
\label{3eq3}
\end{eqnarray}%
where $\lambda =r/5.$ We note that the equations (\ref{3eq3}) constrain
nonholonomic Ricci flows to result in symmetric metrics and that we wrote
them with respect to N--adapted frames. A simple class of solutions can be
constructed for the families of N--anholonomic Einstein spaces when
\begin{equation}
\widehat{R}_{ij}-\lambda g_{ij}=0,\quad \widehat{R}_{ab}-\lambda g_{ab}=0
\label{3eq12}
\end{equation}%
and
\begin{equation}
\frac{\partial }{\partial \chi }g_{ij}=-g_{ab}\frac{\partial }{\partial \chi
}(N_{i}^{a}N_{j}^{b}),\quad \frac{\partial }{\partial \chi }g_{ab}=0.
\label{3eq34}
\end{equation}%
Such equations define some effective Einstein metrics subjected to Ricci
flows under evolution of the N--anholonomic structure $N_{i}^{a}$ correlated
with the evolutions of h--metric $g_{ij}$ but $g_{ab}$ stated for a fixed
value of $\chi _{0}.$

For our further considerations, we need the results of two Corollaries (see
Refs. \cite{3caozhu,3nrf1,3nrf2} for detailed proofs and discussions both
for holonomic and nonholonomic manifolds):

\begin{corollary}
\label{corrff}The evolution, for all $\chi \in \lbrack 0,\chi _{0}),$ of
preferred frames on a N--anholonomic manifold
\begin{equation*}
\ \mathbf{e}_{\alpha }(\chi )=\ \mathbf{e}_{\alpha }^{\ \underline{\alpha }%
}(\chi ,u)\partial _{\underline{\alpha }}
\end{equation*}%
is defined by the coefficients
\begin{eqnarray*}
\ \mathbf{e}_{\alpha }^{\ \underline{\alpha }}(\chi ,u) &=&\left[
\begin{array}{cc}
\ e_{i}^{\ \underline{i}}(\chi ,u) & ~N_{i}^{b}(\chi ,u)\ e_{b}^{\
\underline{a}}(\chi ,u) \\
0 & \ e_{a}^{\ \underline{a}}(\chi ,u)%
\end{array}%
\right] ,\  \\
\mathbf{e}_{\ \underline{\alpha }}^{\alpha }(\chi ,u)\ &=&\left[
\begin{array}{cc}
e_{\ \underline{i}}^{i}=\delta _{\underline{i}}^{i} & e_{\ \underline{i}%
}^{b}=-N_{k}^{b}(\chi ,u)\ \ \delta _{\underline{i}}^{k} \\
e_{\ \underline{a}}^{i}=0 & e_{\ \underline{a}}^{a}=\delta _{\underline{a}%
}^{a}%
\end{array}%
\right]
\end{eqnarray*}%
with
\begin{equation*}
\ g_{ij}(\chi )=\ e_{i}^{\ \underline{i}}(\chi ,u)\ e_{j}^{\ \underline{j}%
}(\chi ,u)\eta _{\underline{i}\underline{j}}\mbox{\ and \ }g_{ab}(\chi )=\
e_{a}^{\ \underline{a}}(\chi ,u)\ e_{b}^{\ \underline{b}}(\chi ,u)\eta _{%
\underline{a}\underline{b}},
\end{equation*}%
where $\eta _{\underline{i}\underline{j}}=diag[\pm 1,...\pm 1]$ and $\eta _{%
\underline{a}\underline{b}}=diag[\pm 1,...\pm 1]$ establish the signature of
$\ \mathbf{g}_{\alpha \beta }^{[0]}(u),$ is given by equations
\begin{equation}
\frac{\partial }{\partial \chi }\mathbf{e}_{\ \underline{\alpha }}^{\alpha
}\ =\ \mathbf{g}^{\alpha \beta }~\widehat{\mathbf{R}}_{\beta \gamma }~\
\mathbf{e}_{\ \underline{\alpha }}^{\gamma }.  \label{3aeq5}
\end{equation}
\end{corollary}

For simplicity, we omit formulas for h- and v--decomposition of (\ref{3aeq5}%
).

\begin{corollary}
\label{corrfsc}The scalar curvature (\ref{3sdccurv})
\begin{equation*}
\overleftrightarrow{\mathbf{R}}\doteqdot \mathbf{g}^{\alpha \beta }\widehat{%
\mathbf{R}}_{\alpha \beta }=g^{ij}\widehat{R}_{ij}+h^{ab}\widehat{S}_{ab}=%
\overrightarrow{R}+\overleftarrow{S}
\end{equation*}
for the canonical d--connection on $TM$ satisfies the evolution equations%
\begin{equation}
\frac{\partial \overrightarrow{R}}{\partial \chi }=\widehat{D}_{i}\widehat{D}%
^{i}\overrightarrow{R}+2\widehat{R}_{ij}\widehat{R}^{ij}\mbox{\ and\  }\frac{%
\partial \overleftarrow{S}}{\partial \chi }=\widehat{D}_{a}\widehat{D}^{a}%
\overleftarrow{S}+2\widehat{S}_{ab}\widehat{S}^{ab}.  \label{3scev}
\end{equation}
\end{corollary}

\begin{proof}
It is similar to that for the Levi Civita connection because on $TM$ the
coefficients of canonical d--connection with respect to N-adapted frames are
the same as those for the Levi Civita but decomposed into h-- and
v--components. We note that on $TM$ the Ricci d--tensor $\widehat{\mathbf{R}}%
_{\alpha \beta }$ is symmetric which does not hold true for a general
nonholonomic manifold or vector bundle, see formulas (\ref{3dricci}). The
evolution equations (\ref{3scev}) consist a particular case of more general
formulas for Ricci flows on N--anholonomic manifolds proved in Theorem 4.1
of Ref. \cite{3nrf1} (on $TM$ the distorsion tensor transforming the $\nabla
$ into $\widehat{D}$ is zero).$\square $
\end{proof}

Finally, in this section, we note that a number of geometric ideas and
methods applied in this section were considered in the approaches to the
geometry of nonholonomic spaces and generalized Finsler--Lagrange geometry
elaborated by the schools of G. Vranceanu and R. Miron and by A. Bejancu in
Romania \cite{3vr1,3vr2,3ma1,3ma2,3bej,3bejf}. We emphasize that this way it
is possible to construct geometric models with metric compatible linear
connections which is important for elaborating standard approaches in modern
(non)commutative gravity and string theory \cite{3vncg,3vsgg}. For Finsler
spaces with nontrivial nonmetricity, for instance, for those defined by the
the Berwald and Chern connections, see details in \cite{3bcs}, the physical
theories with local anisotropy are not imbedded into the class of standard
models. It is also a more cumbersome task to elaborate a theory of Ricci
flows of noncompatible metrics and connection structures.

\section{N--Adapted Curve Flows in Vector Bundles}

We formulate the geometry of curve flows adapted to the nonlinear connection
structures constructed by certain classes of canonical lifts from the base
space and nonholonomic frame deformations resulting into constant curvature
coefficients for the canonical d--connection. The case of tangent bundles
will be emphasized as a special one when both $\nabla $ and $\widehat{%
\mathbf{D}}$ can be torsionless.

\subsection{Non--stretching and N--adapted curve flows}

Let us consider a vector bundle $\mathcal{E}=(E,\pi ,F,M),$ $\dim E=$ $n+m$
(in a particular case, $E=TM,$ when $m=n)$ provided with a d--metric $%
\mathbf{g}=[g,h]$ (\ref{3m1}) and N--connection $N_{i}^{a}$ (\ref{3whitney})
structures. A non--stretching curve $\gamma (\tau ,\mathbf{l})$ on $\mathcal{%
E}\mathbf{,}$ where $\tau $ is a parameter and $\mathbf{l}$ is the arclength
of the curve on $\mathcal{E}\mathbf{,}$ is defined with such evolution
d--vector $\mathbf{Y}=\gamma _{\tau }$ and tangent d--vector $\mathbf{X}%
=\gamma _{\mathbf{l}}$ that
\begin{equation}
\mathbf{g(X,X)=}1.  \label{3nsc}
\end{equation}%
The curve $\gamma (\tau ,\mathbf{l})$ swept out a two--dimensional surface
in $T_{\gamma (\tau ,\mathbf{l})}\mathbf{V}\subset T\mathbf{V.}$ If the
geometric objects evolve as Ricci flows on parameter $\chi ,$ we get an
additional parameter for the geometric objects like connections and metrics
and the non--stretching condition (\ref{3nsc}) transforms into $\mathbf{%
g(X,X,}\chi )=1 $ for a family of d--metrics $~^{\chi }\mathbf{g\doteqdot }$
$\mathbf{g}(\chi )=[g(\chi ),h(\chi )]$ and $N_{i}^{a}(\chi )$ which can be
satisfied by certain families of curves, $\gamma (\tau ,\mathbf{l,}\chi
)\doteqdot ~^{\chi }\gamma (\tau ,\mathbf{l})$ (briefly, we shall write only
$~^{\chi }\gamma =\gamma (\chi ))$ and related curve evolution and tangent
vectors, $\mathbf{X(}\chi \mathbf{)}=~^{\chi }\mathbf{X}=~^{\chi }\gamma _{%
\mathbf{l}}$ and $\mathbf{Y(}\chi \mathbf{)}=~^{\chi }\mathbf{X}=~^{\chi
}\gamma _{\tau }.$

We work with families of N--adapted bases (\ref{3dder}) and (\ref{3ddif})
and the connection 1--forms $\mathbf{\Gamma }_{\ \beta }^{\alpha }(\chi )=%
\mathbf{\Gamma }_{\ \beta \gamma }^{\alpha }(\chi )~\mathbf{e}^{\gamma
}(\chi )$ (equivalently, $~^{\chi }\mathbf{\Gamma }_{\ \beta }^{\alpha
}=~^{\chi }\mathbf{\Gamma }_{\ \beta \gamma }^{\alpha }~^{\chi }\mathbf{e}%
^{\gamma })$ with the coefficients $\mathbf{\Gamma }_{\ \beta \gamma
}^{\alpha }(\chi )=$ $~^{\chi }\mathbf{\Gamma }_{\ \beta \gamma }^{\alpha },$
for the canonical d--connection operator $\mathbf{D}(\chi )\doteqdot $ $%
~^{\chi }\mathbf{D}$ (\ref{3candcon}) acting in the form%
\begin{equation}
\mathbf{D}_{\mathbf{X}}\mathbf{e}_{\alpha }=(\mathbf{X\rfloor \Gamma }%
_{\alpha \ }^{\ \gamma })\mathbf{e}_{\gamma }\mbox{ and }\mathbf{D}_{\mathbf{%
Y}}\mathbf{e}_{\alpha }=(\mathbf{Y\rfloor \Gamma }_{\alpha \ }^{\ \gamma })%
\mathbf{e}_{\gamma },  \label{3part01}
\end{equation}%
where ''$\mathbf{\rfloor "}$ denotes the interior product and the indices
are lowered and raised respectively by the d--metric $\mathbf{g}_{\alpha
\beta }=[g_{ij},h_{ab}]$ and its inverse $\mathbf{g}^{\alpha \beta
}=[g^{ij},h^{ab}].$ \footnote{%
For simplicity, we shall omit the Ricci flow parameter if it does not result
in ambiguities.} We note that $\mathbf{D}_{\mathbf{X}(\chi )}=~^{\chi }%
\mathbf{X}^{\alpha }~~^{\chi }\mathbf{D}_{\alpha }$ is the covariant
derivation operator along curve $\gamma (\tau ,\mathbf{l,}\chi ).$ It is
convenient to orient the N--adapted frames to be parallel respectively to
curves $~^{\chi }\gamma $
\begin{eqnarray}
e^{1} &\doteqdot &h\mathbf{X,}\mbox{ for }i=1,\mbox{ and }e^{\widehat{i}},%
\mbox{ where }h\mathbf{g(}h\mathbf{X,}e^{\widehat{i}}\mathbf{)=}0,
\label{3curvframe} \\
\mathbf{e}^{n+1} &\doteqdot &v\mathbf{X,}\mbox{ for }a=n+1,\mbox{ and }%
\mathbf{e}^{\widehat{a}},\mbox{ where }v\mathbf{g(}v\mathbf{X,\mathbf{e}}^{%
\widehat{a}}\mathbf{)=}0,  \notag
\end{eqnarray}%
for $\widehat{i}=2,3,...n$ and $\widehat{a}=n+2,n+3,...,n+m.$ For such
frames, the covariant derivatives of each ''normal'' d--vectors $~^{\chi }%
\mathbf{e}^{\widehat{\alpha }}$ result into the d--vectors adapted to $%
~^{\chi }\gamma ,$
\begin{eqnarray}
~^{\chi }\mathbf{D}_{\mathbf{X}}~^{\chi }e^{\widehat{i}} &\mathbf{=}&\mathbf{%
-}\rho ^{\widehat{i}}\mathbf{(}u,\chi \mathbf{)\ ~^{\chi }X}\mbox{ and }%
~^{\chi }\mathbf{D}_{h\mathbf{X}}~h~^{\chi }\mathbf{X}=\rho ^{\widehat{i}}%
\mathbf{(}u,\chi \mathbf{)\ ~^{\chi }\mathbf{e}}_{\widehat{i}},
\label{3part02} \\
~^{\chi }\mathbf{D}_{\mathbf{X}}~^{\chi }\mathbf{\mathbf{e}}^{\widehat{a}} &%
\mathbf{=}&\mathbf{-}\rho ^{\widehat{a}}\mathbf{(}u,\chi \mathbf{)\ ~^{\chi
}X}\mbox{ and }~^{\chi }\mathbf{D}_{v\mathbf{X}}~v~^{\chi }\mathbf{X}=\rho ^{%
\widehat{a}}\mathbf{(}u,\chi \mathbf{)\ }~^{\chi }e_{\widehat{a}},  \notag
\end{eqnarray}%
which holds for certain classes of functions $~^{\chi }\rho ^{\widehat{i}%
}=\rho ^{\widehat{i}}\mathbf{(}u,\chi \mathbf{)}$ and $~^{\chi }\rho ^{%
\widehat{a}}=\rho ^{\widehat{a}}\mathbf{(}u,\chi \mathbf{).}$ The formulas (%
\ref{3part01}) and (\ref{3part02}) are distinguished into h-- and
v--components for $~^{\chi }\mathbf{X=}h\mathbf{X}(\chi )+v\mathbf{X}(\chi )$
and $~^{\chi }\mathbf{D=(}h\mathbf{D}(\chi ),v\mathbf{D(\chi ))}$ for $%
~^{\chi }\mathbf{D=\{~^{\chi }\Gamma }_{\ \alpha \beta }^{\gamma }\},$ where
$h\mathbf{D}(\chi )=\{~^{\chi }L_{jk}^{i},~^{\chi }L_{bk}^{a}\}$ and $v%
\mathbf{D(}\chi )\mathbf{=\{}~^{\chi }C_{jc}^{i},~^{\chi }C_{bc}^{a}\}.$

Along any curve $\gamma (\chi ),$ we can move differential forms in a
parallel N--adapted form. For instance, $\mathbf{\Gamma }_{\ \mathbf{X}%
}^{\alpha \beta }\doteqdot \mathbf{X\rfloor \Gamma }_{\ }^{\alpha \beta }$
which for families of d--objects is to be written $~^{\chi }\mathbf{\Gamma }%
_{\ \mathbf{X}}^{\alpha \beta }\doteqdot ~^{\chi }\mathbf{X}(\chi )\mathbf{%
\rfloor ~^{\chi }\Gamma }_{\ }^{\alpha \beta }.$ The algebraic
characterization of such spaces, can be obtained if we perform a frame
transform preserving the decomposition (\ref{3whitney}) to an
orthonormalized basis $\mathbf{e}_{\alpha ^{\prime }},$ when
\begin{equation}
\mathbf{e}_{\alpha }\rightarrow A_{\alpha }^{\ \alpha ^{\prime }}(u)\
\mathbf{e}_{\alpha ^{\prime }},~(~^{\chi }\mathbf{e}_{\alpha }\rightarrow
A_{\alpha }^{\ \alpha ^{\prime }}(u,\chi )\ ~^{\chi }\mathbf{e}_{\alpha
^{\prime }}),  \label{3orthbas}
\end{equation}%
called orthonormal d--basis (family of d--bases). In this case, the
coefficients of the d--metric (\ref{3m1}) transform into the Euclidean ones,
$\mathbf{g}_{\alpha ^{\prime }\beta ^{\prime }}=\delta _{\alpha ^{\prime
}\beta ^{\prime }},$ (we can define such frame transform (\ref{3aux4}) when
the the same constant coefficients for d--metric are generated for all
values of parameter $\chi ;$ for such configurations, we do not emphasize
the labels/dependencies on Ricci flow parameter which are present in
d--connection operators and N--connection coefficients). In distinguished
form, we obtain families of two skew matrices%
\begin{equation*}
\mathbf{~^{\chi }\Gamma }_{h\mathbf{X}}^{i^{\prime }j^{\prime }}\doteqdot h%
\mathbf{X(}\chi \mathbf{)\rfloor ~^{\chi }\Gamma }_{\ }^{i^{\prime
}j^{\prime }}=2\ e_{h\mathbf{X}}^{[i^{\prime }}\ \mathbf{~^{\chi }}\rho
^{j^{\prime }]}
\end{equation*}%
and
\begin{equation*}
\mathbf{~^{\chi }\Gamma }_{v\mathbf{X}}^{a^{\prime }b^{\prime }}\doteqdot v%
\mathbf{X\mathbf{(}}\chi \mathbf{\mathbf{)}\rfloor ~^{\chi }\Gamma }_{\
}^{a^{\prime }b^{\prime }}=2\mathbf{\ e}_{v\mathbf{X}}^{[a^{\prime }}\
\mathbf{~^{\chi }}\rho ^{b^{\prime }]}
\end{equation*}%
where
\begin{equation*}
\ e_{h\mathbf{X}}^{i^{\prime }}\doteqdot g(h\mathbf{X,}e^{i^{\prime }})=[1,%
\underbrace{0,\ldots ,0}_{n-1}]\mbox{ and }\ e_{v\mathbf{X}}^{a^{\prime
}}\doteqdot h(v\mathbf{X,}e^{a^{\prime }})=[1,\underbrace{0,\ldots ,0}%
_{m-1}],
\end{equation*}%
we omitted the Ricci flow label performing the constructions according
Remark \ref{rem11}, and
\begin{equation*}
\mathbf{~^{\chi }\Gamma }_{h\mathbf{X\,}i^{\prime }}^{\qquad j^{\prime }}=%
\left[
\begin{array}{cc}
0 & \mathbf{~^{\chi }}\rho ^{j^{\prime }} \\
-\mathbf{~^{\chi }}\rho _{i^{\prime }} & \mathbf{0}_{[h]}%
\end{array}%
\right] \mbox{ and }\mathbf{\Gamma }_{v\mathbf{X\,}a^{\prime }}^{\qquad
b^{\prime }}=\left[
\begin{array}{cc}
0 & \mathbf{~^{\chi }}\rho ^{b^{\prime }} \\
-\mathbf{~^{\chi }}\rho _{a^{\prime }} & \mathbf{0}_{[v]}%
\end{array}%
\right]
\end{equation*}%
with $\mathbf{0}_{[h]}$ and $\mathbf{0}_{[v]}$ being respectively $%
(n-1)\times (n-1)$ and $(m-1)\times (m-1)$ matrices. The above presented
row--matrices and skew--matrices show that locally an total space of a
vector bundle of dimension $n+m,$ with respect to distinguished
orthonormalized frames are characterized algebraically by couples of unit
vectors in $\mathbb{R}^{n}$ and $\mathbb{R}^{m}$ preserved respectively by
the $SO(n-1)$ and $SO(m-1)$ rotation subgroups of the local N--adapted frame
structure group $SO(n)\oplus SO(m).$ The connection matrices $\mathbf{%
~^{\chi }\Gamma }_{h\mathbf{X\,}i^{\prime }}^{\qquad j^{\prime }}$ and $%
\mathbf{~^{\chi }\Gamma }_{v\mathbf{X\,}a^{\prime }}^{\qquad b^{\prime }}$
belong to the orthogonal complements of the corresponding Lie subalgebras
and algebras, $\mathfrak{so}(n-1)\subset \mathfrak{so}(n)$ and $\mathfrak{so}%
(m-1)\subset \mathfrak{so}(m).$

The torsion (\ref{3tors}) and curvature (\ref{3curv}) tensors can be in
orthonormalized component form with respect to (\ref{3curvframe}) mapped
into a distinguished orthotnomalized dual frame (\ref{3orthbas}),%
\begin{equation}
\mathcal{T}^{\alpha ^{\prime }}\doteqdot \mathbf{D}_{\mathbf{X}}\mathbf{e}_{%
\mathbf{Y}}^{\alpha ^{\prime }}-\mathbf{D}_{\mathbf{Y}}\mathbf{e}_{\mathbf{X}%
}^{\alpha ^{\prime }}+\mathbf{e}_{\mathbf{Y}}^{\beta ^{\prime }}\Gamma _{%
\mathbf{X}\beta ^{\prime }}^{\quad \alpha ^{\prime }}-\mathbf{e}_{\mathbf{X}%
}^{\beta ^{\prime }}\Gamma _{\mathbf{Y}\beta ^{\prime }}^{\quad \alpha
^{\prime }}  \label{3mtors}
\end{equation}%
and
\begin{equation}
\mathcal{R}_{\beta ^{\prime }}^{\;\alpha ^{\prime }}(\mathbf{X,Y})=\mathbf{D}%
_{\mathbf{Y}}\Gamma _{\mathbf{X}\beta ^{\prime }}^{\quad \alpha ^{\prime }}-%
\mathbf{D}_{\mathbf{X}}\Gamma _{\mathbf{Y}\beta ^{\prime }}^{\quad \alpha
^{\prime }}+\Gamma _{\mathbf{Y}\beta ^{\prime }}^{\quad \gamma ^{\prime
}}\Gamma _{\mathbf{X}\gamma ^{\prime }}^{\quad \alpha ^{\prime }}-\Gamma _{%
\mathbf{X}\beta ^{\prime }}^{\quad \gamma ^{\prime }}\Gamma _{\mathbf{Y}%
\gamma ^{\prime }}^{\quad \alpha ^{\prime }},  \label{3mcurv}
\end{equation}%
where $\mathbf{e}_{\mathbf{Y}}^{\alpha ^{\prime }}\doteqdot \mathbf{g}(%
\mathbf{Y},\mathbf{e}^{\alpha ^{\prime }})$ and $\Gamma _{\mathbf{Y}\beta
^{\prime }}^{\quad \alpha ^{\prime }}\doteqdot \mathbf{Y\rfloor }\Gamma
_{\beta ^{\prime }}^{\;\alpha ^{\prime }}=\mathbf{g}(\mathbf{e}^{\alpha
^{\prime }},\mathbf{D}_{\mathbf{Y}}\mathbf{e}_{\beta ^{\prime }})$ define
respectively the N--adapted orthonormalized frame row--matrix and the
canonical d--connection skew--matrix in the flow directs, and $\mathcal{R}%
_{\beta ^{\prime }}^{\;\alpha ^{\prime }}(\mathbf{X,Y})\doteqdot \mathbf{g}(%
\mathbf{e}^{\alpha ^{\prime }},[\mathbf{D}_{\mathbf{X}},$ $\mathbf{D}_{%
\mathbf{Y}}]\mathbf{e}_{\beta ^{\prime }})$ is the curvature matrix. Both
torsion and curvature components can be distinguished in h-- and
v--components like (\ref{3dtors}) and (\ref{3dcurv}), by considering
N--adapted decompositions of type
\begin{equation*}
\mathbf{g}=[g,h],\mathbf{e}_{\beta ^{\prime }}=(\mathbf{e}_{j^{\prime
}},e_{b^{\prime }}),\mathbf{e}^{\alpha ^{\prime }}=(e^{i^{\prime
}},e^{a^{\prime }}),\mathbf{X=}h\mathbf{X}+v\mathbf{X,D=(}h\mathbf{D},v%
\mathbf{D).}
\end{equation*}%
Finally, we note that the matrices for torsion (\ref{3mtors}) and curvature (%
\ref{3mcurv}) can be computed for any families, parametrized by $\chi ,$
metric compatible linear connection like the Levi Civita and the canonical
d--connection. For our purposes, in this work, we are interested to define
such a frame of reference with respect to which the curvature tensor has
constant coefficients and the torsion tensor vanishes.

\subsection{N--anholonomic bundles with constant matrix curvature}

For vanishing N--connection torsion and constant matrix curvature of the
canonical d--connection, we get couples of holonomic Riemannian manifolds
and the equations (\ref{3mtors}) and (\ref{3mcurv}) directly encode couples
of bi--Hamiltonian structures, see details in Refs. \cite%
{3vsh,3avw,3saw,3anc2,3anc3}. A well known class of Riemannian manifolds for
which the frame curvature matrix constant consists of the symmetric spaces $%
M=G/H$ for compact semisimple Lie groups $G\supset H.$ A complete
classification and summary of main results on such spaces are given in Refs. %
\cite{3helag,3kob,3sharpe}.

We suppose that the base manifold is a symmetric space $M=hG/SO(n)$ with the
isotropy subgroup $hH=SO(n)\supset O(n)$ and the typical fiber space to be a
symmetric space $F=vG/SO(m)$ with the isotropy subgroup $vH=SO(m)\supset
O(m).$ This means that $hG=SO(n+1)$ and $vG=SO(m+1)$ which is enough for a
study of real holonomic and nonholonomic manifolds and geometric mechanics
models.\footnote{%
it is necessary to consider $hG=SU(n)$ and $vG=SU(m)$ for the geometric
models with spinor and gauge fields} \

Our aim is to solder in a canonic way (like in the N--connection geometry)
the horizontal and vertical symmetric Riemannian spaces of dimension $n$ and
$m$ with a (total) symmetric Riemannian space $V$ of dimension $n+m,$ when $%
V=G/SO(n+m)$ with the isotropy group $H=SO(n+m)\supset O(n+m)$ and $%
G=SO(n+m+1).$ We note that for the just mentioned horizontal, vertical and
total symmetric Riemannian spaces one exists natural settings to Klein
geometry. For instance, the metric tensor $hg=\{\mathring{g}_{ij}\}$ on a
symmetric Riemannian space $M$ is defined by the Cartan--Killing inner
product $<\cdot ,\cdot >_{h}$ on $T_{x}hG\simeq h\mathfrak{g}$ restricted to
the Lie algebra quotient spaces $h\mathfrak{p=}h\mathfrak{g/}h\mathfrak{h,}$
with $T_{x}hH\simeq h\mathfrak{h,}$ where $h\mathfrak{g=}h\mathfrak{h}\oplus
h\mathfrak{p}$ is stated such that there is an involutive automorphism of $%
hG $ under $hH$ is fixed, i.e. $[h\mathfrak{h,}h\mathfrak{p]}\subseteq $ $h%
\mathfrak{p}$ and $[h\mathfrak{p,}h\mathfrak{p]}\subseteq h\mathfrak{h.}$ In
a similar form, we can define the group spaces and related inner products
and\ Lie algebras,%
\begin{eqnarray}
\mbox{for\ }vg &=&\{\mathring{h}_{ab}\},\;<\cdot ,\cdot
>_{v},\;T_{y}vG\simeq v\mathfrak{g,\;}v\mathfrak{p=}v\mathfrak{g/}v\mathfrak{%
h,}\mbox{ with }  \notag \\
T_{y}vH &\simeq &v\mathfrak{h,}v\mathfrak{g=}v\mathfrak{h}\oplus v\mathfrak{%
p,}\mbox{where }\mathfrak{\;}[v\mathfrak{h,}v\mathfrak{p]}\subseteq v%
\mathfrak{p,\;}[v\mathfrak{p,}v\mathfrak{p]}\subseteq v\mathfrak{h;}  \notag
\\
&&  \label{3algstr} \\
\mbox{for\ }\mathbf{g} &=&\{\mathring{g}_{\alpha \beta }\},\;<\cdot ,\cdot
>_{\mathbf{g}},\;T_{(x,y)}G\simeq \mathfrak{g,\;p=g/h,}\mbox{ with }  \notag
\\
T_{(x,y)}H &\simeq &\mathfrak{h,g=h}\oplus \mathfrak{p,}\mbox{where }%
\mathfrak{\;}[\mathfrak{h,p]}\subseteq \mathfrak{p,\;}[\mathfrak{p,p]}%
\subseteq \mathfrak{h.}  \notag
\end{eqnarray}%
We parametrize the metric structure with constant coefficients on $%
V=G/SO(n+m)$ in the form%
\begin{equation*}
\mathring{g}=\mathring{g}_{\alpha \beta }du^{\alpha }\otimes du^{\beta },
\end{equation*}%
where $u^{\alpha }$ are local coordinates and
\begin{equation}
\mathring{g}_{\alpha \beta }=\left[
\begin{array}{cc}
\mathring{g}_{ij}+\mathring{N}_{i}^{a}N_{j}^{b}\mathring{h}_{ab} & \mathring{%
N}_{j}^{e}\mathring{h}_{ae} \\
\mathring{N}_{i}^{e}\mathring{h}_{be} & \mathring{h}_{ab}%
\end{array}%
\right]  \label{3constans}
\end{equation}%
when trivial, constant, N--connection coefficients are computed $\mathring{N}%
_{j}^{e}=\mathring{h}^{eb}\mathring{g}_{jb}$ for any given sets $\mathring{h}%
^{eb}$ and $\mathring{g}_{jb},$ i.e. from the inverse metrics coefficients
defined respectively on $hG=SO(n+1)$ and by off--blocks $(n\times n)$-- and $%
(m\times m)$--terms of the metric $\mathring{g}_{\alpha \beta }.$ As a
result, we define an equivalent d--metric structure of type (\ref{3m1})
\begin{eqnarray}
\mathbf{\mathring{g}} &=&\ \mathring{g}_{ij}\ e^{i}\otimes e^{j}+\ \mathring{%
h}_{ab}\ \mathbf{\mathring{e}}^{a}\otimes \mathbf{\mathring{e}}^{b},
\label{3m1const} \\
e^{i} &=&dx^{i},\ \;\mathbf{\mathring{e}}^{a}=dy^{a}+\mathring{N}%
_{i}^{e}dx^{i}  \notag
\end{eqnarray}%
defining a trivial $(n+m)$--splitting $\mathbf{\mathring{g}=}\mathring{g}%
\mathbf{\oplus _{\mathring{N}}}\mathring{h}\mathbf{\ }$because all
nonholonomy coefficients $\mathring{W}_{\alpha \beta }^{\gamma }$ and
N--connection curvature coefficients $\mathring{\Omega}_{ij}^{a}$ are zero.
In more general form, we can consider any covariant coordinate transforms of
(\ref{3m1const}) preserving the\ $(n+m)$--splitting resulting in any $%
W_{\alpha \beta }^{\gamma }=0$ (\ref{3anhrel}) and $\Omega _{ij}^{a}=0$ (\ref%
{3ncurv}). It should be noted that even such trivial parametrizations define
algebraic classifications of \ symmetric Riemannian spaces of dimension $n+m$
with constant matrix curvature admitting splitting (by certain algebraic
constraints) into symmetric Riemannian subspaces of dimension $n$ and $m,$
also both with constant matrix curvature and introducing the concept of
N--anholonomic Riemannian space of type $\mathbf{\mathring{V}}=[hG=SO(n+1),$
$vG=SO(m+1),\;\mathring{N}_{i}^{e}].$ One can be considered that such
trivially N--anholonomic group spaces have possess a Lie d--algebra symmetry
$\mathfrak{so}_{\mathring{N}}(n+m)\doteqdot \mathfrak{so}(n)\oplus \mathfrak{%
so}(m).$

The simplest generalization on a vector bundle $\mathbf{\mathring{E}}$ is to
consider nonhlonomic distributions on $V=G/SO(n+m)$ defined locally by
families of N--connection coefficients $N_{i}^{a}(x,y,\chi )$ with
nonvanishing $W_{\alpha \beta }^{\gamma }(\chi )$ and $\Omega _{ij}^{a}(\chi
)$ but with constant d--metric coefficients when
\begin{eqnarray}
\mathbf{g}(\chi ) &=&\ \mathring{g}_{ij}\ e^{i}\otimes e^{j}+\ \mathring{h}%
_{ab}\ ~^{\chi }\mathbf{e}^{a}\otimes ~^{\chi }\mathbf{e}^{b},  \label{3m1b}
\\
e^{i} &=&dx^{i},\ ~^{\chi }\mathbf{e}^{a}=dy^{a}+N_{i}^{a}(x,y,\chi )dx^{i}.
\notag
\end{eqnarray}%
This family of metric is very similar to (\ref{3clgs}) but with the
coefficients $\ \mathring{g}_{ij}\ $\ and $\ \mathring{h}_{ab}$ induced by
the corresponding Lie d--algebra structure $\mathfrak{so}_{\mathring{N}%
}(n+m).$ Such spaces transform into families of N--anholonomic
Riemann--Cartan manifolds $\mathbf{\mathring{V}}_{\mathbf{N}}=[hG=SO(n+1),$ $%
vG=SO(m+1),\;N_{i}^{e}]$ with nontrivial N--connection curvature and induced
d--torsion coefficients of the canonical d--connection (see formulas (\ref%
{3dtors}) computed for constant d--metric coefficients and the canonical
d--connection coefficients in (\ref{3candcon})). One has zero curvature for
the canonical d--connection (in general, such spaces are curved ones with
generic off--diagonal metric (\ref{3m1b}) and nonzero curvature tensor for
the Levi Civita connection).\footnote{%
Introducing, constant values for the d--metric coefficients we get zero
coefficients for the canonical d--connection which in its turn results in
zero values of (\ref{3dcurv}).} This allows us to classify the
N--anholonomic manifolds (and vector bundles) as having the same group and
algebraic structures of couples of symmetric Riemannian spaces of dimension $%
n$ and $m$ but nonholonomically soldered to the symmetric Riemannian space
of dimension $n+m.$ With respect to N--adapted orthonormal bases (\ref%
{3orthbas}), with distinguished h-- and v--subspaces, we obtain the same
inner products and group and Lie algebra spaces as in (\ref{3algstr}).

The classification of N--anholonomic vector bundles is almost similar to
that for symmetric Riemannian spaces if we consider that $n=m$ and try to
model tangent bundles of such spaces, provided with N--connection structure.
For instance, we can take a (semi) Riemannian structure with the
N--connection induced by a absolute energy structure like in (\ref{3cnlce})
and with the canonical d--connection structure (\ref{3candcon}), for $\tilde{%
g}_{ef}=\ \mathring{g}_{ab},$ like in (\ref{3aux4}). A straightforward
computation of the canonical d--connection coefficients\footnote{%
on tangent bundles, such d--connections can be defined to be torsionless}
and of d--curvatures for $\;^{\circ }\tilde{g}_{ij}$ and $\;^{\circ }%
\tilde{N}_{\ j}^{i}$ proves that the nonholonomic Riemanian manifold $\left(
M=SO(n+1)/SO(n),\;^{\circ }\mathcal{L}\right) $ possess constant both zero
canonical d--con\-nec\-ti\-on curvature and torsion but with induced
nontrivial N--connection curvature $\;^{\circ }\tilde{\Omega}_{jk}^{i}.$
Such spaces, being tangent to symmetric Riemannian spaces, are classified
similarly to the Riemannian ones with constant matrix curvature, see (\ref%
{3algstr}) for $n=m$ but provided with a nonholonomic structure induced by
generating function $\;^{\circ }\mathcal{L}.$ We can introduce Ricci flows
on parameter $\chi $ when for certain systems of coordinates the metric
coefficients are constant but satisfy the evolution equations (\ref{3eq12})
and (\ref{3eq34}).

\section{Basic Equations for N--anholonomic Curve Flows}

Introducing N--adapted orthonormalized bases, for N--anholonomic spaces of
dimension $n+n,$ with constant curvatures of the canonical d--connection, we
can derive bi--Hamiltonian and vector soliton structures similarly to \cite%
{3anc3,3anc2,3saw}. In symbolic, abstract index form, the constructions for
nonholonomic vector bundles are similar to those for the Riemannian
symmetric--spaces soldered to Klein geometry. We have to distinguish the
horizontal and vertical components of geometric objects and related
equations. The previous bi--Hamiltonian and solitonic constructions were for
an extrinsic approach soldering the Riemannian symmetric--space geometry to
the Klein geometry \cite{3sharpe}. For the N--anhlonomic spaces of dimension
$n+n,$ with constant d--curvatures, similar constructions hold true but we
have to adapt them to the N--connection structure, see Appendix \ref{apb}.

There is an isomorphism between the real space $\mathfrak{so}(n)$ and the
Lie algebra of $n\times n$ skew--symmetric matrices. This allows us to
establish an isomorphism between $h\mathfrak{p}$ $\simeq \mathbb{R}^{n}$ and
the tangent spaces $T_{x}M=\mathfrak{so}(n+1)/$ $\mathfrak{so}(n)$ of the
Riemannian manifold $M=SO(n+1)/$ $SO(n)$ as described by the following
canonical decomposition
\begin{equation*}
h\mathfrak{g}=\mathfrak{so}(n+1)\supset h\mathfrak{p\in }\left[
\begin{array}{cc}
0 & h\mathbf{p} \\
-h\mathbf{p}^{T} & h\mathbf{0}%
\end{array}%
\right] \mbox{\ for\ }h\mathbf{0\in }h\mathfrak{h=so}(n)
\end{equation*}%
with $h\mathbf{p=\{}p^{i^{\prime }}\mathbf{\}\in }\mathbb{R}^{n}$ being the
h--component of the d--vector $\mathbf{p=(}p^{i^{\prime }}\mathbf{,}%
p^{a^{\prime }}\mathbf{)}$ and $h\mathbf{p}^{T}$ mean the transposition of
the row $h\mathbf{p.}$ The Cartan--Killing inner product on $h\mathfrak{g}$
is stated following the rule%
\begin{eqnarray*}
h\mathbf{p\cdot }h\mathbf{p} &\mathbf{=}&\left\langle \left[
\begin{array}{cc}
0 & h\mathbf{p} \\
-h\mathbf{p}^{T} & h\mathbf{0}%
\end{array}%
\right] ,\left[
\begin{array}{cc}
0 & h\mathbf{p} \\
-h\mathbf{p}^{T} & h\mathbf{0}%
\end{array}%
\right] \right\rangle \\
&\mathbf{\doteqdot }&\frac{1}{2}tr\left\{ \left[
\begin{array}{cc}
0 & h\mathbf{p} \\
-h\mathbf{p}^{T} & h\mathbf{0}%
\end{array}%
\right] ^{T}\left[
\begin{array}{cc}
0 & h\mathbf{p} \\
-h\mathbf{p}^{T} & h\mathbf{0}%
\end{array}%
\right] \right\} ,
\end{eqnarray*}%
where $tr$ denotes the trace of the corresponding product of matrices. This
product identifies canonically $h\mathfrak{p}$ $\simeq \mathbb{R}^{n}$ with
its dual $h\mathfrak{p}^{\ast }$ $\simeq \mathbb{R}^{n}.$ In a similar form,
we can consider
\begin{equation*}
v\mathfrak{g}=\mathfrak{so}(m+1)\supset v\mathfrak{p\in }\left[
\begin{array}{cc}
0 & v\mathbf{p} \\
-v\mathbf{p}^{T} & v\mathbf{0}%
\end{array}%
\right] \mbox{\ for\ }v\mathbf{0\in }v\mathfrak{h=so}(m)
\end{equation*}%
with $v\mathbf{p=\{}p^{a^{\prime }}\mathbf{\}\in }\mathbb{R}^{m}$ being the
v--component of the d--vector $\mathbf{p=(}p^{i^{\prime }}\mathbf{,}%
p^{a^{\prime }}\mathbf{)}$ and define the Cartan--Killing inner product $v%
\mathbf{p\cdot }v\mathbf{p\doteqdot }\frac{1}{2}tr\{...\}.$ In general, in
the tangent bundle of a N--anholonomic manifold, we can consider the
Cartan--Killing N--adapted inner product $\mathbf{p\cdot p=}h\mathbf{p\cdot }%
h\mathbf{p+}v\mathbf{p\cdot }v\mathbf{p.}$

Following the introduced Cartan--Killing parametrizations, we analyze the
flow $\gamma (\tau ,\mathbf{l})$ of a non--stretching curve in $\mathbf{V}_{%
\mathbf{N}}=\mathbf{G}/SO(n)\oplus $ $SO(m).$ Let us introduce a family of
coframes $\mathbf{e}\in T_{\gamma }^{\ast }\mathbf{V}_{\mathbf{N}}\otimes (h%
\mathfrak{p\oplus }v\mathfrak{p}),$ which is a N--adapted $\left( SO(n)%
\mathfrak{\oplus }SO(m)\right) $--parallel basis along $\gamma ,$ and its
associated canonical d--con\-nec\-tion 1--form $~^{\chi }\mathbf{\Gamma }\in
T_{\gamma }^{\ast }\mathbf{V}_{\mathbf{N(\chi )}}\otimes (\mathfrak{so}(n)%
\mathfrak{\oplus so}(m)).$ Such d--objects are respectively parametrized:%
\begin{equation*}
\mathbf{e}_{\mathbf{X}}=\mathbf{e}_{h\mathbf{X}}+\mathbf{e}_{v\mathbf{X}},
\end{equation*}%
for
\begin{equation*}
\mathbf{e}_{h\mathbf{X}}=\gamma _{h\mathbf{X}}\rfloor h\mathbf{e=}\left[
\begin{array}{cc}
0 & (1,\overrightarrow{0}) \\
-(1,\overrightarrow{0})^{T} & h\mathbf{0}%
\end{array}%
\right]
\end{equation*}%
and
\begin{equation*}
\mathbf{e}_{v\mathbf{X}}=\gamma _{v\mathbf{X}}\rfloor v\mathbf{e=}\left[
\begin{array}{cc}
0 & (1,\overleftarrow{0}) \\
-(1,\overleftarrow{0})^{T} & v\mathbf{0}%
\end{array}%
\right] ,
\end{equation*}%
where we write $(1,\overrightarrow{0})\in \mathbb{R}^{n},\overrightarrow{0}%
\in \mathbb{R}^{n-1}$ and $(1,\overleftarrow{0})\in \mathbb{R}^{m},%
\overleftarrow{0}\in \mathbb{R}^{m-1};$%
\begin{equation*}
\ ~^{\chi }\mathbf{\Gamma =}\left[ \ ~^{\chi }\mathbf{\Gamma }_{h\mathbf{X}%
},\ ~^{\chi }\mathbf{\Gamma }_{v\mathbf{X}}\right] ,
\end{equation*}%
for
\begin{equation*}
\ ~^{\chi }\mathbf{\Gamma }_{h\mathbf{X}}\mathbf{=}\ ~^{\chi }\gamma _{h%
\mathbf{X}}\rfloor \ ~^{\chi }\mathbf{L=}\left[
\begin{array}{cc}
0 & (0,\overrightarrow{0}) \\
-(0,\overrightarrow{0})^{T} & \ ~^{\chi }\mathbf{L}%
\end{array}%
\right] \in \mathfrak{so}(n+1),
\end{equation*}%
where
\begin{equation*}
\ ~^{\chi }\mathbf{L=}\left[
\begin{array}{cc}
0 & \ ~^{\chi }\overrightarrow{v} \\
-\ ~^{\chi }\overrightarrow{v}^{T} & h\mathbf{0}%
\end{array}%
\right] \in \mathfrak{so}(n),~\ ~^{\chi }\overrightarrow{v}\in \mathbb{R}%
^{n-1},~h\mathbf{0\in }\mathfrak{so}(n-1),
\end{equation*}%
and
\begin{equation*}
\ ~^{\chi }\mathbf{\Gamma }_{v\mathbf{X}}\mathbf{=}\ ~^{\chi }\gamma _{v%
\mathbf{X}}\rfloor \ ~^{\chi }\mathbf{C=}\left[
\begin{array}{cc}
0 & (0,\overleftarrow{0}) \\
-(0,\overleftarrow{0})^{T} & \ ~^{\chi }\mathbf{C}%
\end{array}%
\right] \in \mathfrak{so}(m+1),
\end{equation*}%
where
\begin{equation*}
\ ~^{\chi }\mathbf{C=}\left[
\begin{array}{cc}
0 & \ ~^{\chi }\overleftarrow{v} \\
-\ ~^{\chi }\overleftarrow{v}^{T} & v\mathbf{0}%
\end{array}%
\right] \in \mathfrak{so}(m),~\ ~^{\chi }\overleftarrow{v}\in \mathbb{R}%
^{m-1},~v\mathbf{0\in }\mathfrak{so}(m-1).
\end{equation*}%
The above parametrizations are fixed in order to preserve the $SO(n)$ and $%
SO(m)$ rotation gauge freedoms on the N--adapted coframe and canonical
d--connec\-ti\-on 1--form, distinguished in h- and v--components.

There are defined decompositions of horizontal $SO(n+1)/$ $SO(n)$ matrices
like%
\begin{eqnarray*}
h\mathfrak{p} &\mathfrak{\ni }&\left[
\begin{array}{cc}
0 & h\mathbf{p} \\
-h\mathbf{p}^{T} & h\mathbf{0}%
\end{array}%
\right] =\left[
\begin{array}{cc}
0 & \left( h\mathbf{p}_{\parallel },\overrightarrow{0}\right) \\
-\left( h\mathbf{p}_{\parallel },\overrightarrow{0}\right) ^{T} & h\mathbf{0}%
\end{array}%
\right] \\
&&+\left[
\begin{array}{cc}
0 & \left( 0,h\overrightarrow{\mathbf{p}}_{\perp }\right) \\
-\left( 0,h\overrightarrow{\mathbf{p}}_{\perp }\right) ^{T} & h\mathbf{0}%
\end{array}%
\right] ,
\end{eqnarray*}%
into tangential and normal parts relative to $\mathbf{e}_{h\mathbf{X}}$ via
corresponding decompositions of h--vectors $h\mathbf{p=(}h\mathbf{\mathbf{p}%
_{\parallel },}h\mathbf{\overrightarrow{\mathbf{p}}_{\perp })\in }\mathbb{R}%
^{n}$ relative to $\left( 1,\overrightarrow{0}\right) ,$ when $h\mathbf{%
\mathbf{p}_{\parallel }}$ is identified with $h\mathfrak{p}_{C}$ and $h%
\mathbf{\overrightarrow{\mathbf{p}}_{\perp }}$ is identified with $h%
\mathfrak{p}_{\perp }=h\mathfrak{p}_{C^{\perp }}.$ In a similar form, it is
possible to decompose vertical $SO(m+1)/$ $SO(m)$ matrices,
\begin{eqnarray*}
v\mathfrak{p} &\mathfrak{\ni }&\left[
\begin{array}{cc}
0 & v\mathbf{p} \\
-v\mathbf{p}^{T} & v\mathbf{0}%
\end{array}%
\right] =\left[
\begin{array}{cc}
0 & \left( v\mathbf{p}_{\parallel },\overleftarrow{0}\right) \\
-\left( v\mathbf{p}_{\parallel },\overleftarrow{0}\right) ^{T} & v\mathbf{0}%
\end{array}%
\right] \\
&&+\left[
\begin{array}{cc}
0 & \left( 0,v\overleftarrow{\mathbf{p}}_{\perp }\right) \\
-\left( 0,v\overleftarrow{\mathbf{p}}_{\perp }\right) ^{T} & v\mathbf{0}%
\end{array}%
\right] ,
\end{eqnarray*}%
into tangential and normal parts relative to $\mathbf{e}_{v\mathbf{X}}$ via
corresponding decompositions of h--vectors $v\mathbf{p=(}v\mathbf{\mathbf{p}%
_{\parallel },}v\overleftarrow{\mathbf{\mathbf{p}}}\mathbf{_{\perp })\in }%
\mathbb{R}^{m}$ relative to $\left( 1,\overleftarrow{0}\right) ,$ when $v%
\mathbf{\mathbf{p}_{\parallel }}$ is identified with $v\mathfrak{p}_{C}$ and
$v\overleftarrow{\mathbf{\mathbf{p}}}\mathbf{_{\perp }}$ is identified with $%
v\mathfrak{p}_{\perp }=v\mathfrak{p}_{C^{\perp }}.$

The canonical d--connection induces matrices decomposed with respect to the
flow direction. In the h--direction, we parametrize%
\begin{equation*}
\mathbf{e}_{h\mathbf{Y}}=\gamma _{\tau }\rfloor h\mathbf{e=}\left[
\begin{array}{cc}
0 & \left( h\mathbf{e}_{\parallel },h\overrightarrow{\mathbf{e}}_{\perp
}\right) \\
-\left( h\mathbf{e}_{\parallel },h\overrightarrow{\mathbf{e}}_{\perp
}\right) ^{T} & h\mathbf{0}%
\end{array}%
\right] ,
\end{equation*}%
when $\mathbf{e}_{h\mathbf{Y}}\in h\mathfrak{p,}\left( h\mathbf{e}%
_{\parallel },h\overrightarrow{\mathbf{e}}_{\perp }\right) \in \mathbb{R}%
^{n} $ and $h\overrightarrow{\mathbf{e}}_{\perp }\in \mathbb{R}^{n-1},$ and
\begin{equation}
\ ~^{\chi }\mathbf{\Gamma }_{h\mathbf{Y}}\mathbf{=}\ ~^{\chi }\gamma _{h%
\mathbf{Y}}\rfloor \ ~^{\chi }\mathbf{L=}\left[
\begin{array}{cc}
0 & (0,\overrightarrow{0}) \\
-(0,\overrightarrow{0})^{T} & h\mathbf{\varpi }_{\tau }(\chi )%
\end{array}%
\right] \in \mathfrak{so}(n+1),  \label{3auxaaa}
\end{equation}%
where
\begin{equation*}
h\mathbf{\varpi }_{\tau }(\chi )\mathbf{=}\left[
\begin{array}{cc}
0 & \overrightarrow{\varpi }(\chi ) \\
-\overrightarrow{\varpi }^{T}(\chi ) & h\mathbf{\Theta }(\chi )%
\end{array}%
\right] \in \mathfrak{so}(n),~\overrightarrow{\varpi }(\chi )\in \mathbb{R}%
^{n-1},~h\mathbf{\Theta (}\chi \mathbf{)\in }\mathfrak{so}(n-1).
\end{equation*}%
In the v--direction, we parametrize
\begin{equation*}
\ ~^{\chi }\mathbf{e}_{v\mathbf{Y}}=\ ~^{\chi }\gamma _{\tau }\rfloor v%
\mathbf{e=}\left[
\begin{array}{cc}
0 & \left( v\mathbf{e}_{\parallel },v\overleftarrow{\mathbf{e}}_{\perp
}\right) \\
-\left( v\mathbf{e}_{\parallel },v\overleftarrow{\mathbf{e}}_{\perp }\right)
^{T} & v\mathbf{0}%
\end{array}%
\right] ,
\end{equation*}%
when $\ ~^{\chi }\mathbf{e}_{v\mathbf{Y}}\in v\mathfrak{p,}\left( v\mathbf{e}%
_{\parallel },v\overleftarrow{\mathbf{e}}_{\perp }\right) \in \mathbb{R}^{m}$
and $v\overleftarrow{\mathbf{e}}_{\perp }\in \mathbb{R}^{m-1},$ and
\begin{equation*}
\ ~^{\chi }\mathbf{\Gamma }_{v\mathbf{Y}}\mathbf{=}\ ~^{\chi }\gamma _{v%
\mathbf{Y}}\rfloor \mathbf{C(}\chi \mathbf{)=}\left[
\begin{array}{cc}
0 & (0,\overleftarrow{0}) \\
-(0,\overleftarrow{0})^{T} & v\mathbf{\varpi }_{\tau }\mathbf{(}\chi \mathbf{%
)}%
\end{array}%
\right] \in \mathfrak{so}(m+1),
\end{equation*}%
where
\begin{equation*}
v\mathbf{\varpi }_{\tau }\mathbf{(}\chi \mathbf{)=}\left[
\begin{array}{cc}
0 & \overleftarrow{\varpi }\mathbf{(}\chi \mathbf{)} \\
-\overleftarrow{\varpi }^{T}\mathbf{(}\chi \mathbf{)} & v\mathbf{\Theta (}%
\chi \mathbf{)}%
\end{array}%
\right] \in \mathfrak{so}(m),~\overleftarrow{\varpi }(\chi \mathbf{)}\in
\mathbb{R}^{m-1},~v\mathbf{\Theta \mathbf{(}\chi )\in }\mathfrak{so}(m-1).
\end{equation*}%
The components $h\mathbf{e}_{\parallel }(\chi )$ and $h\overrightarrow{%
\mathbf{e}}_{\perp }(\chi )$ correspond to the decomposition
\begin{equation*}
\mathbf{~^{\chi }e}_{h\mathbf{Y}}=h\mathbf{g(~^{\chi }\gamma }_{\tau
},~^{\chi }\mathbf{\gamma }_{\mathbf{l}},\chi \mathbf{)~^{\chi }e}_{h\mathbf{%
X}}+\mathbf{(~^{\chi }\gamma }_{\tau })_{\perp }\rfloor h\mathbf{e}_{\perp }%
\mathbf{(}\chi )
\end{equation*}%
into tangential and normal parts relative to $\mathbf{~^{\chi }e}_{h\mathbf{X%
}}.$ In a similar form, one considers $v\mathbf{e}_{\parallel }(\chi )(\chi
) $ and $v\overleftarrow{\mathbf{e}}_{\perp }(\chi )$ corresponding to the
decomposition
\begin{equation*}
\mathbf{~^{\chi }e}_{v\mathbf{Y}}=v\mathbf{g(\gamma }_{\tau },\mathbf{\gamma
}_{\mathbf{l}},\chi \mathbf{)~^{\chi }e}_{v\mathbf{X}}+\mathbf{(~^{\chi
}\gamma }_{\tau })_{\perp }\rfloor v\mathbf{e}_{\perp }(\chi ).
\end{equation*}

Using the above stated matrix parametrizations, we get%
\begin{eqnarray}
\left[ \mathbf{~^{\chi }e}_{h\mathbf{X}},\mathbf{~^{\chi }e}_{h\mathbf{Y}}%
\right] &=&-\left[
\begin{array}{cc}
0 & 0 \\
0 & h\mathbf{e}_{\perp }(\chi )%
\end{array}%
\right] \in \mathfrak{so}(n+1),  \label{3aux41} \\
\mbox{ \ for \ }h\mathbf{e}_{\perp }(\chi ) &=&\left[
\begin{array}{cc}
0 & h\overrightarrow{\mathbf{e}}_{\perp }(\chi ) \\
-(h\overrightarrow{\mathbf{e}}_{\perp }(\chi ))^{T} & h\mathbf{0}%
\end{array}%
\right] \in \mathfrak{so}(n);  \notag \\
\left[ \mathbf{~^{\chi }\Gamma }_{h\mathbf{Y}},\mathbf{~^{\chi }e}_{h\mathbf{%
Y}}\right] &=&-\left[
\begin{array}{cc}
0 & \left( 0,\mathbf{~^{\chi }}\overrightarrow{\varpi }\right) \\
-\left( 0,\mathbf{~^{\chi }}\overrightarrow{\varpi }\right) ^{T} & 0%
\end{array}%
\right] \in h\mathfrak{p}_{\perp };  \notag \\
\left[ \mathbf{~^{\chi }\Gamma }_{h\mathbf{X}},\mathbf{~^{\chi }e}_{h\mathbf{%
Y}}\right] &=&\left[
\begin{array}{cc}
0 & \mathbf{~^{\chi }}\left( \overrightarrow{v}\cdot h\overrightarrow{%
\mathbf{e}}_{\perp },-h\mathbf{e}_{\parallel }\overrightarrow{v}\right) \\
\mathbf{~^{\chi }}\left( -\overrightarrow{v}\cdot h\overrightarrow{\mathbf{e}%
}_{\perp },h\mathbf{e}_{\parallel }\overrightarrow{v}\right) ^{T} & h\mathbf{%
0}%
\end{array}%
\right]  \notag \\
&\in &h\mathfrak{p};  \notag
\end{eqnarray}%
and
\begin{eqnarray}
\left[ \mathbf{~^{\chi }e}_{v\mathbf{X}},\mathbf{~^{\chi }e}_{v\mathbf{Y}}%
\right] &=&-\left[
\begin{array}{cc}
0 & 0 \\
0 & v\mathbf{e}_{\perp }(\chi )%
\end{array}%
\right] \in \mathfrak{so}(m+1),  \label{3aux41a} \\
\mbox{ \ for \ }v\mathbf{e}_{\perp }(\chi ) &=&\left[
\begin{array}{cc}
0 & v\overrightarrow{\mathbf{e}}_{\perp }(\chi ) \\
-(v\overrightarrow{\mathbf{e}}_{\perp })^{T}(\chi ) & v\mathbf{0}%
\end{array}%
\right] \in \mathfrak{so}(m);  \notag \\
\left[ \mathbf{~^{\chi }\Gamma }_{v\mathbf{Y}},\mathbf{~^{\chi }e}_{v\mathbf{%
Y}}\right] &=&-\left[
\begin{array}{cc}
0 & \left( 0,\mathbf{~^{\chi }}\overleftarrow{\varpi }\right) \\
-\left( 0,\mathbf{~^{\chi }}\overleftarrow{\varpi }\right) ^{T} & 0%
\end{array}%
\right] \in v\mathfrak{p}_{\perp };  \notag \\
\left[ \mathbf{~^{\chi }\Gamma }_{v\mathbf{X}},\mathbf{~^{\chi }e}_{v\mathbf{%
Y}}\right] &=&\left[
\begin{array}{cc}
0 & \mathbf{~^{\chi }}\left( \overleftarrow{v}\cdot v\overleftarrow{\mathbf{e%
}}_{\perp },-v\mathbf{e}_{\parallel }\overleftarrow{v}\right) \\
\mathbf{~^{\chi }}\left( -\overleftarrow{v}\cdot v\overleftarrow{\mathbf{e}}%
_{\perp },v\mathbf{e}_{\parallel }\overleftarrow{v}\right) ^{T} & v\mathbf{0}%
\end{array}%
\right]  \notag \\
&\in &v\mathfrak{p}.  \notag
\end{eqnarray}

We can use formulas (\ref{3aux41}) and (\ref{3aux41a}) in order to write the
structure equations (\ref{3mtors}) and (\ref{3mcurv}) in terms of N--adapted
curve flow operators soldered to the geometry Klein N--anholonomic spaces
using the relations (\ref{3aux33}). One obtains respectively the $\mathbf{G}$%
--invariant N--adapted torsion and curvature generated by the canonical
d--connection,
\begin{equation}
\mathbf{T}(\gamma _{\tau },\gamma _{\mathbf{l}})=\left( \mathbf{D}_{\mathbf{X%
}}\gamma _{\tau }-\mathbf{D}_{\mathbf{Y}}\gamma _{\mathbf{l}}\right) \rfloor
\mathbf{e=D}_{\mathbf{X}}\mathbf{e}_{\mathbf{Y}}-\mathbf{D}_{\mathbf{Y}}%
\mathbf{e}_{\mathbf{X}}+\left[ \mathbf{\Gamma }_{\mathbf{X}},\mathbf{e}_{%
\mathbf{Y}}\right] -\left[ \mathbf{\Gamma }_{\mathbf{Y}},\mathbf{e}_{\mathbf{%
X}}\right]  \label{3torscf}
\end{equation}%
and
\begin{equation}
\mathbf{R}(\gamma _{\tau },\gamma _{\mathbf{l}})\mathbf{e=}\left[ \mathbf{D}%
_{\mathbf{X}},\mathbf{D}_{\mathbf{Y}}\right] \mathbf{e=D}_{\mathbf{X}}%
\mathbf{\Gamma }_{\mathbf{Y}}-\mathbf{D}_{\mathbf{Y}}\mathbf{\Gamma }_{%
\mathbf{X}}+\left[ \mathbf{\Gamma }_{\mathbf{X}},\mathbf{\Gamma }_{\mathbf{Y}%
}\right]  \label{3curvcf}
\end{equation}%
where $\mathbf{e}_{\mathbf{X}}\doteqdot \gamma _{\mathbf{l}}\rfloor \mathbf{%
e,}$ $\mathbf{e}_{\mathbf{Y}}\doteqdot \gamma _{\mathbf{\tau }}\rfloor
\mathbf{e,}$ $\mathbf{\Gamma }_{\mathbf{X}}\doteqdot \gamma _{\mathbf{l}%
}\rfloor \mathbf{\Gamma }$ and $\mathbf{\Gamma }_{\mathbf{Y}}\doteqdot
\gamma _{\mathbf{\tau }}\rfloor \mathbf{\Gamma .}$ The formulas (\ref%
{3torscf}) and (\ref{3curvcf}) are equivalent, respectively, to (\ref{3dtors}%
) and (\ref{3dcurv}). In general, $\mathbf{T}(\gamma _{\tau },\gamma _{%
\mathbf{l}})\neq 0$ and $\mathbf{R}(\gamma _{\tau },\gamma _{\mathbf{l}})%
\mathbf{e}$ can not be defined to have constant matrix coefficients with
respect to a N--adapted basis. For N--anholonomic spaces with dimensions $%
n=m,$ we have $^{\mathcal{E}}\mathbf{T}(\gamma _{\tau },\gamma _{\mathbf{l}%
})=0$ and $^{\mathcal{E}}\mathbf{R}(\gamma _{\tau },\gamma _{\mathbf{l}})%
\mathbf{e}$ defined by constant, or vanishing, d--curvature coefficients
(see discussions related to formulas (\ref{3dcurvtb}) and (\ref{3candcontm}%
)). For such cases, we can consider the h-- and v--components of (\ref%
{3torscf}) and (\ref{3curvcf}) in a similar manner as for symmetric
Riemannian spaces but with the canonical d--connection instead of the Levi
Civita one. One obtains, respectively,%
\begin{eqnarray}
0 &=&\left( \mathbf{D}_{h\mathbf{X}}\gamma _{\tau }-\mathbf{D}_{h\mathbf{Y}%
}\gamma _{\mathbf{l}}\right) \rfloor h\mathbf{e}  \label{3torseq} \\
&\mathbf{=}&\mathbf{D}_{h\mathbf{X}}\mathbf{e}_{h\mathbf{Y}}-\mathbf{D}_{h%
\mathbf{Y}}\mathbf{e}_{h\mathbf{X}}+\left[ \mathbf{L}_{h\mathbf{X}},\mathbf{e%
}_{h\mathbf{Y}}\right] -\left[ \mathbf{L}_{h\mathbf{Y}},\mathbf{e}_{h\mathbf{%
X}}\right] ;  \notag \\
0 &=&\left( \mathbf{D}_{v\mathbf{X}}\gamma _{\tau }-\mathbf{D}_{v\mathbf{Y}%
}\gamma _{\mathbf{l}}\right) \rfloor v\mathbf{e}  \notag \\
&\mathbf{=}&\mathbf{D}_{v\mathbf{X}}\mathbf{e}_{v\mathbf{Y}}-\mathbf{D}_{v%
\mathbf{Y}}\mathbf{e}_{v\mathbf{X}}+\left[ \mathbf{C}_{v\mathbf{X}},\mathbf{e%
}_{v\mathbf{Y}}\right] -\left[ \mathbf{C}_{v\mathbf{Y}},\mathbf{e}_{v\mathbf{%
X}}\right] ,  \notag
\end{eqnarray}%
and%
\begin{eqnarray}
h\mathbf{R}(\gamma _{\tau },\gamma _{\mathbf{l}})h\mathbf{e} &\mathbf{=}&%
\left[ \mathbf{D}_{h\mathbf{X}},\mathbf{D}_{h\mathbf{Y}}\right] h\mathbf{e=D}%
_{h\mathbf{X}}\mathbf{L}_{h\mathbf{Y}}-\mathbf{D}_{h\mathbf{Y}}\mathbf{L}_{h%
\mathbf{X}}+\left[ \mathbf{L}_{h\mathbf{X}},\mathbf{L}_{h\mathbf{Y}}\right]
\label{3curveq} \\
v\mathbf{R}(\gamma _{\tau },\gamma _{\mathbf{l}})v\mathbf{e} &\mathbf{=}&%
\left[ \mathbf{D}_{v\mathbf{X}},\mathbf{D}_{v\mathbf{Y}}\right] v\mathbf{e=D}%
_{v\mathbf{X}}\mathbf{C}_{v\mathbf{Y}}-\mathbf{D}_{v\mathbf{Y}}\mathbf{C}_{v%
\mathbf{X}}+\left[ \mathbf{C}_{v\mathbf{X}},\mathbf{C}_{v\mathbf{Y}}\right] .
\notag
\end{eqnarray}

Following the N--adapted curve flow parametrizations (\ref{3aux41}) and (\ref%
{3aux41a}), the equations (\ref{3torseq}) and (\ref{3curveq}) are written
\begin{eqnarray}
0 &=&\mathbf{D}_{h\mathbf{X}}h\mathbf{e}_{\parallel }+\overrightarrow{v}%
\cdot h\overrightarrow{\mathbf{e}}_{\perp },~0=\mathbf{D}_{v\mathbf{X}}v%
\mathbf{e}_{\parallel }+\overleftarrow{v}\cdot v\overleftarrow{\mathbf{e}}%
_{\perp },;  \label{3torseqd} \\
~0 &=&\overrightarrow{\varpi }-h\mathbf{e}_{\parallel }\overrightarrow{v}+%
\mathbf{D}_{h\mathbf{X}}h\overrightarrow{\mathbf{e}}_{\perp },~0=%
\overleftarrow{\varpi }-v\mathbf{e}_{\parallel }\overleftarrow{v}+\mathbf{D}%
_{v\mathbf{X}}v\overleftarrow{\mathbf{e}}_{\perp };  \notag
\end{eqnarray}%
and%
\begin{eqnarray}
\mathbf{D}_{h\mathbf{X}}\overrightarrow{\varpi }-\mathbf{D}_{h\mathbf{Y}}%
\overrightarrow{v}+\overrightarrow{v}\rfloor h\mathbf{\Theta } &\mathbf{=}&h%
\overrightarrow{\mathbf{e}}_{\perp },~\mathbf{D}_{v\mathbf{X}}\overleftarrow{%
\varpi }-\mathbf{D}_{v\mathbf{Y}}\overleftarrow{v}+\overleftarrow{v}\rfloor v%
\mathbf{\Theta =}v\overleftarrow{\mathbf{e}}_{\perp };  \notag \\
\mathbf{D}_{h\mathbf{X}}h\mathbf{\Theta -}\overrightarrow{v}\otimes
\overrightarrow{\varpi }+\overrightarrow{\varpi }\otimes \overrightarrow{v}
&=&0,~\mathbf{D}_{v\mathbf{X}}v\mathbf{\Theta -}\overleftarrow{v}\otimes
\overleftarrow{\varpi }+\overleftarrow{\varpi }\otimes \overleftarrow{v}=0.
\label{3curveqd}
\end{eqnarray}%
The tensor and interior products, for instance, for the h--components, are
defined in the form: $\otimes $ denotes the outer product of pairs of
vectors ($1\times n$ row matrices), producing $n\times n$ matrices $%
\overrightarrow{A}\otimes \overrightarrow{B}=\overrightarrow{A}^{T}%
\overrightarrow{B},$ and $\rfloor $ denotes multiplication of $n\times n$
matrices on vectors ($1\times n$ row matrices); one holds the properties $%
\overrightarrow{A}\rfloor \left( \overrightarrow{B}\otimes \overrightarrow{C}%
\right) =\left( \overrightarrow{A}\cdot \overrightarrow{B}\right)
\overrightarrow{C}$ which is the transpose of the standard matrix product on
column vectors, and $\left( \overrightarrow{B}\otimes \overrightarrow{C}%
\right) \overrightarrow{A}=\left( \overrightarrow{C}\cdot \overrightarrow{A}%
\right) \overrightarrow{B}.$ Here we note that similar formulas hold for the
v--components but, for instance, we have to change, correspondingly, $%
n\rightarrow m$ and $\overrightarrow{A}\rightarrow \overleftarrow{A}.$

The variables $\mathbf{e}_{\parallel }$ and $\mathbf{\Theta ,}$ written in
h-- and v--components, can be expressed corresponding in terms of variables $%
\overrightarrow{v},\overrightarrow{\varpi },h\overrightarrow{\mathbf{e}}%
_{\perp }$ and $\overleftarrow{v},\overleftarrow{\varpi },v\overleftarrow{%
\mathbf{e}}_{\perp }$ (see respectively the first two equations in (\ref%
{3torseqd}) and the last two equations in (\ref{3curveqd})),%
\begin{equation*}
h\mathbf{e}_{\parallel }=-\mathbf{D}_{h\mathbf{X}}^{-1}(\overrightarrow{v}%
\cdot h\overrightarrow{\mathbf{e}}_{\perp }),~v\mathbf{e}_{\parallel }=-%
\mathbf{D}_{v\mathbf{X}}^{-1}(\overleftarrow{v}\cdot v\overleftarrow{\mathbf{%
e}}_{\perp }),
\end{equation*}%
and%
\begin{equation*}
h\mathbf{\Theta =D}_{h\mathbf{X}}^{-1}\left( \overrightarrow{v}\otimes
\overrightarrow{\varpi }-\overrightarrow{\varpi }\otimes \overrightarrow{v}%
\right) ,~v\mathbf{\Theta =D}_{v\mathbf{X}}^{-1}\left( \overleftarrow{v}%
\otimes \overleftarrow{\varpi }-\overleftarrow{\varpi }\otimes
\overleftarrow{v}\right) .
\end{equation*}%
Substituting these values, correspondingly, in the last two equations in (%
\ref{3torseqd}) and in the first two equations in (\ref{3curveqd}), we
express
\begin{equation*}
\overrightarrow{\varpi }=-\mathbf{D}_{h\mathbf{X}}h\overrightarrow{\mathbf{e}%
}_{\perp }-\mathbf{D}_{h\mathbf{X}}^{-1}(\overrightarrow{v}\cdot h%
\overrightarrow{\mathbf{e}}_{\perp })\overrightarrow{v},~\overleftarrow{%
\varpi }=-\mathbf{D}_{v\mathbf{X}}v\overleftarrow{\mathbf{e}}_{\perp }-%
\mathbf{D}_{v\mathbf{X}}^{-1}(\overleftarrow{v}\cdot v\overleftarrow{\mathbf{%
e}}_{\perp })\overleftarrow{v},
\end{equation*}%
contained in the h-- and v--flow equations respectively on $\overrightarrow{v%
}$ and $\overleftarrow{v},$ considered as scalar components when $\mathbf{D}%
_{h\mathbf{Y}}\overrightarrow{v}=\overrightarrow{v}_{\tau }$ and $\mathbf{D}%
_{h\mathbf{Y}}\overleftarrow{v}=\overleftarrow{v}_{\tau },$
\begin{eqnarray}
\overrightarrow{v}_{\tau } &=&\mathbf{D}_{h\mathbf{X}}\overrightarrow{\varpi
}-\overrightarrow{v}\rfloor \mathbf{D}_{h\mathbf{X}}^{-1}\left(
\overrightarrow{v}\otimes \overrightarrow{\varpi }-\overrightarrow{\varpi }%
\otimes \overrightarrow{v}\right) -\overrightarrow{R}h\overrightarrow{%
\mathbf{e}}_{\perp },  \label{3floweq} \\
\overleftarrow{v}_{\tau } &=&\mathbf{D}_{v\mathbf{X}}\overleftarrow{\varpi }-%
\overleftarrow{v}\rfloor \mathbf{D}_{v\mathbf{X}}^{-1}\left( \overleftarrow{v%
}\otimes \overleftarrow{\varpi }-\overleftarrow{\varpi }\otimes
\overleftarrow{v}\right) -\overleftarrow{S}v\overleftarrow{\mathbf{e}}%
_{\perp },  \notag
\end{eqnarray}%
where the scalar curvatures of the canonical d--connection, $\overrightarrow{%
R}$ and $\overleftarrow{S}$ are defined by formulas (\ref{3sdccurv}). For
symmetric Riemannian spaces like $SO(n+1)/SO(n)\simeq S^{n},$ the value $%
\overrightarrow{R}$ is just the scalar curvature $\chi =1,$\ see \cite{3anc3}%
. On tangent bundles, it is possible that $\overrightarrow{R}$ and $%
\overleftarrow{S}$ are certain zero or nonzero constants with the h--part
equivalent to the base scalar curvature.

For Ricci flows of geometric objects, the curve flow evolution equations (%
\ref{3floweq}) contain additional dependencies on parameter $\chi ,$
\begin{eqnarray}
\mathbf{~^{\chi }}\overrightarrow{v}_{\tau } &=&\mathbf{~^{\chi }D}_{h%
\mathbf{X}}\overrightarrow{\varpi }(\chi )-\overrightarrow{R}(\chi )~h%
\overrightarrow{\mathbf{e}}_{\perp }(\chi )  \notag \\
&&-\mathbf{~^{\chi }}\overrightarrow{v}\rfloor \mathbf{~^{\chi }D}_{h\mathbf{%
X}}^{-1}\left( \mathbf{~^{\chi }}\overrightarrow{v}\otimes \overrightarrow{%
\varpi }(\chi )-\overrightarrow{\varpi }(\chi )\otimes \mathbf{~^{\chi }}%
\overrightarrow{v}\right) ,  \notag \\
\mathbf{~^{\chi }}\overleftarrow{v}_{\tau } &=&\mathbf{~^{\chi }D}_{v\mathbf{%
X}}\overleftarrow{\varpi }-\overleftarrow{S}(\chi )v\overleftarrow{\mathbf{e}%
}_{\perp }(\chi ),  \notag \\
&&-\mathbf{~^{\chi }}\overleftarrow{v}\rfloor \mathbf{~^{\chi }D}_{v\mathbf{X%
}}^{-1}\left( \mathbf{~^{\chi }}\overleftarrow{v}\otimes \overleftarrow{%
\varpi }(\chi )-\overleftarrow{\varpi }(\chi )\otimes \mathbf{~^{\chi }}%
\overleftarrow{v}(\chi )\right)  \notag
\end{eqnarray}%
where the scalar curvatures evolve following formulas (\ref{3scev}) for
Ricci flows.

The above presented considerations consist the proof of

\begin{lemma}
\label{lemcfe}There are canonical lifts of Ricci flows from a (semi)
Riemannian manfiold $M$ to $TM$ when certain families of constant constant
curvature matrix coefficients for the canonical d--connections define
families of N--adapted Hamiltonian sympletic operators,
\begin{eqnarray}
h\mathcal{J}(\chi ) &=&\mathbf{~^{\chi }D}_{h\mathbf{X}}+\mathbf{~^{\chi }D}%
_{h\mathbf{X}}^{-1}\left( \mathbf{~^{\chi }}\overrightarrow{v}\cdot \right)
\mathbf{~^{\chi }}\overrightarrow{v},  \label{3sop} \\
v\mathcal{J}(\chi ) &=&\mathbf{~^{\chi }D}_{v\mathbf{X}}+\mathbf{~^{\chi }D}%
_{v\mathbf{X}}^{-1}\left( \mathbf{~^{\chi }}\overleftarrow{v}\cdot \right)
\mathbf{~^{\chi }}\overleftarrow{v},  \notag
\end{eqnarray}%
and cosympletic operators%
\begin{eqnarray}
h\mathcal{H}(\chi ) &\doteqdot &\mathbf{~^{\chi }D}_{h\mathbf{X}}+\mathbf{%
~^{\chi }}\overrightarrow{v}\rfloor \mathbf{~^{\chi }D}_{h\mathbf{X}%
}^{-1}\left( \mathbf{~^{\chi }}\overrightarrow{v}\wedge \right)
\label{3csop} \\
v\mathcal{H}(\chi ) &\doteqdot &\mathbf{~^{\chi }D}_{v\mathbf{X}}+\mathbf{%
~^{\chi }}\overleftarrow{v}\rfloor \mathbf{~^{\chi }D}_{v\mathbf{X}%
}^{-1}\left( \mathbf{~^{\chi }}\overleftarrow{v}\wedge \right) ,  \notag
\end{eqnarray}%
where, for instance, $\overrightarrow{A}\wedge \overrightarrow{B}=%
\overrightarrow{A}\otimes \overrightarrow{B}-\overrightarrow{B}\otimes $ $%
\overrightarrow{A}.$\
\end{lemma}

For any fixed value $\chi =\chi _{0},$ the formulas for this Lemma transform
into similar ones from Ref. \cite{3vsh}. The properties of operators (\ref%
{3sop}) and (\ref{3csop}) are defined by

\begin{theorem}
\label{mr1}The Ricci flows of d--operators $\mathbf{~^{\chi }}\mathcal{J=}%
\left( h\mathcal{J}(\chi ),v\mathcal{J}(\chi )\right) $ and $\mathbf{~^{\chi
}}\mathcal{H=}\left( h\mathcal{H}(\chi ),v\mathcal{H}(\chi )\right) $ $\ $%
are defined respectively by $\left( O(n-1),O(m-1)\right) $--in\-va\-riant
Hamiltonian sympletic and cosympletic d--operators with respect to the
corresponding Hamiltonian d--variables $\left( \mathbf{~^{\chi }}%
\overrightarrow{v},\mathbf{~^{\chi }}\overleftarrow{v}\right) .$ Such
d--operators defines the Hamiltonian form for the curve and Ricci flows
equations on N--anholonomic tangent bundles with constant d--connection
curvature: the h--flows are given by%
\begin{eqnarray}
\mathbf{~^{\chi }}\overrightarrow{v}_{\tau } &=&h\mathcal{H}\left( \mathbf{%
~^{\chi }}\overrightarrow{\varpi },\chi \right) -\overrightarrow{R}(\chi )~h%
\overrightarrow{\mathbf{e}}_{\perp }(\chi )  \notag \\
&=&h\mathfrak{R}\left( h\overrightarrow{\mathbf{e}}_{\perp }(\chi ),\chi
\right) -\overrightarrow{R}(\chi )~h\overrightarrow{\mathbf{e}}_{\perp
}(\chi ),  \notag \\
\mathbf{~^{\chi }}\overrightarrow{\varpi } &=&h\mathcal{J}\left( h%
\overrightarrow{\mathbf{e}}_{\perp }(\chi ),\chi \right) ,  \label{3hhfeq1}
\\
\frac{\partial \overrightarrow{R}}{\partial \chi } &=&\widehat{D}_{i}%
\widehat{D}^{i}\overrightarrow{R}+2\widehat{R}_{ij}\widehat{R}^{ij};  \notag
\end{eqnarray}%
the v--flows are given by
\begin{eqnarray}
\mathbf{~^{\chi }}\overleftarrow{v}_{\tau } &=&v\mathcal{H}\left( \mathbf{%
~^{\chi }}\overleftarrow{\varpi }(\chi ),\chi \right) -\overleftarrow{S}%
(\chi )~v\overleftarrow{\mathbf{e}}_{\perp }(\chi )  \notag \\
&=&v\mathfrak{R}\left( v\overleftarrow{\mathbf{e}}_{\perp }(\chi ),\chi
\right) -\overleftarrow{S}(\chi )~v\overleftarrow{\mathbf{e}}_{\perp }(\chi
),  \notag \\
\mathbf{~^{\chi }}\overleftarrow{\varpi } &=&v\mathcal{J}\left( v%
\overleftarrow{\mathbf{e}}_{\perp }(\chi ),\chi \right) ,  \label{3vhfeq1} \\
\frac{\partial \overleftarrow{S}}{\partial \chi } &=&\widehat{D}_{a}\widehat{%
D}^{a}\overleftarrow{S}+2\widehat{S}_{ab}\widehat{S}^{ab};  \notag
\end{eqnarray}%
where the so--called Ricci flows of heriditary recursion d--operator has the
respective h-- and v--components
\begin{equation}
h\mathfrak{R}(\chi )=h\mathcal{H}(\chi )\circ h\mathcal{J}(\chi )%
\mbox{ \
and \ }v\mathfrak{R}(\chi )=v\mathcal{H}(\chi )\circ v\mathcal{J}(\chi )
\label{3reqop}
\end{equation}%
and $\widehat{D}_{i}\widehat{D}^{i}\overrightarrow{R}=\widehat{D}_{a}%
\widehat{D}^{a}\overleftarrow{S}=0$ and $\widehat{R}_{ij}\widehat{R}%
^{ij}=const,\widehat{S}_{ab}\widehat{S}^{ab}=const$ for lifts to constant
curvature matrices.
\end{theorem}

\begin{proof}
One follows from the Lemma \ref{lemcfe}, Corollaries \ref{corrff} and \ref%
{corrfsc} and (\ref{3floweq}). In a detailed form, for holonomic structures,
it is given in Ref. \cite{3saw} and discussed in \cite{3anc3}. $\square $
\end{proof}

Finally, we note that for any fixed value $\chi _{0}$ we get a Theorem from %
\cite{3vsh}, on curve flows in symmetric Riemannian spaces, which has
generalizations for curve flows on generalized Lagrange and Finsler spaces %
\cite{3avw}.

\section{Curve Flows and Solitonic Hierarchies for Ricci Flows}

The final aim of this paper is to prove that for any nonholonomic Ricci flow
system \ we can define naturally a family of N--adapted bi--Hamiltonian flow
hierarchies inducing anholonomic solitonic configurations.

\subsection{Formulation of the main theorem}

Following a usual solitonic techniques generalized for N--anholonomic
spaces, see details in Ref. \cite{3vsh,3avw,3anc2,3anc3}, the recursion
h--operators from (\ref{3reqop}),%
\begin{eqnarray}
h\mathfrak{R}(\chi ) &=&\mathbf{~^{\chi }D}_{h\mathbf{X}}\left( \mathbf{%
~^{\chi }D}_{h\mathbf{X}}+\mathbf{~^{\chi }D}_{h\mathbf{X}}^{-1}\left(
\mathbf{~^{\chi }}\overrightarrow{v}\cdot \right) \mathbf{~^{\chi }}%
\overrightarrow{v}\right)  \label{3reqoph} \\
&&+\mathbf{~^{\chi }}\overrightarrow{v}\rfloor \mathbf{~^{\chi }D}_{h\mathbf{%
X}}^{-1}\left( \mathbf{~^{\chi }}\overrightarrow{v}\wedge \mathbf{~^{\chi }D}%
_{h\mathbf{X}}\right)  \notag \\
&=&\mathbf{~^{\chi }D}_{h\mathbf{X}}^{2}+|\mathbf{~^{\chi }D}_{h\mathbf{X}%
}|^{2}+\mathbf{~^{\chi }D}_{h\mathbf{X}}^{-1}\left( \mathbf{~^{\chi }}%
\overrightarrow{v}\cdot \right) \mathbf{~^{\chi }}\overrightarrow{v}_{%
\mathbf{l}}-\mathbf{~^{\chi }}\overrightarrow{v}\rfloor \mathbf{~^{\chi }D}%
_{h\mathbf{X}}^{-1}(\mathbf{~^{\chi }}\overrightarrow{v}_{\mathbf{l}}\wedge
),  \notag
\end{eqnarray}%
generate a \ family of horizontal hierarchies of commuting Hamiltonian
vector fields $h\overrightarrow{\mathbf{e}}_{\perp }^{(k)}(\chi )$ starting
from $h\overrightarrow{\mathbf{e}}_{\perp }^{(0)}(\chi )=\mathbf{~^{\chi }}%
\overrightarrow{v}_{\mathbf{l}}$ given by the infinitesimal generator of $%
\mathbf{l}$--translations in terms of arclength $\mathbf{l}$ along the curve
(we use a boldface $\mathbf{l}$ in order to emphasized that the curve is on
a N--anholonomic manifold when the geometric objects are subjected to Ricci
flows). A family of vertical hierarchies of commuting vector fields $v%
\overleftarrow{\mathbf{e}}_{\perp }^{(k)}(\chi )$ starting from $v%
\overleftarrow{\mathbf{e}}_{\perp }^{(0)}(\chi )$ $=\mathbf{~^{\chi }}%
\overleftarrow{v}_{\mathbf{l}}$ is generated by the recursion v--operators%
\begin{eqnarray}
v\mathfrak{R}(\chi ) &=&\mathbf{~^{\chi }D}_{v\mathbf{X}}\left( \mathbf{%
~^{\chi }D}_{v\mathbf{X}}+\mathbf{~^{\chi }D}_{v\mathbf{X}}^{-1}\left(
\mathbf{~^{\chi }}\overleftarrow{v}\cdot \right) \mathbf{~^{\chi }}%
\overleftarrow{v}\right)  \label{3reqopv} \\
&&+\mathbf{~^{\chi }}\overleftarrow{v}\rfloor \mathbf{~^{\chi }D}_{v\mathbf{X%
}}^{-1}\left( \mathbf{~^{\chi }}\overleftarrow{v}\wedge \mathbf{~^{\chi }D}%
_{v\mathbf{X}}\right)  \notag \\
&=&\mathbf{~^{\chi }D}_{v\mathbf{X}}^{2}+|\mathbf{~^{\chi }D}_{v\mathbf{X}%
}|^{2}+\mathbf{~^{\chi }D}_{v\mathbf{X}}^{-1}\left( \mathbf{~^{\chi }}%
\overleftarrow{v}\cdot \right) \mathbf{~^{\chi }}\overleftarrow{v}_{\mathbf{l%
}}-\mathbf{~^{\chi }}\overleftarrow{v}\rfloor \mathbf{~^{\chi }D}_{v\mathbf{X%
}}^{-1}(\mathbf{~^{\chi }}\overleftarrow{v}_{\mathbf{l}}\wedge ).  \notag
\end{eqnarray}%
There are related hierarchies, generated by adjoint operators $\mathfrak{R}%
^{\ast }=(h\mathfrak{R}^{\ast },$ $v\mathfrak{R}^{\ast }),$ of involuntive
Hamiltonian h--covector fields $\overrightarrow{\varpi }^{(k)}=\delta \left(
hH^{(k)}\right) /\delta \overrightarrow{v}$ in terms of Hamiltonians $%
hH=hH^{(k)}(\overrightarrow{v},\overrightarrow{v}_{\mathbf{l}},%
\overrightarrow{v}_{2\mathbf{l}},...)$ starting from $\overrightarrow{\varpi
}^{(0)}=\overrightarrow{v},hH^{(0)}=\frac{1}{2}|\overrightarrow{v}|^{2}$ and
of involutive Hamiltonian v--covector fields $\overleftarrow{\varpi }%
^{(k)}=\delta \left( vH^{(k)}\right) /$ $\delta \overleftarrow{v}$ in terms
of Hamiltonians $vH=vH^{(k)}(\overleftarrow{v},\overleftarrow{v}_{\mathbf{l}%
},\overleftarrow{v}_{2\mathbf{l}},...)$ starting from $\overleftarrow{\varpi
}^{(0)}=\overleftarrow{v},vH^{(0)}=\frac{1}{2}|\overleftarrow{v}|^{2}.$

The relations between different type families of hierarchies are established
correspondingly by formulas%
\begin{equation*}
h\overrightarrow{\mathbf{e}}_{\perp }^{(k)}(\chi )=h\mathcal{H}\left(
\overrightarrow{\varpi }^{(k)},\overrightarrow{\varpi }^{(k+1)},\chi \right)
=h\mathcal{J}\left( h\overrightarrow{\mathbf{e}}_{\perp }^{(k)},\chi \right)
\end{equation*}%
and
\begin{equation*}
v\overleftarrow{\mathbf{e}}_{\perp }^{(k)}(\chi )=v\mathcal{H}\left(
\overleftarrow{\varpi }^{(k)},\overleftarrow{\varpi }^{(k+1)},\chi \right) =v%
\mathcal{J}\left( v\overleftarrow{\mathbf{e}}_{\perp }^{(k)},\chi \right) ,
\end{equation*}%
where $k=0,1,2,....$ All hierarchies (horizontal, vertical and their adjoint
ones) have a typical mKdV scaling symmetry, for instance, $\mathbf{%
l\rightarrow \lambda l}$ and $\mathbf{~^{\chi }}\overrightarrow{v}%
\rightarrow \mathbf{\lambda }^{-1}\mathbf{~^{\chi }}\overrightarrow{v}$
under which the values $h\overrightarrow{\mathbf{e}}_{\perp }^{(k)}(\chi )$
and $hH^{(k)}(\chi )$ have scaling weight $2+2k,$ while $\overrightarrow{%
\varpi }^{(k)}(\chi )$ has scaling weight $1+2k.$ Following the above
presented considerations, we prove

\begin{corollary}
\label{3c2} There are Ricci flow families of N--adapted hierarchies of
distinguished horizontal and vertical commuting bi--Hamiltonian flows,
correspondingly, on $\overrightarrow{v}$ and $\overleftarrow{v}$ associated
to the recursion d--operator (\ref{3reqop}) given by families of $%
O(n-1)\oplus O(m-1)$ --invariant d--vector h--evolution equations,%
\begin{eqnarray*}
\mathbf{~^{\chi }}\overrightarrow{v}_{\tau } &=&h\overrightarrow{\mathbf{e}}%
_{\perp }^{(k+1)}(\chi )-\overrightarrow{R}(\chi )~h\overrightarrow{\mathbf{e%
}}_{\perp }^{(k)}(\chi )=h\mathcal{H}\left( \delta \left( hH^{(k,%
\overrightarrow{R})}(\chi )\right) /\delta \mathbf{~^{\chi }}\overrightarrow{%
v}\right) \\
&=&\left( h\mathcal{J}(\chi )\right) ^{-1}\left( \delta \left( hH^{(k+1,%
\overrightarrow{R})}(\chi )\right) /\delta \mathbf{~^{\chi }}\overrightarrow{%
v}\right) ,
\end{eqnarray*}%
with families of horizontal Hamiltonians
\begin{equation*}
hH^{(k+1,\overrightarrow{R})}(\chi )=hH^{(k+1,\overrightarrow{R})}(\chi )-%
\overrightarrow{R}(\chi )~hH^{(k,\overrightarrow{R})}(\chi ),
\end{equation*}%
and v--evolution equations
\begin{eqnarray*}
\mathbf{~^{\chi }}\overleftarrow{v}_{\tau } &=&v\overleftarrow{\mathbf{e}}%
_{\perp }^{(k+1)}(\chi )-\overleftarrow{S}(\chi )~v\overleftarrow{\mathbf{e}}%
_{\perp }^{(k)}(\chi )=v\mathcal{H}\left( \delta \left( vH^{(k,%
\overleftarrow{S})}(\chi )\right) /\delta \mathbf{~^{\chi }}\overleftarrow{v}%
\right) \\
&=&\left( v\mathcal{J}(\chi )\right) ^{-1}\left( \delta \left( vH^{(k+1,%
\overleftarrow{S})}(\chi )\right) /\delta \mathbf{~^{\chi }}\overleftarrow{v}%
\right) ,
\end{eqnarray*}%
and with \ families of vertical Hamiltonians
\begin{equation*}
vH^{(k+1,\overleftarrow{S})}(\chi )=vH^{(k+1,\overleftarrow{S})}(\chi )-%
\overleftarrow{S}(\chi )~vH^{(k,\overleftarrow{S})}(\chi ),
\end{equation*}
for $k=0,1,2,.....$ The Ricci flows of d--operators $\mathcal{H}(\chi )$ and
$\mathcal{J}$ $(\chi )$ are N--adapted and mutually compatible from which we
can construct a family of alternative (explicit) Hamilton d--operators
\begin{equation*}
~^{a}\mathcal{H}(\chi )\mathcal{=\mathbf{~^{\chi }}H\circ J}(\chi )\circ
\mathbf{~^{\chi }}\mathcal{H=}\mathbf{~^{\chi }}\mathfrak{R\circ }\mathcal{H}%
(\chi ).
\end{equation*}
\end{corollary}

The Main Result of this work is formulated in the form:

\begin{theorem}
\label{3mt} For any vector/tangent bundle with Ricci flows of d--metric
structures, one can be defined a family of hierarchies of bi-Hamiltonian
N--adapted flows of curves $\mathbf{~^{\chi }}\gamma (\tau ,\mathbf{l}%
)=h\gamma (\tau ,\mathbf{l},\chi )+v\gamma (\tau ,\mathbf{l},\chi )$
described by families of geometric nonholonomic map equations:

The $0$ flows are defined as convective (travelling wave) maps%
\begin{eqnarray}
\mathbf{~^{\chi }}\gamma _{\tau } &=&\mathbf{~^{\chi }}\gamma _{\mathbf{l}},%
\mbox{\ distinguished as }  \label{3trmap} \\
\left( h\gamma \right) _{\tau }(\chi ) &=&\left( h\gamma \right) _{h\mathbf{X%
}}(\chi )\mbox{\ and \ }\left( v\gamma \right) _{\tau }(\chi )=\left(
v\gamma \right) _{v\mathbf{X}}(\chi ).  \notag
\end{eqnarray}

There are families of +1 flows defined as Ricci flows of non--stretching
mKdV maps%
\begin{eqnarray}
-\mathbf{~^{\chi }}\left( h\gamma \right) _{\tau } &=&\left[ \mathbf{~^{\chi
}D}_{h\mathbf{X}}^{2} +\frac{3}{2}\left| \mathbf{~^{\chi }D}_{h\mathbf{X}%
}\left( h\gamma \right) _{h\mathbf{X}}(\chi )\right| _{h\mathbf{g}}^{2} %
\right] \left( h\gamma \right) _{h\mathbf{X}}(\chi ),  \label{31map} \\
-\mathbf{~^{\chi }}\left( v\gamma \right) _{\tau } &=& \left[ \mathbf{%
~^{\chi }D}_{v\mathbf{X}}^{2} +\frac{3}{2}\left| \mathbf{~^{\chi }D}_{v%
\mathbf{X}}\left( v\gamma \right) _{v\mathbf{X}}(\chi )\right| _{v\mathbf{g}%
}^{2} \right] \left( v\gamma \right) _{v\mathbf{X}}(\chi ),  \notag
\end{eqnarray}%
and the families of +2,...flows as higher order analogs.

There are also families of -1 flows defined by the kernels of recursion
operators (\ref{3reqoph}) and (\ref{3reqopv}) inducing non--stretching maps%
\begin{equation}
\mathbf{~^{\chi }D}_{h\mathbf{Y}}\left( h\gamma \right) _{h\mathbf{X}}(\chi
)=0\mbox{\ and \
}\mathbf{~^{\chi }D}_{v\mathbf{Y}}\left( v\gamma \right) _{v\mathbf{X}}(\chi
)=0.  \label{3-1map}
\end{equation}
\end{theorem}

Proof is outlined in Appendix \ref{3ssp}.

\subsection{Nonholonomic mKdV and SG hierarchies}

Let us consider some explicit constructions when families of solitonic
hierarchies are derived following the conditions of Theorem \ref{3mt}.

The h--flow and v--flow equations resulting from (\ref{3-1map}) are%
\begin{equation}
\mathbf{~^{\chi }}\overrightarrow{v}_{\tau }=-\overrightarrow{R}(\chi )~h%
\overrightarrow{\mathbf{e}}_{\perp }(\chi )\mbox{ \ and \ }\mathbf{~^{\chi }}%
\overleftarrow{v}_{\tau }=-\overleftarrow{S}(\chi )~v\overleftarrow{\mathbf{e%
}}_{\perp }(\chi ),  \label{3deveq}
\end{equation}%
when, respectively,%
\begin{equation*}
0=\mathbf{~^{\chi }}\overrightarrow{\varpi }=-\mathbf{~^{\chi }D}_{h\mathbf{X%
}}h\overrightarrow{\mathbf{e}}_{\perp }(\chi )+h\mathbf{e}_{\parallel }(\chi
)\mathbf{~^{\chi }}\overrightarrow{v},\mathbf{~^{\chi }D}_{h\mathbf{X}}h%
\mathbf{e}_{\parallel }(\chi )=h\overrightarrow{\mathbf{e}}_{\perp }(\chi
)\cdot \mathbf{~^{\chi }}\overrightarrow{v}
\end{equation*}%
and
\begin{equation*}
0=\mathbf{~^{\chi }}\overleftarrow{\varpi }=-\mathbf{~^{\chi }D}_{v\mathbf{X}%
}v\overleftarrow{\mathbf{e}}_{\perp }(\chi )+v\mathbf{e}_{\parallel }(\chi )%
\mathbf{~^{\chi }}\overleftarrow{v},~\mathbf{~^{\chi }D}_{v\mathbf{X}}v%
\mathbf{e}_{\parallel }(\chi )=v\overleftarrow{\mathbf{e}}_{\perp }(\chi
)\cdot \mathbf{~^{\chi }}\overleftarrow{v}.
\end{equation*}%
The d--flow equations possess horizontal and vertical conservation laws%
\begin{equation*}
\mathbf{~^{\chi }D}_{h\mathbf{X}}\left( (h\mathbf{e}_{\parallel }(\chi
))^{2}+|h\overrightarrow{\mathbf{e}}_{\perp }(\chi )|^{2}\right) =0,
\end{equation*}%
for $(h\mathbf{e}_{\parallel }(\chi ))^{2}+|h\overrightarrow{\mathbf{e}}%
_{\perp }(\chi )|^{2}=<h\mathbf{e}_{\tau }(\chi ),h\mathbf{e}_{\tau }(\chi
)>_{h\mathfrak{p}}=|\left( h\gamma \right) _{\tau }(\chi )|_{h\mathbf{g}%
}^{2},$ and
\begin{equation*}
\mathbf{~^{\chi }D}_{v\mathbf{Y}}\left( (v\mathbf{e}_{\parallel }(\chi
))^{2}+|v\overleftarrow{\mathbf{e}}_{\perp }(\chi )|^{2}\right) =0,
\end{equation*}%
for $(v\mathbf{e}_{\parallel }(\chi ))^{2}+|v\overleftarrow{\mathbf{e}}%
_{\perp }(\chi )|^{2}=<v\mathbf{e}_{\tau }(\chi ),v\mathbf{e}_{\tau }(\chi
)>_{v\mathfrak{p}}=|\left( v\gamma \right) _{\tau }(\chi )|_{v\mathbf{g}%
}^{2}.$ This corresponds to
\begin{equation*}
\mathbf{~^{\chi }D}_{h\mathbf{X}}|\left( h\gamma \right) _{\tau }(\chi )|_{h%
\mathbf{g}}^{2}=0\mbox{ \ and \ }\mathbf{~^{\chi }D}_{v\mathbf{X}}|\left(
v\gamma \right) _{\tau }(\chi )|_{v\mathbf{g}}^{2}=0.
\end{equation*}%
In general, such laws are more sophisticate than those on (semi) Riemannian
spaces because of nonholonomic constraints resulting in non--symmetric Ricci
tensors and different types of identities. But for the geometries modelled
for dimensions $n=m$ with canonical d--connections, we get similar h-- and
v--components of the conservation law equations as on symmetric Riemannian
spaces.\footnote{%
We note that the problem of formulating conservation laws on
N--anholo\-no\-mic spaces (in particular, on nonholonomic vector bundles) in
analyzed in Ref. \cite{3vsgg}.}

It is possible to rescale conformally the variable $\tau $ and obtain $%
|\left( h\gamma \right) _{\tau }(\chi )|_{h\mathbf{g}}^{2}$ $=1$ and (it can
be for other rescale) $|\left( v\gamma \right) _{\tau }(\chi )|_{v\mathbf{g}%
}^{2}=1,$ i.e.
\begin{equation*}
(h\mathbf{e}_{\parallel }(\chi ))^{2}+|h\overrightarrow{\mathbf{e}}_{\perp
}(\chi )|^{2}=1\mbox{ \ and \ }(v\mathbf{e}_{\parallel }(\chi ))^{2}+|v%
\overleftarrow{\mathbf{e}}_{\perp }(\chi )|^{2}=1.
\end{equation*}%
We can express $h\mathbf{e}_{\parallel }(\chi )$ and $h\overrightarrow{%
\mathbf{e}}_{\perp }(\chi )$ in terms of $\mathbf{~^{\chi }}\overrightarrow{v%
}$ and its derivatives and, similarly, we can express $v\mathbf{e}%
_{\parallel }(\chi )$ and $v\overleftarrow{\mathbf{e}}_{\perp }(\chi )$ in
terms of $\mathbf{~^{\chi }}\overleftarrow{v}$ and its derivatives, which
follows from (\ref{3deveq}). The N--adapted wave map equations describing
the -1 flows reduce to a system of two independent nonlocal evolution
equations for the h-- and v--components, parametrized by $\chi ,$%
\begin{eqnarray*}
\mathbf{~^{\chi }}\overrightarrow{v}_{\tau } &=&-\mathbf{~^{\chi }D}_{h%
\mathbf{X}}^{-1}\left( \sqrt{\overrightarrow{R}^{2}(\chi )-|\mathbf{~^{\chi }%
}\overrightarrow{v}_{\tau }|^{2}}\mathbf{~^{\chi }}\overrightarrow{v}\right)
, \\
\mathbf{~^{\chi }}\overleftarrow{v}_{\tau } &=&-\mathbf{~^{\chi }D}_{v%
\mathbf{X}}^{-1}\left( \sqrt{\overleftarrow{S}^{2}(\chi )-|\mathbf{~^{\chi }}%
\overleftarrow{v}_{\tau }|^{2}}\mathbf{~^{\chi }}\overleftarrow{v}\right) .
\end{eqnarray*}%
For N--anholonomic spaces of constant scalar d--curvatures, we can rescale
the equations on $\tau $ to the case when the terms $\overrightarrow{R}%
^{2}(\chi )$ and $\overleftarrow{S}^{2}(\chi )$ are constant, and the
evolution equations transform into a system of hyperbolic d--vector
equations,%
\begin{eqnarray}
\mathbf{~^{\chi }D}_{h\mathbf{X}}(\mathbf{~^{\chi }}\overrightarrow{v}_{\tau
}) &=&-\sqrt{1-|\mathbf{~^{\chi }}\overrightarrow{v}_{\tau }|^{2}}\mathbf{%
~^{\chi }}\overrightarrow{v},  \label{3heq} \\
\mathbf{~^{\chi }D}_{v\mathbf{X}}(\mathbf{~^{\chi }}\overleftarrow{v}_{\tau
}) &=&-\sqrt{1-|\mathbf{~^{\chi }}\overleftarrow{v}_{\tau }|^{2}}~\mathbf{%
~^{\chi }}\overleftarrow{v},  \notag
\end{eqnarray}%
where $\mathbf{~^{\chi }D}_{h\mathbf{X}}=\partial _{h\mathbf{l}}$ and $%
\mathbf{~^{\chi }D}_{v\mathbf{X}}=\partial _{v\mathbf{l}}$ are usual partial
derivatives on direction $\mathbf{l=}h\mathbf{l+}v\mathbf{l}$ with $\mathbf{%
~^{\chi }}\overrightarrow{v}_{\tau }$ and $\mathbf{~^{\chi }}\overleftarrow{v%
}_{\tau }$ considered as scalar functions for the covariant derivatives $%
\mathbf{~^{\chi }D}_{h\mathbf{X}}$ and $\mathbf{~^{\chi }D}_{v\mathbf{X}}$
defined by the canonical d--connection. It also follows that $h%
\overrightarrow{\mathbf{e}}_{\perp }(\chi )$ and $v\overleftarrow{\mathbf{e}}%
_{\perp }(\chi )$ obey corresponding vector sine--Gordon (SG) equations%
\begin{equation}
\left( \sqrt{(1-|h\overrightarrow{\mathbf{e}}_{\perp }(\chi )|^{2})^{-1}}%
~\partial _{h\mathbf{l}}(h\overrightarrow{\mathbf{e}}_{\perp }(\chi
))\right) _{\tau }=-h\overrightarrow{\mathbf{e}}_{\perp }(\chi )
\label{3sgeh}
\end{equation}%
and
\begin{equation}
\left( \sqrt{(1-|v\overleftarrow{\mathbf{e}}_{\perp }(\chi )|^{2})^{-1}}%
~\partial _{v\mathbf{l}}(v\overleftarrow{\mathbf{e}}_{\perp }(\chi ))\right)
_{\tau }=-v\overleftarrow{\mathbf{e}}_{\perp }(\chi ).  \label{3sgev}
\end{equation}

The above presented formulas and Corollary \ref{3c2} imply

\begin{conclusion}
The Ricci flow families of recursion d--operators $\mathbf{~^{\chi }}%
\mathfrak{R}=(h\mathfrak{R}(\chi )\mathfrak{,}h\mathfrak{R}(\chi ))$ (\ref%
{3reqop}), see (\ref{3reqoph}) and (\ref{3reqopv}), generate two hierarchies
of vector mKdV symmetries: the first one is horizontal,
\begin{eqnarray}
\mathbf{~^{\chi }}\overrightarrow{v}_{\tau }^{(0)} &=&\mathbf{~^{\chi }}%
\overrightarrow{v}_{h\mathbf{l}},~\mathbf{~^{\chi }}\overrightarrow{v}_{\tau
}^{(1)}=h\mathfrak{R}(\overrightarrow{v}_{h\mathbf{l}},\chi )=\mathbf{%
~^{\chi }}\overrightarrow{v}_{3h\mathbf{l}}+\frac{3}{2}|\overrightarrow{v}%
(\chi )|^{2}~\mathbf{~^{\chi }}\overrightarrow{v}_{h\mathbf{l}},  \notag \\
\mathbf{~^{\chi }}\overrightarrow{v}_{\tau }^{(2)} &=&h\mathfrak{R}^{2}(%
\overrightarrow{v}_{h\mathbf{l}},\chi )=\mathbf{~^{\chi }}\overrightarrow{v}%
_{5h\mathbf{l}}+\frac{5}{2}\left( |\overrightarrow{v}(\chi )|^{2}~\mathbf{%
^{\chi }}\overrightarrow{v}_{2h\mathbf{l}}\right) _{h\mathbf{l}}  \notag \\
&&+\frac{5}{2}\left( (|\overrightarrow{v}(\chi )|^{2})_{h\mathbf{l~}h\mathbf{%
l}}+|\mathbf{~^{\chi }}\overrightarrow{v}_{h\mathbf{l}}|^{2}+\frac{3}{4}|%
\overrightarrow{v}(\chi )|^{4}\right) \mathbf{~^{\chi }}\overrightarrow{v}_{h%
\mathbf{l}}  \notag \\
&&-\frac{1}{2}|\mathbf{~^{\chi }}\overrightarrow{v}_{h\mathbf{l}}|^{2}~%
\overrightarrow{v}(\chi ),...,  \label{3mkdv1a}
\end{eqnarray}%
with all such terms commuting with the -1 flow
\begin{equation}
(\mathbf{~^{\chi }}\overrightarrow{v}_{\tau })^{-1}=h\overrightarrow{\mathbf{%
e}}_{\perp }(\chi )  \label{3mkdv1b}
\end{equation}%
associated to the vector SG equation (\ref{3sgeh}); the second one is
vertical,
\begin{eqnarray}
\mathbf{~^{\chi }}\overleftarrow{v}_{\tau }^{(0)} &=&\mathbf{~^{\chi }}%
\overleftarrow{v}_{v\mathbf{l}},\mathbf{~^{\chi }}\overleftarrow{v}_{\tau
}^{(1)}=v\mathfrak{R}(\overleftarrow{v}_{v\mathbf{l}},\chi )=\mathbf{~^{\chi
}}\overleftarrow{v}_{3v\mathbf{l}}+\frac{3}{2}|\overleftarrow{v}(\chi )|^{2}~%
\mathbf{~^{\chi }}\overleftarrow{v}_{v\mathbf{l}},  \notag \\
\mathbf{~^{\chi }}\overleftarrow{v}_{\tau }^{(2)} &=&v\mathfrak{R}^{2}(%
\overleftarrow{v}_{v\mathbf{l}},\chi )=\mathbf{~^{\chi }}\overleftarrow{v}%
_{5v\mathbf{l}}+\frac{5}{2}\left( |\overleftarrow{v}(\chi )|^{2}\mathbf{%
~^{\chi }}\overleftarrow{v}_{2v\mathbf{l}}\right) _{v\mathbf{l}}  \notag \\
&&+\frac{5}{2}\left( (|\overleftarrow{v}(\chi )|^{2})_{v\mathbf{l~}v\mathbf{l%
}}+|\mathbf{~^{\chi }}\overleftarrow{v}_{v\mathbf{l}}|^{2}+\frac{3}{4}|%
\overleftarrow{v}(\chi )|^{4}\right) \mathbf{~^{\chi }}\overleftarrow{v}_{v%
\mathbf{l}}  \notag \\
&&-\frac{1}{2}|\mathbf{~^{\chi }}\overleftarrow{v}_{v\mathbf{l}}|^{2}~%
\overleftarrow{v}(\chi ),...,  \label{3mkdv2a}
\end{eqnarray}%
with all such terms commuting with the -1 flow
\begin{equation}
(\mathbf{~^{\chi }}\overleftarrow{v}_{\tau })^{-1}=v\overleftarrow{\mathbf{e}%
}_{\perp }(\chi )  \label{3mkdv2b}
\end{equation}%
associated to the vector SG equation (\ref{3sgev}).
\end{conclusion}

In its turn, using the above Conclusion, we derive that the family of
adjoint d--operators $\mathfrak{R}^{\ast }(\chi )=\mathbf{~^{\chi }}\mathcal{%
J\circ H}(\chi )$ generates corresponding families of horizontal hierarchies
of Hamiltonians (for simplicity, here we omit labels/dependences on $\chi $),%
\begin{eqnarray}
hH^{(0)} &=&\frac{1}{2}|\overrightarrow{v}|^{2},~hH^{(1)}=-\frac{1}{2}|%
\overrightarrow{v}_{h\mathbf{l}}|^{2}+\frac{1}{8}|\overrightarrow{v}|^{4},
\label{3hhh} \\
hH^{(2)} &=&\frac{1}{2}|\overrightarrow{v}_{2h\mathbf{l}}|^{2}-\frac{3}{4}|%
\overrightarrow{v}|^{2}~|\overrightarrow{v}_{h\mathbf{l}}|^{2}-\frac{1}{2}%
\left( \overrightarrow{v}\cdot \overrightarrow{v}_{h\mathbf{l}}\right) +%
\frac{1}{16}|\overrightarrow{v}|^{6},...,  \notag
\end{eqnarray}%
and of vertical hierarchies of Hamiltonians%
\begin{eqnarray}
vH^{(0)} &=&\frac{1}{2}|\overleftarrow{v}|^{2},~vH^{(1)}=-\frac{1}{2}|%
\overleftarrow{v}_{v\mathbf{l}}|^{2}+\frac{1}{8}|\overleftarrow{v}|^{4},
\label{3hhv} \\
vH^{(2)} &=&\frac{1}{2}|\overleftarrow{v}_{2v\mathbf{l}}|^{2}-\frac{3}{4}|%
\overleftarrow{v}|^{2}~|\overleftarrow{v}_{v\mathbf{l}}|^{2}-\frac{1}{2}%
\left( \overleftarrow{v}\cdot \overleftarrow{v}_{v\mathbf{l}}\right) +\frac{1%
}{16}|\overleftarrow{v}|^{6},...,  \notag
\end{eqnarray}%
all of which are conserved densities for respective horizontal and vertical
-1 flows and determining higher conservation laws for the corresponding
hyperbolic equations (\ref{3sgeh}) and (\ref{3sgev}).

The above presented horizontal equations (\ref{3sgeh}), (\ref{3mkdv1a}), (%
\ref{3mkdv1b}) and (\ref{3hhh}) and of vertical equations (\ref{3sgev}), (%
\ref{3mkdv2a}), (\ref{3mkdv2b}) and (\ref{3hhv}) have similar mKdV scaling
symmetries but on different parameters $\lambda _{h}$ and $\lambda _{v}$
because, in general, there are two independent values of scalar curvatures $%
\overrightarrow{R}$ and $\overleftarrow{S},$ see (\ref{3sdccurv}). The
horizontal scaling symmetries are $h\mathbf{l\rightarrow }\lambda _{h}h%
\mathbf{l,}\overrightarrow{v}\rightarrow \left( \lambda _{h}\right) ^{-1}%
\overrightarrow{v}$ and $\tau \rightarrow \left( \lambda _{h}\right)
^{1+2k}, $ for $k=-1,0,1,2,...$ For the vertical scaling symmetries, one has
$v\mathbf{l\rightarrow }\lambda _{v}v\mathbf{l,}\overleftarrow{v}\rightarrow
\left( \lambda _{v}\right) ^{-1}\overleftarrow{v}$ and $\tau \rightarrow
\left( \lambda _{v}\right) ^{1+2k},$ for $k=-1,0,1,2,...$

Finally, we disucss how exact solutions for the Einstein equations%
\begin{equation}
~^{\chi }\widehat{\mathbf{R}}_{\alpha \beta }=~^{\chi }\widehat{\lambda }~%
\mathbf{g}_{\alpha \beta }  \label{nheeq}
\end{equation}%
can be extracted from Ricci flows of solitonic hierarchies (such spaces with
nonhomogeneous effective constants $~^{\chi }\widehat{\lambda }(u)$ were
examined in details in Refs. \cite{3vsgg,3vncg}, for exact solutions in
gravity, and \cite{3nrf2,3nrfes,3nrft1,3nrft2}, for exact solutions and
applications in physics of the nonholonomic Ricci flow theory).

\begin{corollary}
Ricci flows of nonhomogeneous Einstein spaces, of signature $\left(
-+...+,-+...+\right) $ defined by (\ref{nheeq}) are can be extracted from
solitonic hierarchies satisfying the conditions of Theorem \ref{3mt} by
certain classes of constraints for
\begin{equation}
\overrightarrow{R}(\chi )=\left( n-1\right) ~^{\chi }\widehat{\lambda }=%
\overleftarrow{S}(\chi )=\left( m-1\right) ~^{\chi }\widehat{\lambda }
\label{fe1}
\end{equation}
solving the equations
\begin{equation}
\frac{\partial \widehat{\lambda }(\chi )}{\partial \chi }-\left[ \widehat{%
\lambda }(\chi )\right] ^{2}=\widehat{D}_{i}\widehat{D}^{i}~\widehat{\lambda
}(\chi )=\widehat{D}_{a}\widehat{D}^{a}~\widehat{\lambda }(\chi )
\label{fe2}
\end{equation}%
and evolution of N--adapted frames stated by formulas%
\begin{equation}
\frac{\partial \mathbf{e}_{~\underline{\alpha }}^{\alpha }}{\partial \chi }=%
\widehat{\lambda }(\chi )\mathbf{e}_{~\underline{\alpha }}^{\alpha }.
\label{fe3}
\end{equation}
\end{corollary}

\begin{proof}
The equations (\ref{fe2}) folow from equations (\ref{3hhfeq1}) and (\ref%
{3vhfeq1}), see also (\ref{3scev}), when hold true (\ref{fe1}). The equation
(\ref{fe3}) is a consequence of (\ref{3aeq5}). Similar proofs can be
provided for different signatures of metrics. $\square $
\end{proof}

The vacuum equations in Einstein gravity can be parametrized as solitonic
hierarchies with constant on $\chi $ N--anholonomic frames, see (\ref{fe3}).
The corresponding SG hierarchies are defined as solutions of the equations (%
\ref{3heq}) and (\ref{3sgeh}) and (\ref{3sgev}), with constant scalar
curvature  and frame coefficients.

\section{Conclusion}

In summary, we have considered a geometric formalism of encoding  general
(semi) Riemannian metrics and their lifts to tangent bundles  into families
of nonholonomic hierarchies of bi--Hamiltonian structures and  related
solitonic equations derived for curve flows on tangent spaces.  Towards this
ends, we have applied a programme of study that is based on prior works on
nonholonomic Ricci \cite{3nrf1,3nrf2} and curve flows \cite{3vsh,3avw}. The
premise of this methodology is that one can  derive solitonic hierarchies
for non--stretching curve flows on  constant curvature Riemannian spaces %
\cite{3chou1,3mbsw,3saw,3anc3}.  The validity of this approach was
substantiated by the encoding  into solitonic hierarchies of arbitrary
(semi) Riemannian and  Finsler--Lagrange metrics \cite{3vsh,3avw}.

Our analysis was completed by explicit constructions related to solitonic
encoding of Ricci  flow equations \cite{3ham1,3ham2,3per1,3per2,3per3}.
First of all, we note that it was not possible to  perform such
constructions working only with the Levi Civita  connection because, in
general, the curvature tensor for this  connection can not be parametrized
by constant coefficients.  The idea was to re--define equivalently the
geometric objects and basic Ricci flow and field equations for  other
classes of linear connections generated on vector/ tangent  bundles and/or
(semi) Riemannian manifolds enabled with  certain classes of preferred
frames with associated nonlinear  connection (N--connection) structure. In
particular, we  elaborated such lifts of (semi) Riemannian metrics to the
tangent bundle when in a canonical form there are defined  metric
structures, a class of N--connections and distinguished  (d) connections
when with respect to certain classes of  N--adapted frames it is possible to
get constant curvature  coefficients and zero torsion. Such geometric models
and their  lifts to different classes of fibred spaces  are related to
effective Lagrangians and, inversely, any  regular geometric mechanics can
be encoded into geometric  objects on nonholonomic Riemannian manifolds/
tangent bundles.

Secondly, we considered Ricci flows of geometric objects and  fundamental
evolution/field equations. We found that a  nonholonomically constrained
flow of (semi) Riemannian metrics  can result in generalized
Finsler--Lagrange configurations which  motivates both application of such
Finsler geometry methods for usual Riemannian  spaces if moving frames are
introduced into consideration and  a deep study of nonholonomic structures
with associated  N--connections soldering couples of Klein spaces and/or
constant curvature Riemannian spaces.

Third, we argued that the geometry Riemann and  Finsler--Lagrange spa\-ces
can be encoded into bi--Hamilton  structures and nonholonomic solitonic
equations and anticipated  that such curve flow -- solitonic hierarchies can
be constructed in a similar manner for exact solutions of
Einstein--Yang--Mills--Dirac equations, derived following the anholonomic
frame method, in noncommutative generalizations of gravity and geometry and
possible quantum models based on nonholonomic Lagrange--Fedosov manifolds.
In spite of the fact that there are a number of conceptual and technical
difficulties (such as the  physical meaning of the general N--connections,
additionally to the preferred  frame systems and nonlinear generalizations
of the usual linear connections,  explicit relations of the Ricci flows to
renormalization  group flows in quantum gravity models, cumbersome geometric
analysis calculus...)  the outcome of such approach is almost obvious that
we can encode  nonlinear fundamental  field/evolution equations in terms of
corresponding vector solitonic  equations, their hierarchies and
conservation laws.

Finally, we tried not only to encode some geometric information  about
metrics, connections and frames into solitons but also  formulated certain
criteria when from such solitonic equations  we can extract vacuum
gravitational spaces or certain more general  solutions for the Ricci flow
equations and related Einstein  spaces or Euler--Lagrange equations in
geometric mechanics. For  more general classes of solutions and extensions
to quantum  gravity, noncommutative geometry, these are desirable purposes
for  further investigations.

\vskip4pt \textbf{Acknowledgement: } The author is grateful to Prof. Barbara
Keyfitz for the possibility to communicate his work at the Fields Institute
Colloquium/ Seminar.

\appendix

\setcounter{equation}{0} \renewcommand{\theequation}
{A.\arabic{equation}} \setcounter{subsection}{0}
\renewcommand{\thesubsection}
{A.\arabic{subsection}}

\section{The Geometry of N--anholonomic Vector Bundles}

We denote by $\pi ^{\top }:TE\rightarrow TM$ the differential of map $\pi
:E\rightarrow M$ defined by fiber preserving morphisms of the tangent
bundles $TE$ and $TM.$ The kernel of $\pi ^{\top }$ is just the vertical
subspace $vE$ with a related inclusion mapping $i:vE\rightarrow TE.$

\begin{definition}
A nonlinear connection (N--connection) $\mathbf{N}$ on a vector bundle $%
\mathcal{E}$ \ is defined by the splitting on the left of an exact sequence
\begin{equation*}
0\rightarrow vE\overset{i}{\rightarrow }TE\rightarrow TE/vE\rightarrow 0,
\end{equation*}%
i. e. by a morphism of submanifolds $\mathbf{N:\ \ }TE\rightarrow vE$ such
that $\mathbf{N\circ }i$ is the unity in $vE.$
\end{definition}

In an equivalent form, we can say that a N--connection is defined by a
Whitney sum of conventional horizontal (h) subspace, $\left( hE\right) ,$
and vertical (v) subspace, $\left( vE\right) ,$
\begin{equation}
TE=hE\oplus vE.  \label{3whitney}
\end{equation}%
This sum defines a nonholonomic (equivalently, anholonomic, or nonitegrable)
distribution of horizontal and vertical subspaces on $TE\mathbf{.}$ Locally,
a N--connection is defined by its coefficients $N_{i}^{a}(u),$%
\begin{equation*}
\mathbf{N}=N_{i}^{a}(u)dx^{i}\otimes \frac{\partial }{\partial y^{a}}.
\end{equation*}%
The well known class of linear connections consists on a particular subclass
with the coefficients being linear on $y^{a},$ i.e., $N_{i}^{a}(u)=\Gamma
_{bj}^{a}(x)y^{b}.$\footnote{%
We use ''boldface'' symbols if it is necessary to emphasize that any space
and/or geometrical objects are provided/adapted to a\ N--connection
structure, or with the coefficients computed with respect to N--adapted
frames.}

\begin{remark}
A bundle space, or a a manifold, is called nonholonomic if it provided with
a nonholonomic distribution (see historical details and summary of results
in \cite{3bejf}). In particular case, when the nonholonomic distribution is
of type (\ref{3whitney}), such spaces are called N--anholonomic \cite{3vsgg}.
\end{remark}

Any N--connection $\mathbf{N}=\left\{ N_{i}^{a}(u)\right\} $ may be
characterized by a N--adapted frame (vielbein) structure $\mathbf{e}_{\nu
}=(e_{i},e_{a}),$ where
\begin{equation}
\mathbf{e}_{i}=\frac{\partial }{\partial x^{i}}-N_{i}^{a}(u)\frac{\partial }{%
\partial y^{a}}\mbox{ and
}e_{a}=\frac{\partial }{\partial y^{a}},  \label{3dder}
\end{equation}%
and the dual frame (coframe) structure $\mathbf{e}^{\mu }=(e^{i},\mathbf{e}%
^{a}),$ where
\begin{equation}
e^{i}=dx^{i}\mbox{ and }\mathbf{e}^{a}=dy^{a}+N_{i}^{a}(u)dx^{i}.
\label{3ddif}
\end{equation}

For any N--connection, we can introduce its N--connection curvature
\begin{equation*}
\mathbf{\Omega }=\frac{1}{2}\Omega _{ij}^{a}\ d^{i}\wedge d^{j}\otimes
\partial _{a},
\end{equation*}%
with the coefficients defined as the Neijenheuse tensor,%
\begin{equation}
\Omega _{ij}^{a}=\mathbf{e}_{[j}N_{i]}^{a}=\mathbf{e}_{j}N_{i}^{a}-\mathbf{e}%
_{i}N_{j}^{a}=\frac{\partial N_{i}^{a}}{\partial x^{j}}-\frac{\partial
N_{j}^{a}}{\partial x^{i}}+N_{i}^{b}\frac{\partial N_{j}^{a}}{\partial y^{b}}%
-N_{j}^{b}\frac{\partial N_{i}^{a}}{\partial y^{b}}.  \label{3ncurv}
\end{equation}

The vielbeins (\ref{3ddif}) satisfy the nonholonomy (equivalently,
anholonomy) relations
\begin{equation}
\lbrack \mathbf{e}_{\alpha },\mathbf{e}_{\beta }]=\mathbf{e}_{\alpha }%
\mathbf{e}_{\beta }-\mathbf{e}_{\beta }\mathbf{e}_{\alpha }=W_{\alpha \beta
}^{\gamma }\mathbf{e}_{\gamma }  \label{3anhrel}
\end{equation}%
with (antisymmetric) nontrivial anholonomy coefficients $W_{ia}^{b}=\partial
_{a}N_{i}^{b}$ and $W_{ji}^{a}=\Omega _{ij}^{a}.$

The geometric objects can be defined in a form adapted to a N--connection
structure, following decompositions being invariant under parallel
transports preserving the splitting (\ref{3whitney}). In this case we call
them to be distinguished (by the connection structure), i.e. d--objects. For
instance, a vector field $\mathbf{X}\in T\mathbf{V}$ \ is expressed
\begin{equation*}
\mathbf{X}=(hX,\ vX),\mbox{ \ or \ }\mathbf{X}=X^{\alpha }\mathbf{e}_{\alpha
}=X^{i}\mathbf{e}_{i}+X^{a}e_{a},
\end{equation*}%
where $hX=X^{i}\mathbf{e}_{i}$ and $vX=X^{a}e_{a}$ state, respectively, the
adapted to the N--connection structure horizontal (h) and vertical (v)
components of the vector (which following Refs. \cite{3ma1,3ma2} is called a
distinguished vector, in brief, d--vector). In a similar fashion, the
geometric objects on $\mathbf{V},$ for instance, tensors, spinors,
connections, ... are called respectively d--tensors, d--spinors,
d--connections if they are adapted to the N--connection splitting (\ref%
{3whitney}).

\begin{definition}
A distinguished connection (in brief, d--connection) $\mathbf{D}=(h\mathbf{D}%
,v\mathbf{D})$ is a linear connection preserving under parallel transports
the nonholonomic decomposition (\ref{3whitney}).
\end{definition}

The N--adapted components $\mathbf{\Gamma }_{\ \beta \gamma }^{\alpha }$ of
a d--connection $\mathbf{D}_{\alpha }=(\mathbf{e}_{\alpha }\rfloor \mathbf{D}%
)$ are defined by the equations
\begin{equation}
\mathbf{D}_{\alpha }\mathbf{e}_{\beta }=\mathbf{\Gamma }_{\ \alpha \beta
}^{\gamma }\mathbf{e}_{\gamma },\mbox{\ or \ }\mathbf{\Gamma }_{\ \alpha
\beta }^{\gamma }\left( u\right) =\left( \mathbf{D}_{\alpha }\mathbf{e}%
_{\beta }\right) \rfloor \mathbf{e}^{\gamma }.  \label{3dcon1}
\end{equation}%
The N--adapted splitting into h-- and v--covariant derivatives is stated by
\begin{equation*}
h\mathbf{D}=\{\mathbf{D}_{k}=\left( L_{jk}^{i},L_{bk\;}^{a}\right) \},%
\mbox{
and }\ v\mathbf{D}=\{\mathbf{D}_{c}=\left( C_{jk}^{i},C_{bc}^{a}\right) \},
\end{equation*}%
where, by definition, $L_{jk}^{i}=\left( \mathbf{D}_{k}\mathbf{e}_{j}\right)
\rfloor e^{i},$ $L_{bk}^{a}=\left( \mathbf{D}_{k}e_{b}\right) \rfloor
\mathbf{e}^{a},$ $C_{jc}^{i}=\left( \mathbf{D}_{c}\mathbf{e}_{j}\right)
\rfloor e^{i},$ $C_{bc}^{a}=\left( \mathbf{D}_{c}e_{b}\right) \rfloor
\mathbf{e}^{a}.$ The components $\mathbf{\Gamma }_{\ \alpha \beta }^{\gamma
}=\left( L_{jk}^{i},L_{bk}^{a},C_{jc}^{i},C_{bc}^{a}\right) $ completely
define a d--connection $\mathbf{D}$ on $\mathbf{E}.$

The simplest way to perform N--adapted computations is to use differential
forms. For instance, starting with the d--connection 1--form,
\begin{equation}
\mathbf{\Gamma }_{\ \beta }^{\alpha }=\mathbf{\Gamma }_{\ \beta \gamma
}^{\alpha }\mathbf{e}^{\gamma },  \label{3dconf}
\end{equation}%
with the coefficients defined with respect to N--elongated frames (\ref%
{3ddif}) and (\ref{3dder}), the torsion of a d--connection,
\begin{equation}
\mathcal{T}^{\alpha }\doteqdot \mathbf{De}^{\alpha }=d\mathbf{e}^{\alpha
}+\Gamma _{\ \beta }^{\alpha }\wedge \mathbf{e}^{\beta },  \label{3tors}
\end{equation}%
is characterized by (N--adapted) d--torsion components,
\begin{eqnarray}
T_{\ jk}^{i} &=&L_{\ jk}^{i}-L_{\ kj}^{i},\ T_{\ ja}^{i}=-T_{\ aj}^{i}=C_{\
ja}^{i},\ T_{\ ji}^{a}=\Omega _{\ ji}^{a},\   \notag \\
T_{\ bi}^{a} &=&-T_{\ ib}^{a}=\frac{\partial N_{i}^{a}}{\partial y^{b}}-L_{\
bi}^{a},\ T_{\ bc}^{a}=C_{\ bc}^{a}-C_{\ cb}^{a}.  \label{3dtors}
\end{eqnarray}%
For d--connection structures on $TM,$ we have to identify indices in the
form $i\leftrightarrows a,j\leftrightarrows b,...$ and the components of N--
and d--connections, for instance, $N_{i}^{a}\leftrightarrows N_{i}^{j}$ and $%
L_{\ jk}^{i}\leftrightarrows L_{\ bk}^{a},C_{\ ja}^{i}\leftrightarrows C_{\
ca}^{b}\leftrightarrows C_{\ jk}^{i}.$

\begin{definition}
A distinguished metric (in brief, d--metric) on a vector bundle $\mathbf{E}$
is a usual second rank metric tensor $\mathbf{g=}g\mathbf{\oplus _{N}}h,$
equivalently,
\begin{equation}
\mathbf{g}=\ g_{ij}(x,y)\ e^{i}\otimes e^{j}+\ h_{ab}(x,y)\ \mathbf{e}%
^{a}\otimes \mathbf{e}^{b},  \label{3m1}
\end{equation}%
adapted to the N--connection decomposition (\ref{3whitney}).
\end{definition}

With respect to a coordinate basis, the metric $\mathbf{g}$ (\ref{3m1}) can
be written in the form
\begin{equation}
\ \mathbf{g}=\underline{g}_{\alpha \beta }\left( u\right) du^{\alpha
}\otimes du^{\beta }  \label{3metr}
\end{equation}%
where%
\begin{equation}
\underline{g}_{\alpha \beta }=\left[
\begin{array}{cc}
g_{ij}+N_{i}^{a}N_{j}^{b}h_{ab} & N_{j}^{e}g_{ae} \\
N_{i}^{e}g_{be} & g_{ab}%
\end{array}%
\right] .  \label{3ansatz}
\end{equation}%
From the class of arbitrary d--connections $\mathbf{D}$ on $\mathbf{V,}$ one
distinguishes those which are metric compatible (metrical) satisfying the
condition%
\begin{equation}
\mathbf{Dg=0}  \label{3metcomp}
\end{equation}%
including all h- and v-projections $D_{j}g_{kl}=0,$ $D_{a}g_{kl}=0,$ $%
D_{j}h_{ab}=0,$ $D_{a}h_{bc}=0.$ For d--metric structures on $\mathbf{%
V\simeq }TM,$ with $g_{ij}=h_{ab},$ the conditions of vanishing
''nonmetricity'' (\ref{3metcomp}) transform into%
\begin{equation}
h\mathbf{D(}g\mathbf{)=}0\mbox{\  and\ }v\mathbf{D(}h\mathbf{)=}0,
\label{3metcompt}
\end{equation}%
i.e. $D_{j}g_{kl}=0$ and $D_{a}g_{kl}=0.$

There are two types of preferred linear connections uniquely determined by a
generic off--diagonal metric structure with $n+m$ splitting, see $\mathbf{g}%
=g\oplus _{N}h$:

\begin{enumerate}
\item The Levi Civita connection $\nabla =\{\Gamma _{\beta \gamma }^{\alpha
}\}$ is by definition torsionless, $~\ _{\shortmid }\mathcal{T}=0,$ and
satisfies the metric compatibility condition$,\nabla \mathbf{g}=0.$

\item The canonical d--connection $\widehat{\mathbf{\Gamma }}_{\ \alpha
\beta }^{\gamma }=\left( \widehat{L}_{jk}^{i},\widehat{L}_{bk}^{a},\widehat{C%
}_{jc}^{i},\widehat{C}_{bc}^{a}\right) $ is also metric compatible, i. e. $%
\widehat{\mathbf{D}}\mathbf{g}=0,$ but the torsion vanishes only on h-- and
v--subspaces, i.e. $\widehat{T}_{jk}^{i}=0$ and $\widehat{T}_{bc}^{a}=0,$
for certain nontrivial values of $\widehat{T}_{ja}^{i},\widehat{T}_{bi}^{a},%
\widehat{T}_{ji}^{a}.$ For simplicity, we omit hats on symbols and write,
for simplicity, $L_{jk}^{i}$ instead of $\widehat{L}_{jk}^{i},$ $T_{ja}^{i}$
instead of $\widehat{T}_{ja}^{i}$ and so on...but preserve the general
symbols $\widehat{\mathbf{D}}$ and $\widehat{\mathbf{\Gamma }}_{\ \alpha
\beta }^{\gamma }.$
\end{enumerate}

With respect to N--adapted frames (\ref{3dder}) and (\ref{3ddif}), we can
verify that the requested properties for $\widehat{\mathbf{D}}$ on $\mathbf{E%
}$ are satisfied if
\begin{eqnarray}
L_{jk}^{i} &=&\frac{1}{2}g^{ir}\left( \mathbf{e}_{k}g_{jr}+\mathbf{e}%
_{j}g_{kr}-\mathbf{e}_{r}g_{jk}\right) ,  \label{3candcon} \\
L_{bk}^{a} &=&e_{b}(N_{k}^{a})+\frac{1}{2}h^{ac}\left( \mathbf{e}%
_{k}h_{bc}-h_{dc}\ e_{b}N_{k}^{d}-h_{db}\ e_{c}N_{k}^{d}\right) ,  \notag \\
C_{jc}^{i} &=&\frac{1}{2}g^{ik}e_{c}g_{jk},\ C_{bc}^{a}=\frac{1}{2}%
h^{ad}\left( e_{c}h_{bd}+e_{c}h_{cd}-e_{d}h_{bc}\right) .  \notag
\end{eqnarray}%
For $\mathbf{E}=TM,$ the canonical d--connection $\ \mathbf{\tilde{D}}=(h%
\tilde{D},v\tilde{D})$ can be defined in torsionless form\footnote{%
i.e. it has the same coefficients as the Levi Civita connection with respect
to N--elongated bases (\ref{3dder}) and (\ref{3ddif})} with the coefficients
$\Gamma _{\ \beta \gamma }^{\alpha }=(L_{\ jk}^{i},L_{bc}^{a}),$
\begin{eqnarray}
L_{\ jk}^{i} &=&\frac{1}{2}g^{ih}(\mathbf{e}_{k}g_{jh}+\mathbf{e}_{j}g_{kh}-%
\mathbf{e}_{h}g_{jk}),  \label{3candcontm} \\
C_{\ bc}^{a} &=&\frac{1}{2}h^{ae}(e_{c}h_{be}+e_{b}h_{ce}-e_{e}h_{bc}).
\notag
\end{eqnarray}

The curvature of a d--connection $\mathbf{D,}$
\begin{equation}
\mathcal{R}_{~\beta }^{\alpha }\doteqdot \mathbf{D\Gamma }_{\ \beta
}^{\alpha }=d\mathbf{\Gamma }_{\ \beta }^{\alpha }-\mathbf{\Gamma }_{\ \beta
}^{\gamma }\wedge \mathbf{\Gamma }_{\ \gamma }^{\alpha },  \label{3curv}
\end{equation}%
splits into six types of N--adapted components with respect to (\ref{3dder})
and (\ref{3ddif}),
\begin{equation*}
\mathbf{R}_{~\beta \gamma \delta }^{\alpha }=\left(
R_{~hjk}^{i},R_{~bjk}^{a},P_{~hja}^{i},P_{~bja}^{c},S_{~jbc}^{i},S_{~bdc}^{a}\right) ,
\end{equation*}%
\begin{eqnarray}
R_{\ hjk}^{i} &=&\mathbf{e}_{k}L_{\ hj}^{i}-\mathbf{e}_{j}L_{\ hk}^{i}+L_{\
hj}^{m}L_{\ mk}^{i}-L_{\ hk}^{m}L_{\ mj}^{i}-C_{\ ha}^{i}\Omega _{\ kj}^{a},
\label{3dcurv} \\
R_{\ bjk}^{a} &=&\mathbf{e}_{k}L_{\ bj}^{a}-\mathbf{e}_{j}L_{\ bk}^{a}+L_{\
bj}^{c}L_{\ ck}^{a}-L_{\ bk}^{c}L_{\ cj}^{a}-C_{\ bc}^{a}\Omega _{\ kj}^{c},
\notag \\
P_{\ jka}^{i} &=&e_{a}L_{\ jk}^{i}-D_{k}C_{\ ja}^{i}+C_{\ jb}^{i}T_{\
ka}^{b},~P_{\ bka}^{c}=e_{a}L_{\ bk}^{c}-D_{k}C_{\ ba}^{c}+C_{\ bd}^{c}T_{\
ka}^{c},  \notag \\
S_{\ jbc}^{i} &=&e_{c}C_{\ jb}^{i}-e_{b}C_{\ jc}^{i}+C_{\ jb}^{h}C_{\
hc}^{i}-C_{\ jc}^{h}C_{\ hb}^{i},  \notag \\
S_{\ bcd}^{a} &=&e_{d}C_{\ bc}^{a}-e_{c}C_{\ bd}^{a}+C_{\ bc}^{e}C_{\
ed}^{a}-C_{\ bd}^{e}C_{\ ec}^{a}.  \notag
\end{eqnarray}

Contracting respectively the components, $\mathbf{R}_{\alpha \beta
}\doteqdot \mathbf{R}_{\ \alpha \beta \tau }^{\tau },$ one computes the h-
v--components of the Ricci d--tensor (there are four N--adapted components)
\begin{equation}
R_{ij}\doteqdot R_{\ ijk}^{k},\ \ R_{ia}\doteqdot -P_{\ ika}^{k},\
R_{ai}\doteqdot P_{\ aib}^{b},\ S_{ab}\doteqdot S_{\ abc}^{c}.
\label{3dricci}
\end{equation}%
The scalar curvature is defined by contracting the Ricci d--tensor with the
inverse metric $\mathbf{g}^{\alpha \beta },$
\begin{equation}
\overleftrightarrow{\mathbf{R}}\doteqdot \mathbf{g}^{\alpha \beta }\mathbf{R}%
_{\alpha \beta }=g^{ij}R_{ij}+h^{ab}S_{ab}=\overrightarrow{R}+\overleftarrow{%
S}.  \label{3sdccurv}
\end{equation}

If $\mathbf{E=}TM,$ there are only three classes of d--curvatures,%
\begin{eqnarray}
R_{\ hjk}^{i} &=&\mathbf{e}_{k}L_{\ hj}^{i}-\mathbf{e}_{j}L_{\ hk}^{i}+L_{\
hj}^{m}L_{\ mk}^{i}-L_{\ hk}^{m}L_{\ mj}^{i}-C_{\ ha}^{i}\Omega _{\ kj}^{a},
\label{3dcurvtb} \\
P_{\ jka}^{i} &=&e_{a}L_{\ jk}^{i}-\mathbf{D}_{k}C_{\ ja}^{i}+C_{\
jb}^{i}T_{\ ka}^{b},  \notag \\
S_{\ bcd}^{a} &=&e_{d}C_{\ bc}^{a}-e_{c}C_{\ bd}^{a}+C_{\ bc}^{e}C_{\
ed}^{a}-C_{\ bd}^{e}C_{\ ec}^{a},  \notag
\end{eqnarray}%
where all indices $a,b,...,i,j,...$ run the same values and, for instance, $%
C_{\ bc}^{e}\rightarrow $ $C_{\ jk}^{i},...$

\setcounter{equation}{0} \renewcommand{\theequation}
{B.\arabic{equation}} \setcounter{subsection}{0}
\renewcommand{\thesubsection}
{B.\arabic{subsection}}

\section{N--anholonomic Klein Spaces}

\label{apb}

There are Ricci flow families of two Hamiltonian variables given by the
principal normals $\;^{h}\nu $ and $\;^{v}\nu ,$ respectively, in the
horizontal and vertical subspaces, defined by the canonical d--connections $%
~^{\chi }\mathbf{D}=(h~^{\chi }\mathbf{D},v~^{\chi }\mathbf{D}),$ see
formulas (\ref{3curvframe}) and (\ref{3part02}),
\begin{equation*}
\;^{h}\nu (\chi )\doteqdot ~^{\chi }\mathbf{D}_{h\mathbf{X}}h\mathbf{X}=\nu
^{\widehat{i}}(\chi )~^{\chi }\mathbf{\mathbf{e}}_{\widehat{i}}%
\mbox{\ and \
}\;^{v}\nu (\chi )\doteqdot ~^{\chi }\mathbf{D}_{v\mathbf{X}}v\mathbf{X}=\nu
^{\widehat{a}}(\chi )~e_{\widehat{a}}.
\end{equation*}%
This normal d--vectors $\mathbf{v}(\chi )=(\;^{h}\nu (\chi ),$ $\;^{v}\nu
(\chi )),$ with components of type
\begin{equation*}
\mathbf{\nu }^{\alpha }(\chi )=(\nu ^{i}(\chi ),\;\nu ^{a}(\chi ))=(\nu
^{1}(\chi ),\nu ^{\widehat{i}}(\chi ),\nu ^{n+1}(\chi ),\nu ^{\widehat{a}%
}(\chi )),
\end{equation*}
are oriented in the tangent directions of curves $\gamma .$ There is also
the principal normal d--vectors $\mathbf{\varpi }(\chi )=(\;^{h}\varpi (\chi
),\;^{v}\varpi (\chi ))$ with components of type $\mathbf{\varpi }^{\alpha
}(\chi )=(\varpi ^{i}(\chi ),$ $\;\varpi ^{a}(\chi ))=(\varpi ^{1}(\chi
),\varpi ^{\widehat{i}}(\chi ),\varpi ^{n+1}(\chi ),\varpi ^{\widehat{a}%
}(\chi ))$ in the flow directions, with
\begin{equation*}
\;^{h}\varpi (\chi )\doteqdot ~^{\chi }\mathbf{D}_{h\mathbf{Y}}h\mathbf{X=}%
\varpi ^{\widehat{i}}(\chi )~^{\chi }\mathbf{\mathbf{e}}_{\widehat{i}%
},\;^{v}\varpi (\chi )\doteqdot ~^{\chi }\mathbf{D}_{v\mathbf{Y}}v\mathbf{X}%
=\varpi ^{\widehat{a}}(\chi )e_{\widehat{a}},
\end{equation*}%
representing a Hamiltonian d--covector field. We can consider that the
normal part of the flow d--vector
\begin{equation*}
\mathbf{h}_{\perp }(\chi )\doteqdot \mathbf{Y}_{\perp }(\chi )=h^{\widehat{i}%
}(\chi )~^{\chi }\mathbf{\mathbf{e}}_{\widehat{i}}+h^{\widehat{a}}(\chi )e_{%
\widehat{a}}
\end{equation*}%
represents a Hamiltonian d--vector field. For such configurations, we can
consider parallel N--adapted frames $~^{\chi }\mathbf{e}_{\alpha ^{\prime
}}=(~^{\chi }\mathbf{e}_{i^{\prime }},e_{a^{\prime }})$ when the
h--variables $\nu ^{\widehat{i^{\prime }}}(\chi ),$ $\varpi ^{\widehat{%
i^{\prime }}}(\chi ),h^{\widehat{i^{\prime }}}(\chi )$ are respectively
encoded in the top row of the horizontal canonical d--connection matrices $%
~^{\chi }\mathbf{\Gamma }_{h\mathbf{X\,}i^{\prime }}^{\qquad j^{\prime }}$
and $~^{\chi }\mathbf{\Gamma }_{h\mathbf{Y\,}i^{\prime }}^{\qquad j^{\prime
}}$ and in the row matrix $\left( ~^{\chi }\mathbf{e}_{\mathbf{Y}%
}^{i^{\prime }}\right) _{\perp }\doteqdot ~^{\chi }\mathbf{e}_{\mathbf{Y}%
}^{i^{\prime }}-~^{\chi }g_{\parallel }\;~^{\chi }\mathbf{e}_{\mathbf{X}%
}^{i^{\prime }}$ where $~^{\chi }g_{\parallel }\doteqdot ~^{\chi }g(h\mathbf{%
Y,}h\mathbf{X})$ is the tangential h--part of the flow d--vector. A similar
encoding holds for v--variables $\nu ^{\widehat{a^{\prime }}}(\chi ),\varpi
^{\widehat{a^{\prime }}}(\chi ),h^{\widehat{a^{\prime }}}(\chi )$ in the top
row of the vertical canonical d--connection matrices \ $~^{\chi }\mathbf{%
\Gamma }_{v\mathbf{X\,}a^{\prime }}^{\qquad b^{\prime }}$ and $~^{\chi }%
\mathbf{\Gamma }_{v\mathbf{Y\,}a^{\prime }}^{\qquad b^{\prime }}$ and in the
row matrix $\left( ~^{\chi }\mathbf{e}_{\mathbf{Y}}^{a^{\prime }}\right)
_{\perp }\doteqdot ~^{\chi }\mathbf{e}_{\mathbf{Y}}^{a^{\prime }}-~^{\chi
}h_{\parallel }\;~^{\chi }\mathbf{e}_{\mathbf{X}}^{a^{\prime }}$ where $%
~^{\chi }h_{\parallel }\doteqdot ~^{\chi }h(v\mathbf{Y,}v\mathbf{X})$ is the
tangential v--part of the flow d--vector. In a compact form of notations, we
shall write $\mathbf{v}^{\alpha ^{\prime }}(\chi )$ and $\mathbf{\varpi }%
^{\alpha ^{\prime }}(\chi )$ where the primed small Greek indices $\alpha
^{\prime },\beta ^{\prime },...$ will denote both N--adapted and then
orthonormalized components of geometric objects (like d--vectors,
d--covectors, d--tensors, d--groups, d--algebras, d--matrices) admitting
further decompositions into h-- and v--components defined as nonintegrable
distributions of such objects.

With respect to N--adapted orthonormalized frames, the geometry of
N--anholonomic manifolds is defined algebraically, on their tangent bundles,
by couples of horizontal and vertical Klein geometries considered in \cite%
{3sharpe} and for bi--Hamiltonian soliton constructions in \cite{3anc2}. The
N--connection structure induces a N--anholonomic Klein space stated by two
left--invariant $h\mathfrak{g}$-- and $v\mathfrak{g}$--valued Maurer--Cartan
form on the Lie d--group $\mathbf{G}=(h\mathbf{G},v\mathbf{G})$ is
identified with the zero--curvature canonical d--connection 1--form $\;^{%
\mathbf{G}}\mathbf{\Gamma }=\{\;^{\mathbf{G}}\mathbf{\Gamma }_{\ \beta
^{\prime }}^{\alpha ^{\prime }}\},$ where
\begin{equation*}
\;^{\mathbf{G}}\mathbf{\Gamma }_{\ \beta ^{\prime }}^{\alpha ^{\prime }}=\;^{%
\mathbf{G}}\mathbf{\Gamma }_{\ \beta ^{\prime }\gamma ^{\prime }}^{\alpha
^{\prime }}\mathbf{e}^{\gamma ^{\prime }}=\;^{h\mathbf{G}}L_{\;j^{\prime
}k^{\prime }}^{i^{\prime }}\mathbf{e}^{k^{\prime }}+\;^{v\mathbf{G}%
}C_{\;j^{\prime }k^{\prime }}^{i^{\prime }}e^{k^{\prime }}.
\end{equation*}%
For trivial N--connection structure in vector bundles with the base and
typical fiber spaces being symmetric Riemannian spaces, we can consider that
$\;^{h\mathbf{G}}L_{\;j^{\prime }k^{\prime }}^{i^{\prime }}$ and $\;^{v%
\mathbf{G}}C_{\;j^{\prime }k^{\prime }}^{i^{\prime }}$ are the coefficients
of the Cartan connections $\;^{h\mathbf{G}}L$ and $\;^{v\mathbf{G}}C,$
respectively for the $h\mathbf{G}$ and $v\mathbf{G,}$ both with vanishing
curvatures, i.e. with
\begin{equation*}
d\;^{\mathbf{G}}\mathbf{\Gamma +}\frac{1}{2}\mathbf{[\;^{\mathbf{G}}\mathbf{%
\Gamma ,}\;^{\mathbf{G}}\mathbf{\Gamma }]=0}
\end{equation*}%
and h-- and v--components, $d\;^{h\mathbf{G}}\mathbf{L}+\frac{1}{2}\mathbf{%
[\;^{h\mathbf{G}}L,\;^{h\mathbf{G}}L]}=0$ and $d\;^{v\mathbf{G}}\mathbf{C}+%
\frac{1}{2}[\;^{v\mathbf{G}}\mathbf{C},$ $\;^{v\mathbf{G}}\mathbf{C}]=0,$
where $d$ denotes the total derivatives on the d--group manifold $\mathbf{G}%
=h\mathbf{G}\oplus v\mathbf{G}$ or their restrictions on $h\mathbf{G}$ or $v%
\mathbf{G.}$ We can consider that $\;^{\mathbf{G}}\mathbf{\Gamma }$ defines
the so--called Cartan d--connection for nonintegrable N--connection
structures, see details and supersymmetric/ noncommutative developments in %
\cite{3vncg,3vsgg}.

Through the Lie d--algebra decompositions $\mathfrak{g}=h\mathfrak{g}\oplus v%
\mathfrak{g,}$ for the horizontal splitting: $h\mathfrak{g}=\mathfrak{so}%
(n)\oplus h\mathfrak{p,}$ when $[h\mathfrak{p},h\mathfrak{p}]\subset
\mathfrak{so}(n)$ and $[\mathfrak{so}(n),h\mathfrak{p}]\subset h\mathfrak{p;}
$ for the vertical splitting $v\mathfrak{g}=\mathfrak{so}(m)\oplus v%
\mathfrak{p,}$ when $[v\mathfrak{p},v\mathfrak{p}]\subset \mathfrak{so}(m)$
and $[\mathfrak{so}(m),v\mathfrak{p}]\subset v\mathfrak{p,}$ the Cartan
d--connection determines an N--anholonomic Riemannian structure on the
nonholonomic bundle $\mathbf{\mathring{E}}=[hG=SO(n+1),$ $%
vG=SO(m+1),\;N_{i}^{e}].$ For $n=m,$ and canonical d--objects
(N--connection, d--metric, d--connection, ...) derived from (\ref{3m1b}), or
any N--anholonomic space with constant d--curvatures, the Cartan
d--connection transform just in the canonical d--connection (\ref{3candcontm}%
). It is possible to consider a quotient space with distinguished structure
group $\mathbf{V}_{\mathbf{N}}=\mathbf{G}/SO(n)\oplus $ $SO(m)$ regarding $%
\mathbf{G}$ as a principal $\left( SO(n)\oplus SO(m)\right) $--bundle over $%
\mathbf{\mathring{E}},$ which is a N--anholonomic bundle. In this case, we
can always fix a local section of this bundle and pull--back $\;^{\mathbf{G}}%
\mathbf{\Gamma }$ to give a $\left( h\mathfrak{g}\oplus v\mathfrak{g}\right)
$--valued 1--form $^{\mathfrak{g}}\mathbf{\Gamma }$ in a point $u\in \mathbf{%
\mathring{E}}.$ Any change of local sections define $SO(n)\oplus $ $SO(m)$
gauge transforms of the canonical d--connection $^{\mathfrak{g}}\mathbf{%
\Gamma ,}$ all preserving the nonholonomic decomposition (\ref{3whitney}).

There are involutive automorphisms $h\sigma =\pm 1$ and $v\sigma =\pm 1,$
respectively, of $h\mathfrak{g}$ and $v\mathfrak{g,}$ defined that $%
\mathfrak{so}(n)$ (or $\mathfrak{so}(m)$) is eigenspace $h\sigma =+1$ (or $%
v\sigma =+1)$ and $h\mathfrak{p}$ (or $v\mathfrak{p}$) is eigenspace $%
h\sigma =-1$ (or $v\sigma =-1).$ Taking into account the existing
eigenspaces, when the symmetric parts $\mathbf{\Gamma \doteqdot }\frac{1}{2}%
\left( ^{\mathfrak{g}}\mathbf{\Gamma +}\sigma \left( ^{\mathfrak{g}}\mathbf{%
\Gamma }\right) \right),$ with respective h- and v--splitting, $\mathbf{%
L\doteqdot }\frac{1}{2}\left( ^{h\mathfrak{g}}\mathbf{L+}h\sigma \left( ^{h%
\mathfrak{g}}\mathbf{L}\right) \right) $ and $\mathbf{C\doteqdot }\frac{1}{2}%
(^{v\mathfrak{g}}\mathbf{C}+h\sigma (^{v\mathfrak{g}}\mathbf{C})),$ defines
a $\left( \mathfrak{so}(n)\oplus \mathfrak{so}(m)\right) $--valued
d--connection 1--form, we can construct N--adapted decompositions.

The antisymmetric part $\mathbf{e\doteqdot }\frac{1}{2}\left( ^{\mathfrak{g}}%
\mathbf{\Gamma -}\sigma \left( ^{\mathfrak{g}}\mathbf{\Gamma }\right)
\right),$ with h- and v--splitting, $h\mathbf{e\doteqdot }$ $\frac{1}{2}%
\left( ^{h\mathfrak{g}}\mathbf{L-}h\sigma \left( ^{h\mathfrak{g}}\mathbf{L}%
\right) \right) $ and $v\mathbf{e\doteqdot }\frac{1}{2}(^{v\mathfrak{g}}%
\mathbf{C}-h\sigma (^{v\mathfrak{g}}\mathbf{C})),$ defines a $\left( h%
\mathfrak{p}\oplus v\mathfrak{p}\right) $--valued N--adapted coframe for the
Cartan--Killing inner product $<\cdot ,\cdot >_{\mathfrak{p}}$ on $T_{u}%
\mathbf{G}\simeq h\mathfrak{g}\oplus v\mathfrak{g}$ restricted to $T_{u}%
\mathbf{V}_{\mathbf{N}}\simeq \mathfrak{p.}$ This inner product,
distinguished into h- and v--components, provides a d--metric structure of
type $\mathbf{g}=[g,h]$ (\ref{3m1}), where $g=<h\mathbf{e\otimes }h\mathbf{e}%
>_{h\mathfrak{p}}$ and $h=<v\mathbf{e\otimes }v\mathbf{e}>_{v\mathfrak{p}}$
on $\mathbf{V}_{\mathbf{N}}=\mathbf{G}/SO(n)\oplus $ $SO(m).$

We generate a $\mathbf{G(}=h\mathbf{G}\oplus v\mathbf{G)}$--invariant
d--derivatives $\mathbf{D}$ whose restriction to the tangent space $T\mathbf{%
V}_{\mathbf{N}}$ for any N--anholonomic curve flow $\gamma (\tau ,\mathbf{l,}%
\chi )$ in $\mathbf{V}_{\mathbf{N}}=\mathbf{G}/SO(n)\oplus $ $SO(m)$ is
defined via%
\begin{equation}
~^{\chi }\mathbf{D}_{\mathbf{X}}~^{\chi }\mathbf{e=}\left[ ~^{\chi }\mathbf{e%
},\gamma _{\mathbf{l}}\rfloor ~^{\chi }\mathbf{\Gamma }\right]
\mbox{\ and \
}~^{\chi }\mathbf{D}_{\mathbf{Y}}~^{\chi }\mathbf{e=}\left[ ~^{\chi }\mathbf{%
e},\gamma _{\mathbf{\tau }}\rfloor ~^{\chi }\mathbf{\Gamma }\right] ,
\label{3aux33}
\end{equation}%
admitting further h- and v--decompositions. The derivatives $~^{\chi }%
\mathbf{D}_{\mathbf{X}}$ and $~^{\chi }\mathbf{D}_{\mathbf{Y}}$ in (\ref%
{3aux33}) are equivalent to those considered in (\ref{3part01}) and obey the
Cartan structure equations (\ref{3mtors}) and (\ref{3mcurv}). For the
canonical d--connections, a large class of N--anholonomic spaces of
dimension $n=m,$ the d--torsions are zero and the d--curvatures are with
constant coefficients.

Let $~^{\chi }\mathbf{e}^{\alpha ^{\prime }}=(e^{i^{\prime }},~^{\chi }%
\mathbf{e}^{a^{\prime }})$ be a family of N--adapted orthonormalized
coframes identified with the $\left( h\mathfrak{p}\oplus v\mathfrak{p}%
\right) $--valued coframe $\mathbf{e}$ in a fixed orthonormal basis for $%
\mathfrak{p=}h\mathfrak{p}\oplus v\mathfrak{p\subset }h\mathfrak{g}\oplus v%
\mathfrak{g.}$ Considering the kernel/ cokernel of Lie algebra
multiplications in the h- and v--subspaces, respectively, $\left[ \mathbf{e}%
_{h\mathbf{X}},\cdot \right] _{h\mathfrak{g}}$ and $\left[ \mathbf{e}_{v%
\mathbf{X}},\cdot \right] _{v\mathfrak{g}},$ we can decompose the coframes
into parallel and perpendicular parts with respect to $\mathbf{e}_{\mathbf{X}%
}.$ We write
\begin{equation*}
~^{\chi }\mathbf{e=(e}_{C}=h\mathbf{e}_{C}+v\mathbf{e}_{C},\mathbf{e}%
_{C^{\perp }}=h\mathbf{e}_{C^{\perp }}+v\mathbf{e}_{C^{\perp }}\mathbf{),}
\end{equation*}%
for $\mathfrak{p(}=h\mathfrak{p}\oplus v\mathfrak{p)}$--valued mutually
orthogonal d--vectors $\mathbf{e}_{C}$ \ and $\mathbf{e}_{C^{\perp }},$ when
there are satisfied the conditions $\left[ \mathbf{e}_{\mathbf{X}},\mathbf{e}%
_{C}\right] _{\mathfrak{g}}=0$ but $\left[ \mathbf{e}_{\mathbf{X}},\mathbf{e}%
_{C^{\perp }}\right] _{\mathfrak{g}}\neq 0;$ such conditions can be stated
in h- and v--component form, respectively, $\left[ h\mathbf{e}_{\mathbf{X}},h%
\mathbf{e}_{C}\right] _{h\mathfrak{g}}=0,$ $\left[ h\mathbf{e}_{\mathbf{X}},h%
\mathbf{e}_{C^{\perp }}\right] _{h\mathfrak{g}}\neq 0$ and $\left[ v\mathbf{e%
}_{\mathbf{X}},v\mathbf{e}_{C}\right] _{v\mathfrak{g}}=0,$ $\left[ v\mathbf{e%
}_{\mathbf{X}},v\mathbf{e}_{C^{\perp }}\right] _{v\mathfrak{g}}\neq 0.$ One
holds also the algebraic decompositions
\begin{eqnarray*}
T_{u}\mathbf{V}_{\mathbf{N}} &\simeq &\mathfrak{p=}h\mathfrak{p}\oplus v%
\mathfrak{p}=\mathfrak{g=}h\mathfrak{g}\oplus v\mathfrak{g}/\mathfrak{so}%
(n)\oplus \mathfrak{so}(m), \\
\mathfrak{p}\mathfrak{=p} &&_{C}\oplus \mathfrak{p}_{C^{\perp }}=\left( h%
\mathfrak{p}_{C}\oplus v\mathfrak{p}_{C}\right) \oplus \left( h\mathfrak{p}%
_{C^{\perp }}\oplus v\mathfrak{p}_{C^{\perp }}\right) ,
\end{eqnarray*}%
with $\mathfrak{p}_{\parallel }\subseteq \mathfrak{p}_{C}$ and $\mathfrak{p}%
_{C^{\perp }}\subseteq \mathfrak{p}_{\perp },$ where $\left[ \mathfrak{p}%
_{\parallel },\mathfrak{p}_{C}\right] =0,$ $<\mathfrak{p}_{C^{\perp }},%
\mathfrak{p}_{C}>=0,$ but $\left[ \mathfrak{p}_{\parallel },\mathfrak{p}%
_{C^{\perp }}\right] \neq 0$ (i.e. $\mathfrak{p}_{C}$ is the centralizer of $%
\mathbf{e}_{\mathbf{X}}$ in $\mathfrak{p=}h\mathfrak{p}\oplus v\mathfrak{%
p\subset }h\mathfrak{g}\oplus v\mathfrak{g);}$ in h- \ and v--components,
one have $h\mathfrak{p}_{\parallel }\subseteq h\mathfrak{p}_{C}$ and $h%
\mathfrak{p}_{C^{\perp }}\subseteq h\mathfrak{p}_{\perp },$ where $\left[ h%
\mathfrak{p}_{\parallel },h\mathfrak{p}_{C}\right] =0,$ $<h\mathfrak{p}%
_{C^{\perp }},h\mathfrak{p}_{C}>=0,$ but $\left[ h\mathfrak{p}_{\parallel },h%
\mathfrak{p}_{C^{\perp }}\right] \neq 0$ (i.e. $h\mathfrak{p}_{C}$ is the
centralizer of $\mathbf{e}_{h\mathbf{X}}$ in $h\mathfrak{p\subset }h%
\mathfrak{g)}$ and $v\mathfrak{p}_{\parallel }\subseteq v\mathfrak{p}_{C}$
and $v\mathfrak{p}_{C^{\perp }}\subseteq v\mathfrak{p}_{\perp },$ where $%
\left[ v\mathfrak{p}_{\parallel },v\mathfrak{p}_{C}\right] =0,$ $<v\mathfrak{%
p}_{C^{\perp }},v\mathfrak{p}_{C}>=0,$ but $\left[ v\mathfrak{p}_{\parallel
},v\mathfrak{p}_{C^{\perp }}\right] \neq 0$ (i.e. $v\mathfrak{p}_{C}$ is the
centralizer of $\mathbf{e}_{v\mathbf{X}}$ in $v\mathfrak{p\subset }v%
\mathfrak{g).}$ Using the canonical d--connection derivative $\mathbf{D}_{%
\mathbf{X}}$ of a d--covector perpendicular (or parallel) to $\mathbf{e}_{%
\mathbf{X}},$ we get a new d--vector which is parallel (or perpendicular) to
$\mathbf{e}_{\mathbf{X}},$ i.e. $\mathbf{D}_{\mathbf{X}}\mathbf{e}_{C}\in
\mathfrak{p}_{C^{\perp }}$ (or $\mathbf{D}_{\mathbf{X}}\mathbf{e}_{C^{\perp
}}\in \mathfrak{p}_{C});$ in h- \ and v--components such formulas are
written $\mathbf{D}_{h\mathbf{X}}h\mathbf{e}_{C}\in h\mathfrak{p}_{C^{\perp
}}$ (or $\mathbf{D}_{h\mathbf{X}}h\mathbf{e}_{C^{\perp }}\in h\mathfrak{p}%
_{C})$ and $\mathbf{D}_{v\mathbf{X}}v\mathbf{e}_{C}\in v\mathfrak{p}%
_{C^{\perp }}$ (or $\mathbf{D}_{v\mathbf{X}}v\mathbf{e}_{C^{\perp }}\in v%
\mathfrak{p}_{C}).$ All such d--algebraic relations can be written in
N--anholonomic manifolds and canonical d--connection settings, for instance,
using certain relations of type
\begin{equation*}
~^{\chi }\mathbf{D}_{\mathbf{X}}(~^{\chi }\mathbf{e}^{\alpha ^{\prime
}})_{C}=~^{\chi }\mathbf{v}_{~\beta ^{\prime }}^{\alpha ^{\prime }}(~^{\chi }%
\mathbf{e}^{\beta ^{\prime }})_{C^{\perp }}\mbox{ \ and \ }~^{\chi }\mathbf{D%
}_{\mathbf{X}}(~^{\chi }\mathbf{e}^{\alpha ^{\prime }})_{C^{\perp
}}=-~^{\chi }\mathbf{v}_{~\beta ^{\prime }}^{\alpha ^{\prime }}(~^{\chi }%
\mathbf{e}^{\beta ^{\prime }})_{C},
\end{equation*}%
for some antisymmetric d--tensors $~^{\chi }\mathbf{v}^{\alpha ^{\prime
}\beta ^{\prime }}=-~^{\chi }\mathbf{v}^{\beta ^{\prime }\alpha ^{\prime }}.$
We get a N--adapted $\left( SO(n)\oplus SO(m)\right) $--parallel frame
defining a generalization of the concept of Riemannian parallel frame on
N--adapted manifolds whenever $\mathfrak{p}_{C}$ is larger than $\mathfrak{p}%
_{\parallel }.$ Substituting $~^{\chi }\mathbf{e}^{\alpha ^{\prime
}}=(e^{i^{\prime }},~^{\chi }\mathbf{e}^{a^{\prime }})$ into the last
formulas and considering h- and v--components, we define $SO(n)$--parallel
and $SO(m)$--parallel frames (for simplicity we omit these formulas when the
Greek small letter indices are split into Latin small letter h- and
v--indices).

The final conclusion of this section is that the Cartan structure equations
on hypersurfaces swept out by nonholonomic curve flows on N--anholonomic
spaces with constant matrix curvature for the canonical d--connection
geometrically encode \ two $O(n-1)$-- and $O(m-1)$--invariant, respectively,
horizontal and vertical bi--Hamiltonian operators. This holds true if the
distinguished by N--connection freedom of the d--group action $SO(n)\oplus
SO(m)$ on $~^{\chi }\mathbf{e}$ and $~^{\chi }\mathbf{\Gamma }$ is used to
fix them to be a N--adapted parallel coframe and its associated canonical
d--connection 1--form is related to the canonical covariant derivative on
N--anholonomic manifolds.

\setcounter{equation}{0} \renewcommand{\theequation}
{C.\arabic{equation}} \setcounter{subsection}{0}
\renewcommand{\thesubsection}
{C.\arabic{subsection}}

\section{Proof of the Main Theorem}

\label{3ssp} We provide a proof of Theorem \ref{3mt} for the horizontal
Ricci and curve flows (similar results were published in \cite{3avw} and %
\cite{3anc2}, respectively, for Lagrange--Finsler and symmetric Riemannian
spaces). The vertical constructions are similar but with respective changing
of h-- variables / objects into v- variables/ objects.

One obtains a vector mKdV equation up to a convective term, which can be
absorbed by redefinition of coordinates, defining the +1 flow for $h%
\overrightarrow{\mathbf{e}}_{\perp }(\chi )=~^{\chi }\overrightarrow{v}_{%
\mathbf{l}},$%
\begin{equation*}
~^{\chi }\overrightarrow{v}_{\tau }=~^{\chi }\overrightarrow{v}_{3\mathbf{l}%
}+\frac{3}{2}|\overrightarrow{v}(\chi )|^{2}-\overrightarrow{R}(\chi
)~^{\chi }\overrightarrow{v}_{\mathbf{l}},
\end{equation*}%
when the $+(k+1)$ flow gives a vector mKdV equation of higher order $3+2k$
on $\overrightarrow{v}$ and there is a $0$ h--flow $\overrightarrow{v}_{\tau
}=\overrightarrow{v}_{\mathbf{l}}$ arising from $h\overrightarrow{\mathbf{e}}%
_{\perp }=0$ and $h\overrightarrow{\mathbf{e}}_{\parallel }=1$ belonging
outside the hierarchy generated by $h\mathfrak{R}(\chi ).$ Such flows
correspond to N--adapted horizontal motions of the curve $~^{\chi }\gamma
(\tau ,\mathbf{l})=h\gamma (\tau ,\mathbf{l},\chi )+v\gamma (\tau ,\mathbf{l}%
,\chi ),$ given by
\begin{equation*}
\left( h\gamma \right) _{\tau }(\chi )=f\left( \left( h\gamma \right) _{h%
\mathbf{X}}(\chi ),~^{\chi }\mathbf{D}_{h\mathbf{X}}\left( h\gamma \right)
_{h\mathbf{X}}(\chi ),~^{\chi }\mathbf{D}_{h\mathbf{X}}^{2}\left( h\gamma
\right) _{h\mathbf{X}}(\chi ),...\right)
\end{equation*}%
subject to the non--stretching condition $|\left( h\gamma \right) _{h\mathbf{%
X}}(\chi )|_{h\mathbf{g}}=1,$ when the equation of motion is to be derived
from the identifications
\begin{equation*}
\left( h\gamma \right) _{\tau }(\chi )\longleftrightarrow \mathbf{e}_{h%
\mathbf{Y}}(\chi ),~^{\chi }\mathbf{D}_{h\mathbf{X}}\left( h\gamma \right)
_{h\mathbf{X}}(\chi )\longleftrightarrow ~^{\chi }\mathcal{D}_{h\mathbf{X}%
}~^{\chi }\mathbf{e}_{h\mathbf{X}}=\left[ ~^{\chi }\mathbf{L}_{h\mathbf{X}%
},~^{\chi }\mathbf{e}_{h\mathbf{X}}\right]
\end{equation*}%
and so on, which maps the constructions from the tangent space of the curve
to the space $h\mathfrak{p}.$ For such identifications, we have
\begin{eqnarray*}
\left[ ~^{\chi }\mathbf{L}_{h\mathbf{X}},~^{\chi }\mathbf{e}_{h\mathbf{X}}%
\right] &=&-\left[
\begin{array}{cc}
0 & \left( 0,~^{\chi }\overrightarrow{v}\right) \\
-\left( 0,~^{\chi }\overrightarrow{v}\right) ^{T} & h\mathbf{0}%
\end{array}%
\right] \in h\mathfrak{p}, \\
\left[ ~^{\chi }\mathbf{L}_{h\mathbf{X}},\left[ ~^{\chi }\mathbf{L}_{h%
\mathbf{X}},~^{\chi }\mathbf{e}_{h\mathbf{X}}\right] \right] &=&-\left[
\begin{array}{cc}
0 & \left( |~^{\chi }\overrightarrow{v}|^{2},\overrightarrow{0}\right) \\
-\left( |~^{\chi }\overrightarrow{v}|^{2},\overrightarrow{0}\right) ^{T} & h%
\mathbf{0}%
\end{array}%
\right]
\end{eqnarray*}%
and so on, see similar calculus in (\ref{3aux41}). At the next step, stating
for the +1 h--flow
\begin{equation*}
h\overrightarrow{\mathbf{e}}_{\perp }(\chi )=~^{\chi }\overrightarrow{v}_{%
\mathbf{l}}\mbox{
and }h\overrightarrow{\mathbf{e}}_{\parallel }(\chi )=-~^{\chi }\mathbf{D}_{h%
\mathbf{X}}^{-1}\left( ~^{\chi }\overrightarrow{v}\cdot ~^{\chi }%
\overrightarrow{v}_{\mathbf{l}}\right) =-\frac{1}{2}|\overrightarrow{v}(\chi
)|^{2},
\end{equation*}%
we compute
\begin{eqnarray*}
~^{\chi }\mathbf{e}_{h\mathbf{Y}} &=&\left[
\begin{array}{cc}
0 & \left( h\mathbf{e}_{\parallel },h\overrightarrow{\mathbf{e}}_{\perp
}\right) (\chi ) \\
-\left( h\mathbf{e}_{\parallel },h\overrightarrow{\mathbf{e}}_{\perp
}\right) ^{T}(\chi ) & h\mathbf{0}%
\end{array}%
\right] \\
&=&-\frac{1}{2}|\overrightarrow{v}(\chi )|^{2}\left[
\begin{array}{cc}
0 & \left( 1,\overrightarrow{\mathbf{0}}\right) \\
-\left( 0,\overrightarrow{\mathbf{0}}\right) ^{T} & h\mathbf{0}%
\end{array}%
\right] \\
&&+\left[
\begin{array}{cc}
0 & \left( 0,~^{\chi }\overrightarrow{v}_{h\mathbf{X}}\right) \\
-\left( 0,~^{\chi }\overrightarrow{v}_{h\mathbf{X}}\right) ^{T} & h\mathbf{0}%
\end{array}%
\right] \\
&=&~^{\chi }\mathbf{D}_{h\mathbf{X}}\left[ ~^{\chi }\mathbf{L}_{h\mathbf{X}%
},~^{\chi }\mathbf{e}_{h\mathbf{X}}\right] +\frac{1}{2}\left[ ~^{\chi }%
\mathbf{L}_{h\mathbf{X}},\left[ ~^{\chi }\mathbf{L}_{h\mathbf{X}},~^{\chi }%
\mathbf{e}_{h\mathbf{X}}\right] \right] \\
&=&-~^{\chi }\mathcal{D}_{h\mathbf{X}}\left[ ~^{\chi }\mathbf{L}_{h\mathbf{X}%
},~^{\chi }\mathbf{e}_{h\mathbf{X}}\right] -\frac{3}{2}|\overrightarrow{v}%
(\chi )|^{2}~^{\chi }\mathbf{e}_{h\mathbf{X}}.
\end{eqnarray*}%
Following above presented identifications related to the first and second
terms, when
\begin{eqnarray*}
|\overrightarrow{v}(\chi )|^{2} &=&<\left[ ~^{\chi }\mathbf{L}_{h\mathbf{X}%
},~^{\chi }\mathbf{e}_{h\mathbf{X}}\right] ,\left[ ~^{\chi }\mathbf{L}_{h%
\mathbf{X}},~^{\chi }\mathbf{e}_{h\mathbf{X}}\right] >_{h\mathfrak{p}} \\
&\longleftrightarrow &h\mathbf{g}\left( ~^{\chi }\mathbf{D}_{h\mathbf{X}%
}\left( h\gamma \right) _{h\mathbf{X}}(\chi ),~^{\chi }\mathbf{D}_{h\mathbf{X%
}}\left( h\gamma \right) _{h\mathbf{X}}(\chi ),\chi \right) \\
&=&\left| ~^{\chi }\mathbf{D}_{h\mathbf{X}}\left( h\gamma \right) _{h\mathbf{%
X}}(\chi )\right| _{h\mathbf{g}}^{2},
\end{eqnarray*}%
we can identify $~^{\chi }\mathcal{D}_{h\mathbf{X}}\left[ ~^{\chi }\mathbf{L}%
_{h\mathbf{X}},~^{\chi }\mathbf{e}_{h\mathbf{X}}\right] $ to $~^{\chi }%
\mathbf{D}_{h\mathbf{X}}^{2}\left( h\gamma \right) _{h\mathbf{X}}(\chi )$
and write
\begin{equation*}
-~^{\chi }\mathbf{e}_{h\mathbf{Y}}\longleftrightarrow ~^{\chi }\mathbf{D}_{h%
\mathbf{X}}^{2}\left( h\gamma \right) _{h\mathbf{X}}(\chi )+\frac{3}{2}%
\left| ~^{\chi }\mathbf{D}_{h\mathbf{X}}\left( h\gamma \right) _{h\mathbf{X}%
}(\chi )\right| _{h\mathbf{g}}^{2}~\left( h\gamma \right) _{h\mathbf{X}%
}(\chi )
\end{equation*}%
which is just the first equation (\ref{31map}) in the Theorem \ref{3mt}
defining a family of non--stretching mKdV map h--equations induced by the
h--part of the family of canonical d--connections.

Using the adjoint representation $ad\left( \cdot \right) $ acting in the Lie
algebra $h\mathfrak{g}=h\mathfrak{p}\oplus \mathfrak{so}(n),$ with
\begin{equation*}
ad\left( \left[ ~^{\chi }\mathbf{L}_{h\mathbf{X}},~^{\chi }\mathbf{e}_{h%
\mathbf{X}}\right] \right) ~^{\chi }\mathbf{e}_{h\mathbf{X}}=\left[
\begin{array}{cc}
0 & \left( 0,\overrightarrow{\mathbf{0}}\right) \\
-\left( 0,\overrightarrow{\mathbf{0}}\right) ^{T} & \overrightarrow{\mathbf{v%
}}(\chi )%
\end{array}%
\right] \in \mathfrak{so}(n+1),
\end{equation*}%
where%
\begin{equation*}
\overrightarrow{\mathbf{v}}(\chi )=-\left[
\begin{array}{cc}
0 & \overrightarrow{v}(\chi ) \\
-\overrightarrow{v}^{T}(\chi ) & h\mathbf{0}%
\end{array}%
\in \mathfrak{so}(n)\right] ,
\end{equation*}%
and the derived (applying $ad\left( \left[ ~^{\chi }\mathbf{L}_{h\mathbf{X}%
},~^{\chi }\mathbf{e}_{h\mathbf{X}}\right] \right) $ again )
\begin{eqnarray*}
ad\left( \left[ ~^{\chi }\mathbf{L}_{h\mathbf{X}},~^{\chi }\mathbf{e}_{h%
\mathbf{X}}\right] \right) ^{2}~^{\chi }\mathbf{e}_{h\mathbf{X}} &=&-|%
\overrightarrow{v}(\chi )|^{2}\left[
\begin{array}{cc}
0 & \left( 1,\overrightarrow{\mathbf{0}}\right) \\
-\left( 1,\overrightarrow{\mathbf{0}}\right) ^{T} & \mathbf{0}%
\end{array}%
\right] \\
&=&-|\overrightarrow{v}(\chi )|^{2}~^{\chi }\mathbf{e}_{h\mathbf{X}},
\end{eqnarray*}%
the equation (\ref{31map}) can be represented in alternative form
\begin{eqnarray*}
-\left( h\gamma \right) _{\tau }(\chi ) &=&~^{\chi }\mathbf{D}_{h\mathbf{X}%
}^{2}\left( h\gamma \right) _{h\mathbf{X}}(\chi ) \\
&&-\frac{3}{2}\overrightarrow{R}^{-1}(\chi )~ad\left( ~^{\chi }\mathbf{D}_{h%
\mathbf{X}}\left( h\gamma \right) _{h\mathbf{X}}(\chi )\right) ^{2}~\left(
h\gamma \right) _{h\mathbf{X}}(\chi ),
\end{eqnarray*}%
which is more convenient for analysis of higher order flows on $%
\overrightarrow{v}(\chi )$ subjected to higher--order geometric partial
differential equations. Here we note that the $0$ flow one $\overrightarrow{v%
}(\chi )$ correspond to just convective (linear travelling h--wave but
subjected to certain nonholonomic constraints) map equations (\ref{3trmap}).

Now we consider -1 flows contained in the family of h--hierarchies derived
from the property that $h\overrightarrow{\mathbf{e}}_{\perp }(\chi )$ is
annihilated by the h--operator $h\mathcal{J}(\chi )$ and mapped into $h%
\mathfrak{R}(h\overrightarrow{\mathbf{e}}_{\perp })(\chi )=0.$ This states
that $h\mathcal{J}(h\overrightarrow{\mathbf{e}}_{\perp })(\chi )=~^{\chi }%
\overrightarrow{\varpi }=0.$ Such properties together with (\ref{3auxaaa})
and equations (\ref{3floweq}) imply $~^{\chi }\mathbf{L}_{\tau }=0$ and
hence $h\mathcal{D}_{\tau }\mathbf{e}_{h\mathbf{X}}(\chi )=[~^{\chi }\mathbf{%
L}_{\tau },~^{\chi }\mathbf{e}_{h\mathbf{X}}]=0$ for $h\mathcal{D}_{\tau
}(\chi )=h\mathbf{D}_{\tau }(\chi )+[~^{\chi }\mathbf{L}_{\tau },\cdot ].$
We obtain the equation of motion for the h--component of curve, $h\gamma
(\tau ,\mathbf{l}),$ following the correspondences $~^{\chi }\mathbf{D}_{h%
\mathbf{Y}}\longleftrightarrow h\mathcal{D}_{\tau }(\chi )$ and $h\gamma _{%
\mathbf{l}}(\chi )\longleftrightarrow ~^{\chi }\mathbf{e}_{h\mathbf{X}},$%
\begin{equation*}
~^{\chi }\mathbf{D}_{h\mathbf{Y}}\left( h\gamma (\tau ,\mathbf{l,}\chi
)\right) =0,
\end{equation*}%
which is just the first equation in (\ref{3-1map}).

Finally, we note that the formulas for the v--components, stated by Theorem %
\ref{3mt} can be derived in a similar form by respective substitution in the
the above proof of the h--operators and h--variables into v--ones, for
instance, $h\gamma \rightarrow v\gamma ,$ $h\overrightarrow{\mathbf{e}}%
_{\perp }\rightarrow v\overleftarrow{\mathbf{e}}_{\perp },$ $\overrightarrow{%
v}\rightarrow \overleftarrow{v},\overrightarrow{\varpi }\rightarrow
\overleftarrow{\varpi },\mathbf{D}_{h\mathbf{X}}\rightarrow \mathbf{D}_{v%
\mathbf{X}},$ $\mathbf{D}_{h\mathbf{Y}}\rightarrow \mathbf{D}_{v\mathbf{Y}},%
\mathbf{L\rightarrow C,}\overrightarrow{R}\rightarrow \overleftarrow{S},h%
\mathcal{D\rightarrow }v\mathcal{D},$ $h\mathfrak{R\rightarrow }v\mathfrak{R,%
}h\mathcal{J\rightarrow }v\mathcal{J}$,...where, for simplicity, we omit
parametric dependencies on $\chi .$

\end{document}